\documentclass[leqno,a4paper,11pt]{article}
\usepackage{amsmath,amsthm,amssymb}
\usepackage[utf8]{inputenc}
\usepackage[T1]{fontenc}

\usepackage[a4paper]{geometry} 
\usepackage{hyperref}

\usepackage{enumerate}
\usepackage{bbm}
\usepackage{todonotes}
\usepackage{mathrsfs}

\usepackage{longtable}
\usepackage{ctable}

\usepackage{tikz}

\usepackage{nicefrac}

\usepackage{cleveref}
\crefname{Th}{Theorem}{Theorems}
\crefname{Prop}{Proposition}{Propositions}
\crefname{Lemma}{Lemma}{Lemmas}
\crefname{Cor}{Corollary}{Corollaries}
\crefname{Thx}{Theorem}{Theorems}
\crefname{Propx}{Proposition}{Propositions}
\crefname{Remark}{Remark}{Remarks}
\crefname{Def}{Definition}{Definitions}
\crefname{Example}{Example}{Examples}
\crefname{Question}{Question}{Questions}
\crefname{section}{Section}{Sections}

\usepackage[sort&compress,numbers]{natbib}

\usepackage{centernot}
\providecommand{\noopsort}[1]{} %necessary for the bib file to work

\DeclareMathOperator{\lcm}{lcm}

\newtheorem{Th}{Theorem}[section]
\newtheorem{Prop}[Th]{Proposition}
\newtheorem{Lemma}[Th]{Lemma}
\newtheorem{Cor}[Th]{Corollary}
\newtheorem{Thx}{Theorem}

\newtheorem{Propx}[Thx]{Proposition}

\newtheoremstyle{named}{}{}{\itshape}{}{\bfseries}{.}{.5em}{\thmnote{#3}#1}
\theoremstyle{named}
\newtheorem*{namedtheorem}{}

\theoremstyle{definition}
\newtheorem{Remark}[Th]{Remark}
\newtheorem{Def}{Definition}[section]
\newtheorem{Example}[Th]{Example}
\newtheorem{Question}[Th]{Question}

\newcommand{\divides}{\mid}
\newcommand{\ndivides}{\nmid}

\newcommand{\un}{\underline}
\newcommand{\mob}{\boldsymbol{\mu}}
\newcommand{\cf}{{\cal F}}
\newcommand{\cm}{{\cal M}}
\newcommand{\ov}{\overline}

\newcommand{\Z}{{\mathbb{Z}}}
\newcommand{\N}{{\mathbb{N}}}
\newcommand{\PP}{{\mathbb{P}}}

\newcommand{\vep}{\varepsilon}
\newcommand{\va}{\varphi}
\newcommand{\sA}{\mathscr{A}}
\newcommand{\sB}{\mathscr{B}}
\newcommand{\raz}{\mathbbm{1}}
\newcommand{\bdelta}{\boldsymbol{\delta}}

\newcommand{\sa}{{\underline{s},\underline{a}}}
\newcommand*\samethanks[1][\value{footnote}]{\footnotemark[#1]}

\title{$\mathscr{B}$-free sets and dynamics}
\author{A. Bartnicka\thanks{Research supported by Narodowe Centrum Nauki grant  UMO-2014/15/B/ST1/03736.} \and S. Kasjan\samethanks \and J. Ku\l{}aga-Przymus\samethanks \and M. Lema\'{n}czyk\samethanks}

\begin{document}
\bibliographystyle{siam}

\maketitle

\abstract{Let $\mathscr{B}\subset \N$ and let $\eta\in \{0,1\}^\Z$ be the characteristic function of the set $\mathcal{F}_\mathscr{B}:=\Z\setminus\bigcup_{b\in\mathscr{B}}b\Z$ of $\mathscr{B}$-free numbers. Consider the subshift $(S,X_\eta)$, where $X_\eta$ is the closure of the orbit of $\eta$ under the left shift $S$. In case when $\mathscr{B}=\{p^2 : p\text{ is prime}\}$ the  dynamics of $(S,X_\eta)$ was studied by Sarnak in 2010. This special case and some generalizations, including the case ($\ast$) of $\mathscr{B}$ infinite, pairwise coprime with $\sum_{b\in\mathscr{B}}1/b<\infty$, were discussed by several authors. We continue this line of research for a general $\mathscr{B}$.

The main difference between the general case and the ($\ast$) case is that we may have $X_\eta\subsetneq X_\mathscr{B}:=\{x\in \{0,1\}^\Z : |\text{supp }x\bmod b|\leq b-1 \text{ for each }b\in\mathscr{B}\}$, i.e.\ $X_\eta$ no longer has a characterization in terms of admissible sequences, while in the ($\ast$) case $X_\eta=X_{\mathscr{B}}$. Moreover, $X_\eta$ may not be hereditary (heredity of $X\subset \{0,1\}^\Z$ means that if $x\in X$ and $y\leq x$ coordinatewise then $y\in X$), i.e.\ $X_\eta\neq\widetilde{X}_\eta$, where $\widetilde{X}_\eta$ is the smallest hereditary subshift containing $X_\eta$.

We show that $\eta$ is a quasi-generic point for some natural $S$-invariant measure $\nu_\eta$ on $X_\eta$. We solve the problem of proximality by showing first that $X_\eta$ has a unique minimal subset (to which each point has to be proximal). Moreover, this unique minimal subsystem is a Toeplitz dynamical system which relates the theory of $\mathscr{B}$-free shifts and Toeplitz shifts. We prove that a $\mathscr{B}$-free system is proximal if and only if $\mathscr{B}$ contains an infinite coprime subset.

For other results, including the solution of the problem of invariant measures, a class of sets $\mathscr{B}$, larger than the class given by ($\ast$), which is crucial for us is that of taut sets: $\mathscr{B}$ is taut whenever $\bdelta(\mathcal{F}_\mathscr{B})<\bdelta(\mathcal{F}_{\mathscr{B}\setminus \{b\}})$ for each $b\in\mathscr{B}$ ($\bdelta$ stands for the logarithmic density). We give a characterization of taut sets $\mathscr{B}$ in terms of the support of the corresponding measure $\nu_\eta$. Moreover, for any $\mathscr{B}$ there exists a taut $\mathscr{B}'$ with $\nu_\eta=\nu_{\eta'}$. For taut sets $\mathscr{B},\mathscr{B}'$, we have $\mathscr{B}=\mathscr{B}'$ if and only if $X_\mathscr{B}=X_{\mathscr{B}'}$.

A special role played by $\mathscr{B}$-free systems for taut $\mathscr{B}$ is seen in the following result: For each $\mathscr{B}$ there is a taut $\mathscr{B}'$ such that $(S,\widetilde{X}_{\eta'})$ is a subsystem of $(S,\widetilde{X}_\eta)$ and $\widetilde{X}_{\eta'}$ is a quasi-attractor. In particular, all invariant measures for $(S,\widetilde{X}_\eta)$ are supported by $\widetilde{X}_{\eta'}$.

The system $(S,\widetilde{X}_\eta)$ is shown to be intrinsically ergodic for an arbitrary $\mathscr{B}$. Moreover, we provide a description of all probability invariant measures for $(S,\widetilde{X}_\eta)$.  We prove that the topological entropies of $(S,\widetilde{X}_\eta)$ and $(S,X_\mathscr{B})$ are the same and equal to $\ov{d}(\mathcal{F}_\mathscr{B})$.

We also show that for a subclass of taut $\mathscr{B}$-free systems, namely those for which $\mathscr{B}$ has light tails, i.e.\ $\ov{d}(\sum_{b>K}b\Z)\to 0$,  proximality is the same as heredity.

Finally, we give some applications in number theory on gaps between consecutive $\mathscr{B}$-free numbers. We also apply our results to the set of abundant numbers (positive integers that are smaller than the sum of their proper divisors).}

\tableofcontents

\section{Introduction}
\subsection{Motivation}
\paragraph{Sets of multiples}
For a subset $\sB\subset \N:=\{1,2,\dots\}$, we consider its \emph{set of multiples} $\cm_\sB:=\bigcup_{b\in\sB}b\Z$ and the associated set of \emph{$\sB$-free numbers} $\cf_\sB:=\Z\setminus \cm_\sB$. The interest in sets of multiples was initiated in the 1930s by the study of one particular example, namely, the set of \emph{abundant numbers}, i.e.\ of $n\in\Z$ for which $|n|$ is smaller than the sum of its (positive) proper divisors. In~\cite{bessel1929zahlentheorie}, Bessel-Hagen asked whether the set of abundant numbers has asymptotic density and the positive answer was given independently by Davenport~\cite{Davenport:1933aa}, Chowla~\cite{zbMATH03014412} and Erd\"os~\cite{MR1574879}.  
Nowadays, abundant numbers are still of a certain interest in number theory (see, e.g., the recent works \cite{MR3254753,MR3189967,MR2134854}). 

The works of Davenport, Chowla and Erd\"os triggered various questions on general sets of multiples. In particular, the natural question whether all sets of multiples have asymptotic density was answered negatively by Besicovitch~\cite{MR1512943}. On the other hand, Davenport and Erd\"os~\cite{Davenport1936,MR0043835} showed that $\cm_\sB$ (equivalently, $\cf_\sB$)  always has logarithmic density equal to the lower density. Moreover, in many cases, $\cm_\mathscr{B}$ does have density,, e.g., when
\begin{equation}\label{settingerdosa}
\mathscr{B} \text{ is pairwise coprime and } \sum_{b\in\mathscr{B}}\nicefrac{1}{b}<\infty,
\end{equation}
see, e.g.,~\cite{MR687978}.\footnote{This setting was first studied by Erd\"os~\cite{MR0207673}.} Following~\cite{MR1414678},
all sets $\mathscr{B}\subset\N$ for which $\cm_{\mathscr{B}}$ has density are called \emph{Besicovitch}. 

An important example of a Besicovitch set is
\begin{equation}\label{bf1}
\mathscr{B}=\{p^2  :  p \in\mathscr{P}\},
\end{equation}
where $\mathscr{P}$ denotes the set of primes. Here,~\eqref{settingerdosa} is clearly satisfied. The set $\cf_\sB$ is called the set of \emph{square-free} integers and its density equals $\nicefrac{6}{\pi^2}$, see, e.g., \cite{MR568909}. The characteristic function of $\mathcal{F}_\mathscr{B}$ is the square $\mob^2$ of the M\"obius function $\mob$ extended to $\Z$ in the natural way:  $\mob(-n)=\mob(n)$. (Recall that $\mob(n)=(-1)^k$ when $n$ is a product of $k\geq 1$ distinct primes, $\mob(1)=1$ and $\mob(n)=0$ if $n\in\N$ is not square-free.)

With each set $\cf_{\mathscr{B}}$ of $\mathscr{B}$-free numbers, we associate three natural subshifts\footnote{By a {\em subshift}, we mean a dynamical system $(S,X)$, where $X\subset\{0,1\}^{\Z}$ is closed, $S$-invariant and $S$ stands for the left shift.}
$$
X_{\eta}\subset\widetilde{X}_\eta \subset X_{\mathscr{B}},
$$
with the first and the third defined in the following way:
\begin{itemize}
\item
\emph{$\mathscr{B}$-free subshift} $(S,X_\eta)$, where $X_\eta$ is the closure of the orbit $\mathcal{O}_S(\eta):=\{S^m\eta : m\in\Z\}$ of $\eta=\raz_{\cf_{\mathscr{B}}}\in\{0,1\}^{\Z}$,
\item
\emph{$\mathscr{B}$-admissible subshift} $(S,X_{\mathscr{B}})$, where $X_\mathscr{B}$ is the set of $\mathscr{B}$-{\em admissible} sequences, i.e.\ of $x\in\{0,1\}^{\Z}$ such that, for each $b\in\sB$, the support ${\rm supp}\,x:=\{n\in\Z : x(n)=1\}$ of $x$ taken modulo $b$ is a proper subset of $\Z/b\Z$.\footnote{Admissible blocks and subsets of integers are defined in a similar way. Notice that $X_{\mathscr{B}}$ is closed as the $\mathscr{B}$-admissibility of $x$ is equivalent to the $\mathscr{B}$-admissibility of all finite subsets of ${\rm supp}\,x$. Clearly, $\eta$ is $\mathscr{B}$-admissible.}
\end{itemize}
Notice that the subshift $(S,X_{\mathscr{B}})$ is {\em hereditary}, i.e.\ whenever $x\in X_{\mathscr{B}}$ and $y\leq x$ coordinatewise, then $y\in X_{\mathscr{B}}$. Finally, we consider 
\begin{itemize}
\item
the subshift $(S,\widetilde{X}_\eta)$, where $\widetilde{X}_\eta$ is defined to be the smallest hereditary subshift containing $X_\eta$.
\end{itemize}

\paragraph{Relations with number theory}
Consider two more examples. Let
\begin{equation}\label{bf2}
\mathscr{B}:=\{pq : p,q\in \mathscr{P}\} \text{ and }\mathscr{B}':=\mathscr{P}.
\end{equation}
Then $\cf_{\mathscr{B}}=\mathscr{P}\cup (-\mathscr{P})\cup\{-1,1\}$ and $\cf_{\mathscr{B}'}=\{-1,1\}$. Let $\eta:=\raz_{\cf_{\mathscr{B}}}$, $\eta':=\raz_{\cf_{\mathscr{B}'}}$. Clearly, $X_{\eta'}\subsetneq\widetilde{X}_{\eta'}\subsetneq X_{\mathscr{P}}$.\footnote{$X_\mathscr{P}$ is uncountable, see \cref{nieprz}.} 

Recall the following famous number-theoretical conjectures:
\begin{namedtheorem}[Prime $k$-Tuples Conjecture]
For each $k\geq1$ and each $\mathscr{P}$-admissible subset $\{a_1,\ldots,a_k\}\subset\N\cup\{0\}$, there exist infinitely many $n\in\N$ such that $\{a_1+n,\ldots,a_k+n\}\subset\mathscr{P}$.
\end{namedtheorem}
Note that the set $\{0,2\}$ is $\mathscr{P}$-admissible and the Prime $k$-Tuples Conjecture in this case is the Twin Prime Conjecture. Note also that if, for some $p\in\mathscr{P}$, we have $\{a_i \bmod p : 1\leq i\leq k\}=\Z/p\Z$  and $ \{a_1+n,\ldots,a_k+n\}\subset\mathscr{P}$ then $n=p-a_i$  for some $1\leq i\leq k$, whence the set of $n\in\N$ such that $ \{a_1+n,\ldots,a_k+n\}\subset\mathscr{P}$ is finite.

\begin{Remark}
It is  not hard to see that the Prime $k$-Tuples Conjecture is equivalent to $X_{\mathscr{P}}\subset \widetilde{X}_\eta$. Indeed, for the necessity, we need to show that if a block $B\in\{0,1\}^s$ is $\mathscr{P}$-admissible then there is a block $B'\in \{0,1\}^s$ appearing on $\eta$ such that $B\leq B'$. The existence of such a $B'$ follows directly from the Prime $k$-Tuples Conjecture. Conversely, let $F=\{a_1,\ldots,a_k\}$ be $\mathscr{P}$-admissible. Take $i_0\geq 1$ large enough, so that $2|F|<p_{i_0+1}$. Then the sets $F\cup (F+kp_1\ldots p_{i_0})$, $k\geq1$, are also $\mathscr{P}$-admissible. These sets, for each $k\geq1$, correspond to some blocks $C_k$ appearing in $X_\mathscr{P}$. By assumption, this implies the existence of $C'_k$ on $\eta$ with $C_k\leq C'_k$, $k\geq 1$. It follows that we have $n,m\in\Z$ such that $F+n,F+m\subset\mathscr{P}$ with $|n-m|$ arbitrarily large, and the Prime $k$-Tuples Conjecture follows.
\end{Remark}

\begin{namedtheorem}[Dickson's Conjecture, \cite{dickson}]
Let $a_i\in \Z$, $b_i\in\N$, $1\leq i\leq k$. If for each $p\in\mathscr{P}$ there exists $n\in\N$ such that $p\ndivides \prod_{1\leq i\leq k}(b_in+a_i)$ then there are infinitely many $n\in\N$ such that $b_in+a_i\in\mathscr{P}$ for $1\leq i\leq k$.
\end{namedtheorem}
Note that if $b_i=1$ for $1\leq i\leq k$, the condition that for each $p\in\mathscr{P}$ there exists $n\in\N$ such that $p\ndivides \prod_{1\leq i\leq k}(b_in+a_i)$ is equivalent to the $\mathscr{P}$-admissibility of $\{a_1,\dots,a_k\}$. 

\begin{Remark}
The following consquence of Dickson's conjecture (more specifically, of its special case when $b_i=1$ for $1\leq i\leq k$) was pointed to us by Professor A.~Schinzel, see $C_{13}$ in~\cite{MR0130203}:
$$
\parbox{0.8\textwidth}{
If $a_1,\dots,a_k \in [-n,n]\cap \Z$ and $\{a_1,\dots,a_k\}$ is $\mathscr{P}$-admissible then, for infinitely many $x\in \N$, we have $[x-n,x+n]\cap\mathscr{P}=\{x+a_i : i=1,\dots,k\}$.
}
$$
This can be rephrased as $X_\mathscr{P}\subset X_\eta$.
\end{Remark}

\paragraph{Dynamical approach}
The above suggests that the sets of multiples and the associated subshifts are difficult to study in full generality. Thus, it seems natural to put first some restrictions on $\mathscr{B}$ and then try to relax them to see which from the previous results ``survive''. Sarnak in his seminal paper~\cite{sarnak-lectures} suggested to study dynamical properties of the \emph{square-free subshift} $(S,X_{\mob^2})$.  He formulated a certain program, in particular, announcing the following results:
\begin{enumerate}[(i)]
\item $\mob^2$ is \emph{generic} for an ergodic $S$-invariant measure $\nu_{\mob^2}$ on $\{0,1\}^\Z$ such that the corresponding measure-theoretical dynamical system $(S,X_{\mob^2},\nu_{\mob^2})$ has zero Kolmogorov entropy,\label{sa1}
\item the topological entropy of $(S,X_{\mob^2})$ is equal to $\nicefrac{6}{\pi^2}$,\label{sa2}
\item  $X_{\mob^2}=X_\mathscr{B}$, where $\mathscr{B}=\{ p^2 : p\in \mathscr{P}\}$,\label{sa3}
\item $(S,X_{\mob^2})$ is \emph{proximal},\label{sa4}
\item $(S,X_{\mob^2})$ has a non-trivial \emph{topological joining} with a rotation on a compact Abelian group\label{sa5}
\end{enumerate}
(we will explain the notions appearing in \eqref{sa1}-\eqref{sa5} later). Today, complete proofs of these facts are available; Sarnak's program has also been studied for some natural generalizations of $(S,X_{\mob^2})$, see~\cite{MR3055764,Peckner:2012fk,Huck:2014aa,Ab-Le-Ru,MR3070541,MR3296562,Baake:2015aa,Ba-Ku}.\footnote{Cf.\ also~\cite{MR0485769,MR1778906} for the harmonic analysis viewpoint.} In particular, in~\cite{Ab-Le-Ru}, Abdalaoui, Lema\'nczyk and de la Rue cover the counterparts of \eqref{sa1}-\eqref{sa3} from Sarnak's list for $\mathscr{B}\subset \N$ satisfying~\eqref{settingerdosa}. In this case, by \eqref{sa3}, we have $X_\eta=\widetilde{X}_\eta=X_\mathscr{B}$.

As we have already mentioned, we intend to relax the assumptions~\eqref{settingerdosa} on $\mathscr{B}$ and tackle similar problems to \eqref{sa1}-\eqref{sa5}.\footnote{This problem was posed during the conference \emph{Ergodic Theory and Dynamical Systems} in Toru\'n, Poland 2014 by M.\ Boshernitzan.} It is all the more important, since $X_{\mathscr{B}'}\subset X_{\mathscr{B}}$ whenever $\mathscr{B}\subset \mathscr{B}'\subset \N$. In other words, any $(S,X_{\mathscr{B}})$ has subsystems of the form $(S,X_{\mathscr{B}'})$ for certain sets $\mathscr{B}'\subset \N$ whose elements are no longer pairwise coprime. (Another way to obtain a natural subsystem of $(S,X_\mathscr{B})$ is to choose $b' \divides b$ for each $b\in\mathscr{B}$ and then note that $X_{\mathscr{B}'}\subset X_\mathscr{B}$, where $\mathscr{B}'=\{b' :b\in\mathscr{B}\}$.) In particular, this applies to the square-free case. As a matter of fact, the square-free subshift constains $X_\mathscr{B}$ whenever $\{p^2 : p\in\mathscr{P}\}\subset \sB \subset\{pq : p,q\in\mathscr{P}\}$, cf.\ \eqref{bf1} and~\eqref{bf2}. 

Recall also that in~\cite{MR3356811} a description of all invariant measures for $(S,X_{\mathscr{B}})$ was given in case~\eqref{settingerdosa}. Moreover, under the same assumptions, $(S,X_{\mathscr{B}})$ was proved to be \emph{intrinsically ergodic} (this means that the system has only one invariant measure $\nu$ such that the Kolmogorov entropy of $(S,X_{\mathscr{B}},\nu)$ is equal to the topological entropy of $(S,X_{\mathscr{B}})$).\footnote{The intrinsic ergodicity of $(S,X_{\mob^2})$ was proved in~\cite{Peckner:2012fk}.}

The present paper seems to be the first attempt to deal with Sarnak's list \eqref{sa1}-\eqref{sa5} and the problem of invariant measures in the general case when $\mathscr{B}\subset \N$, i.e.\ when we drop the assumption~\eqref{settingerdosa}. Sometimes, we put certain restrictions on $\mathscr{B}$. In particular, we deal with $\mathscr{B}$ that:
\begin{itemize}
\item
are \emph{thin}, i.e.\ $\sum_{b\in\mathscr{B}}\nicefrac{1}{b}<\infty$,
\item
have  \emph{light tails}, i.e.\ $\ov{d}(\sum_{b>K}b\Z)<\vep$ for $K$ large enough.
\end{itemize}
Each thin $\mathscr{B}$ has light tails and if $\mathscr{B}$ is pairwise coprime, these two notions coincide. Moreover, light tail sets are Besicovitch. A more subtle notion, which turns out to be crucial in our studies, is that of \emph{tautness}~\cite{MR1414678}: 
\begin{itemize}
\item $\mathscr{B}$ is \emph{taut} when $\bdelta(\mathcal{M}_{\mathscr{B}\setminus \{b\}})<\bdelta(\mathcal{M}_\mathscr{B})$ for each $b\in\mathscr{B}$.\footnote{Symbol $\bdelta$ stands for the logarithmic density.}
\end{itemize}
Any \emph{primitive} set $\sB$ (i.e.\ such that, for $b,b'\in\sB$, we have $b\ndivides b'$) with light tails is taut.

The main difference between the general situation and the setting~\eqref{settingerdosa} is that $X_\eta$ has no longer a characterization in terms of admissible sequences, i.e.\ it may happen that 
the $\mathscr{B}$-admissible subshift $(S,X_\mathscr{B})$ is strictly larger than the $\mathscr{B}$-free subshift $(S,X_\eta)$. What is more, while $X_\mathscr{B}$ is always hereditary, $X_\eta$ need not be so,  and, as we have already seen by inspecting the case $\mathscr{B}=\mathscr{P}$, we may even have $X_\eta\subsetneq \widetilde{X}_\eta\subsetneq X_\mathscr{B}$. On the other hand, there are many similarities or analogies between~\eqref{settingerdosa} and the general case.

\subsection{Main results}\label{smain}
Our main results can be divided into three groups:
\begin{enumerate}[(I)]
\item\label{GI}
structural results,
\item\label{GII}
results on invariant measures and entropy,
\item\label{GIII}
number theoretical results.
\end{enumerate}
The results from groups \eqref{GI} and \eqref{GII} are closely related to one another, whereas the results from group \eqref{GIII} are mostly consequences of the results from \eqref{GI} and \eqref{GII}.

\subsubsection{Structural results}\label{sestru} This group of results contains both topological and measure-theoretical results. Namely, we have:
\begin{Thx}\label{TTA}
For any $\mathscr{B}\subset \N$, the subshift $(S,X_\eta)$ has a unique minimal subset. Moreover, this subset is the orbit closure of a Toeplitz sequence.
\end{Thx}
\begin{Remark}
\cref{TTA} is an extension of \eqref{sa4} from Sarnak's program.
\end{Remark}\noindent
As a consequence of \cref{TTA}, we obtain the following result:
\begin{Cor}\label{TTAwn}
For any $\mathscr{B}\subset \N$, each point $x\in X_\eta$ is proximal to a point in the orbit closure of a Toeplitz sequence.
\end{Cor}\noindent
Moreover, as an immediate consequence of \cref{TTA} and \cref{TTAwn}, we have:
\begin{Cor}\label{TTAwn1}
Let $\mathscr{B}\subset \N$. Then $(S,X_\eta)$ is minimal if and only if $(S,X_\eta)$ is a Toeplitz system.
\end{Cor}\noindent
We also give a simple characterization of those $\mathscr{B}\subset \N$, for which the unique minimal subset of $(S,X_\eta)$ is a singleton:
\begin{Thx}\label{TTB}
Let $\sB\subset \N$. The following conditions are equivalent:
\begin{itemize}
\item
the unique minimal subset of $(S,X_\eta)$ is a singleton,
\item
$\{(\dots,0,0,0,\dots)\}$ is the unique minimal subset of $(S,X_\eta)$,
\item
$(S,X_\eta)$ is proximal,
\item
$\mathscr{B}$ contains an infinite pairwise coprime subset.
\end{itemize}
\end{Thx}
It turns out that measure-theoretic properties of the subshift $(S,\widetilde{X}_\eta)$ strongly depend on the notion of tautness. We have:
\begin{Thx}\label{TTC}
For any $\mathscr{B}\subset \N$, there exists a unique taut set $\mathscr{B}'\subset \N$ such that $\cf_{\sB'}\subset \cf_{\sB}$, $\widetilde{X}_{\eta'}\subset \widetilde{X}_\eta$ and $\mathcal{P}(S,\widetilde{X}_\eta)=\mathcal{P}(S,\widetilde{X}_{\eta'})$.\footnote{Given a topological dynamical system $(T,X)$, by $\mathcal{P}(T,X)$ we denote the set of all probability Borel $T$-invariant measures on $X$.}
\end{Thx}\noindent
Equivalently, \cref{TTC} can be rewritten as follows:
\begin{Cor}\label{TTC1}
For any $\mathscr{B}\subset \N$, there exists a unique taut set $\mathscr{B}'\subset \N$ such that $\cf_{\sB'}\subset \cf_\sB$ and any point $x\in \widetilde{X}_\eta$ is attracted to $\widetilde{X}_{\eta'}$ along a sequence of integers of density $1$:
$$
\lim_{n\to\infty, n\not\in E_x} d(T^nx,\widetilde{X}_{\eta'})=0,\text{ where }d(E_x)=0.
$$
\end{Cor}\noindent
A key ingredient in the proof of \cref{TTC} is the description of all invariant measures on $\widetilde{X}_\eta$, see \cref{OOG} below.\footnote{It follows from \cref{OOG} that in order to prove \cref{TTC}, it suffices to construct a taut set $\mathscr{B}'$ such that $\nu_{\eta'}=\nu_\eta$.}

If~\eqref{settingerdosa} is satisfied, then, as shown in~\cite{Ab-Le-Ru}, we have $X_\eta=\widetilde{X}_\eta=X_\mathscr{B}$, cf.\ \eqref{sa3} in Sarnak's program. In general, this need not be the case. However, we have:
\begin{Thx}\label{TTD}
Let $\mathscr{B}\subset \N$. If $\mathscr{B}$ has light tails and contains an infinite, pairwise coprime subset then $X_\eta=\widetilde{X}_\eta$.
\end{Thx}\noindent
In other words, for primitive $\mathscr{B}$ with light tails, the proximality of $(S,X_\eta)$ is equivalent to the heredity of $X_\eta$. Since every $\mathscr{B}$ that is primitive and has light tails, is taut, a natural question arises whether the assertion of \cref{TTD} remains true for all taut $\mathscr{B}\subset \N$. We leave this question open, conjecturing that the answer is positive.

\subsubsection{Results on invariant measures and entropy}\label{se122}
\begin{Propx}\label{OOE}
For any $\sB\subset \N$, $\eta=\raz_{\cf_\sB}$ is a quasi-generic point for a natural ergodic $S$-invariant measure $\nu_\eta$ on $\{0,1\}^\Z$. In particular, $\nu_\eta(X_\eta)=1$. Moreover, if $\sB$ is Besicovitch then $\eta$ is generic for $\nu_\eta$. 
\end{Propx}
\begin{Remark}
\cref{OOE} means that, for some $(N_k)$, we have the weak convergence $\frac{1}{N_k}\sum_{n\leq N_k}\delta_{S^n\eta}\to \nu_\eta$.
Recall that in case~\eqref{settingerdosa}, this convergence holds along $(N_k)$ with $N_k=k$, see~\cite{Ab-Le-Ru} (i.e.\ $\eta$ is {\em generic} in this case). Recall also that in \eqref{settingerdosa}, $\sB$ is Besicovitch.

We call $\nu_\eta$ the Mirsky measure (in the square-free case the frequencies of blocks on $\eta$ were first studied by Mirsky~\cite{MR0021566,MR0028334}).
\end{Remark}

\begin{Thx}\label{IZOIZO}
Suppose that $\sB\subset \N$ is taut. Then $(S,X_\eta,\nu_\eta)$ is isomorphic to $(T,G,\PP)$, where $G$ is the closure of ${\{(n,n,\dots)\in \prod_{k\geq 1}\Z/b_k\Z : n\in \Z\}}$ in $\prod_{k\geq 1}\Z/b_k\Z$ and $Tg=g+(1,1,\dots)$. In particular, $(S,X_\eta,\nu_\eta)$ has zero entropy.
\end{Thx}

\begin{Remark}
\cref{OOE}, together with \cref{IZOIZO}, extends \eqref{sa1} from Sarnak's program.
\end{Remark}

\begin{Thx}\label{OOF}
If $\mathscr{B}\subset \N$ has light tails then $X_\eta$ is the topological support of~$\nu_\eta$.
\end{Thx}

\begin{Thx}\label{tautychar}
Let $Y:=\{x\in \{0,1\}^\Z : |\text{supp }y \bmod b|=b-1\text{ for each }b\in\sB\}$. For $\mathscr{B}\subset \N$ infinite (and primitive), the following conditions are equivalent:
\begin{enumerate}[(a)]
\item\label{tautycharA}
$\mathscr{B}$ is taut,
\item\label{tautycharB}
$\mathcal{P}(S,Y\cap \widetilde{X}_\eta)\neq \emptyset$,
\item\label{tautycharC}
$\nu_\eta(Y\cap X_\eta)=1$.
\end{enumerate}
\end{Thx}

\begin{Thx}\label{OOG}
For any $\mathscr{B}\subset \N$ and any $\nu\in\mathcal{P}(S,\widetilde{X}_\eta)$, there exists $\rho\in\mathcal{P}(S\times S,X_\eta\times \{0,1\}^\Z)$ whose projection onto the first coordinate equals $\nu_\eta$ and such that $M_\ast(\rho)=\nu$, where $M\colon X_\eta\times \{0,1\}^\Z\to \widetilde{X}_\eta$ stands for the coordinatewise multiplication.
\end{Thx}
\begin{Thx}\label{OOH}
For any $\mathscr{B}\subset \N$, the subshift $(S,\widetilde{X}_\eta)$ is intrinsically ergodic.
\end{Thx}
An important tool here, which can be also of independent interest, is the following result:
\begin{Propx}\label{OOI}
For any $\mathscr{B}\subset \N$, we have $h_{top}(S,\widetilde{X}_\eta)=h_{top}(S,X_\mathscr{B})=\bdelta(\cf_{\mathscr{B}})$.
\end{Propx}
\begin{Remark}
\cref{OOI} is an extension of \eqref{sa2} from Sarnak's program (recall that the density of square-free numbers equals $\nicefrac{6}{\pi^2}$, see, e.g., \cite{MR568909}).
\end{Remark}
The last entropy result we would like to highlight here is the following immediate consequence of \cref{TTC} and the variational principle:
\begin{Cor}\label{OOIwn}
For any $\mathscr{B}\subset \N$, there exists a taut set $\mathscr{B}'\subset \N$ such that $\cf_{\sB'}\subset \cf_{\sB}$ and $h_{top}(S,\widetilde{X}_\eta)=h_{top}(S,\widetilde{X}_{\eta'})$.
\end{Cor}

\subsubsection{Number theoretical results}\label{se123}
\paragraph{General consequences}
Our first result in this section shows, in particular, that a taut set $\sB$ is determined by the family of $\sB$-admissible subsets.
\begin{Thx}\label{OOJ}
Suppose that $\mathscr{B},\mathscr{B}'\subset \N$ are taut. Then the following conditions are equivalent:
\begin{enumerate}[(a)]
\item $\mathscr{B}=\mathscr{B}'$,
\item $\cm_\sB=\cm_{\sB'}$,
\item $X_\mathscr{B} = X_{\mathscr{B}'}$,
\item $\widetilde{X}_\eta = \widetilde{X}_{\eta'}$,
\item $X_\eta=X_{\eta'}$,
\item
$\nu_\eta=\nu_{\eta'}$,
\item
$\mathcal{P}(S,\widetilde{X}_\eta)= \mathcal{P}(S,\widetilde{X}_{\eta'})$.
\end{enumerate}
\begin{Remark}
\cref{OOJ} extends an analogous result from~\cite{MR3356811}, where it was shown that $X_\sB=X_{\sB'}$ is equivalent to $\sB=\sB'$ for $\mathscr{B},\mathscr{B}'\subset\N$ satisfying \eqref{settingerdosa}.
\end{Remark}

\end{Thx}\noindent
As an immediate consequence of \cref{OOE} and \cref{OOF}, we obtain:
\begin{Cor}\label{OOK}
If $\mathscr{B}\subset \N$ has light tails, $F,M\subset \N$ are finite sets such that $F\subset \cf_\sB$, $M\subset \cm_\sB$ then the density of the set of $n\in\N$ such that $F+n\subset \cf_\sB$, $M+n\subset \cm_\sB$ is positive.
\end{Cor}\noindent
\begin{Propx}\label{OOL}
Suppose that $\sB\subset \N$ has light tails and contains an infinite coprime subset $\sB'$. Denote by $(n_j)$ the sequence of consecutive $\mathscr{B}$-free numbers. Then
$$
\limsup_{j\to \infty} \inf_{0\leq k\leq K}(n_{j+k+1}-n_{j+k})=\infty \text{ for any }K\geq 1.
$$
\end{Propx}

\paragraph{Consequences for abundant numbers}
\begin{Cor}\label{VFVF}
Suppose that $A,D\subset \N$ are finite sets, consisting of abundant and non-abundant numbers, respectively. Then the density of $n\in\N$ such that $A+n$ and $D+n$ consist of abundant and deficient numbers, respectively, is of positive density.
\end{Cor}\noindent
\begin{Cor}\label{VFVF1}
The set of $n\in\N$ such that the numbers $n+1,n+2,\dots,n+5$ are deficient has positive density.
\end{Cor}
\begin{Cor}\label{VFVF2}
Denote by $(n_j)$ the sequence of consecutive deficient numbers. Then, for any $K\geq 1$,
$$
\limsup_{j\to \infty} \inf_{0\leq k\leq K}(n_{j+k+1}-n_{j+k})=\infty.
$$
\end{Cor}
\begin{Cor}\label{VFVF3}
Let $\eta:=1-\raz_{\mathbf{A}}\in \{0,1\}^\Z$, where $\mathbf{A}$ is the set of abundant integers. Then $X_\eta=\widetilde{X}_\eta$, in particular $(S,X_\eta)$ is proximal. Moreover, $(S,X_\eta)$ is intrinsically ergodic and we have $h_{top}(S,X_\eta)=1-d(\mathbf{A})$.
\end{Cor}
It remains an open question whether $X_\eta=X_{\sB_\mathbf{A}}$.

\subsection{`Map' of the paper}
In this section we include a table that can be used to locate within the paper the proofs of the main results listed in \cref{smain}.

\setlength{\LTpre}{0pt}
\setlength{\LTpost}{-10pt}

\begin{center}
\begin{longtable}{|c|c|c|}
\hline
\textbf{Result} & \textbf{Proof} & \textbf{Main tools} \\
\specialrule{.2em}{.1em}{.1em} 
\endfirsthead
\multicolumn{3}{l}%
{\textit{Continued from previous page}} \\
\hline
\textbf{Result} & \textbf{Proof} & \textbf{Important tools}\\
\specialrule{.2em}{.1em}{.1em} 
\endhead
\hline \multicolumn{3}{r}{\textit{Continued on next page}} \\
\endfoot
\hline
\endlastfoot
	\cref{TTA} & \cref{proofTTA} & \cref{dlapod}, \cref{topl}\\
	\hline
	\cref{TTAwn} & \cref{proofTTA} & \cref{TTA}, \cref{ausel1}\\
	\hline
	\cref{TTAwn1} & \cref{sestru} & \cref{TTA}, \cref{TTAwn} \\
	\hline
	\cref{TTB} & \cref{TTBsekcja} & Chinese Remainder Theorem\\
	\hline
	\cref{TTC} & \begin{tabular}{@{}c@{}} \cref{tautyuberall}, \\\cref{seatraktor}\end{tabular}  &\begin{tabular}{@{}c@{}}\cref{beh2} and \cref{beh3},  \\	\cref{tautywszystko}, \cref{OOG} and \cref{OOJ}\end{tabular} \\
	\hline
	\cref{TTC1} &\cref{seatraktor} &\cref{TTC}, \cref{lemmaatractor}\\
	\specialrule{.2em}{.1em}{.1em} 
	\cref{TTD} & \cref{s6} & \cref{new1}\\
	\hline
	\cref{OOE} &  \cref{s3} & \cref{da-er} \\
	\hline
	\cref{IZOIZO} & \cref{sewymierne} &  \cref{notry}, \cref{notrynotry}\\
	\hline
	\cref{OOF} &\cref{s6}  & \cref{new1}, \cref{new1a}\\
	\hline
	\cref{tautychar} &\cref{s8} & \cref{TTC}, \cref{OOE}\\
	\hline
	\cref{OOG} & \cref{s10} & \cref{twY}, \cref{miaryall}\\
	\hline
	\cref{OOH} &\begin{tabular}{@{}c@{}} \cref{s9},\\\cref{revi2}\end{tabular} &  \begin{tabular}{@{}c@{}} \cref{jakzbenjim}, \cref{TTC}\\ and the variational principle \end{tabular}\\
	\hline
	\cref{OOI} & \cref{seOOI} & \cref{lerog}\\
	\hline
	\cref{OOIwn} & \cref{se122} & \cref{TTC} and the variational principle 	\\
	\specialrule{.2em}{.1em}{.1em} 
	\cref{OOJ} &\begin{tabular}{@{}c@{}} \cref{Jfirst},\\ \cref{Jfirst1} \end{tabular} &\begin{tabular}{@{}c@{}} \cref{Dirichlet}, \cref{druggi2} \\ \cref{OOG}, \cref{OOI}\end{tabular} \\
	\hline
	\cref{OOK} & \cref{se123} &\cref{OOE}, \cref{OOF} \\
	\hline
	\cref{OOL} & \cref{conse} & \cref{TTD}, \cref{OOE}, \cref{OOF} \\
	\hline
	\cref{VFVF} &\cref{seabu} &\cref{lm121} and \cref{OOK} \\
	\hline
	\cref{VFVF1} &\cref{seabu} & \cref{VFVF} \\
	\hline
	\cref{VFVF2} & \cref{seabu}& \cref{OOL}, \cref{lm121}, \cref{infi}\\
	\hline
	\cref{VFVF3} &\cref{seabu} &\begin{tabular}{@{}c@{}} \cref{lm121}, \cref{infi}, \cref{TTB},\\ \cref{TTD}, \cref{OOH},  \cref{OOI}\end{tabular}
\end{longtable}
\end{center}

\section{Preliminaries}

\subsection{Topological dynamics: basic notions}
\begin{Def}
A \emph{topological dynamical system} is a pair $(T,X)$, where $X$ is a compact space endowed with a metric $d$ and $T$ is a homeomorphism of $X$. We denote by $\mathcal{O}_T(x)$ the orbit of $x\in X$ under $T$, i.e.\ $\mathcal{O}_T(x)=\{T^nx :n\in\Z\}$.
\end{Def}

\begin{Def}
We say that $(T,X)$ is \emph{transitive} if it has a dense orbit. A point $x\in X$ is called \emph{transitive} if $\mathcal{O}_T(x)$ is dense in $X$.
\end{Def}
\begin{Remark}\label{tranzyrem}
Recall that $(T,X)$ is transitive if and only if, for any open sets $U,V\subset X$, there exists $n\in\Z$ such that $T^{-n}U\cap V\neq\emptyset$.
\end{Remark}

\begin{Def}
A point $x\in X$ is called \emph{recurrent} if, for any open set $U\ni x$, there exists $n\neq 0$ such that $T^nx\in U$.
\end{Def}

\begin{Def}
A dynamical system $(T,X)$ is called \emph{topologically weakly mixing} if $(T\times T,X\times X)$ is transitive.
\end{Def}

\begin{Def}
A \emph{minimal set} $M\subset X$ is a non-empty, closed, $T$-invariant set that is minimal with respect to these properties. Equivalently, $M\subset X$ is minimal if for any $x\in M$, we have $\ov{\mathcal{O}_T(x)}=M$. If $M=X$ then $T$ is called \emph{minimal}. A point $x\in X$ is called \emph{minimal} if $(T,\ov{\mathcal{O}_T(x)})$ is minimal.
\end{Def}

\begin{Def}
Let $(T,X)$ be a topological dynamical system. A subset $C\subset X$ is called \emph{wandering} whenever the sets $T^n C$, $n\in\Z$, are pairwise disjoint.
\end{Def}

Given a topological dynamical system $(T,X)$, by $\mathcal{P}(T,X)$ we will denote the set of all Borel probability $T$-invariant measures on $X$ and by $\mathcal{P}^e(T,X)$ the subset of $\mathcal{P}(T,X)$ of ergodic measures (cf.\ \cref{ergody}).
\begin{Def}
If $\mathcal{P}(T,X)$ is a singleton, we say that $(T,X)$ is \emph{uniquely ergodic}.
\end{Def}

\begin{Def}
We say that $x\in X$ is \emph{generic} for $\mu\in\mathcal{P}(T,X)$ if the ergodic theorem holds for $T$ at $x$ for any continuous function $f\in C(X)$: $\frac{1}{N}\sum_{n\leq N}f(T^nx)\to \int f\ d\mu$. 
\end{Def}

\begin{Remark}
In any uniquely ergodic systems all points are generic for the unique invariant measure.
\end{Remark}
\begin{Example}\label{monoerg}
Consider $(T,G)$, where $G$ is a compact Abelian group and $Tg=g+g_0$ for some $g_0\in G$. If $(T,G)$ is {minimal} then it is uniquely ergodic and Haar measure $\PP$ is the unique member of $\mathcal{P}(T,G)$. In particular, all points $g\in G$ are generic for~$\PP$.
\end{Example}

\begin{Def}[see, e.g., \cite{MR1958753}]
A topological dynamical system $(T,X)$ is called \emph{equicontinuous} if the family of maps $\{T^n : n\in\Z\}$ is equicontinuous. Every topological dynamical system has the largest equicontinuous factor, which is called the \emph{maximal equicontinuous factor}.
\end{Def}
\begin{Remark}
All compact Abelian group rotations $(T,G)$ are equicontinuous.
\end{Remark}

\begin{Example}
Let $A$ be a finite set and let $S \colon A^\Z\to A^\Z$ be the left shift, i.e., $S((x_n)_{n\in\Z})=(y_n)_{n\in\Z}$, where $y_n=x_{n+1}$ for each $n\in\Z$. Let $X\subset A^\Z$ be closed and $S$-invariant. We then say that $(S,X)$ is a \emph{subshift}. 
\end{Example}

\begin{Def}\label{top}
We say that $x\in \{0,1\}^\Z$ is a \emph{Toeplitz sequence} whenever for any $n\in\Z$ there exists $d_n\in\N$ such that $x(n+k\cdot d_n)=x(n)$ for any $k\in\Z$. A subshift $(S,Z)$, $Z\subset \{0,1\}^\Z$ is said to be \emph{Toeplitz} if $Z=\ov{\mathcal{O}_S(y)}$ for some Toeplitz sequence $y\in \{0,1\}^\Z$.
\end{Def}

\begin{Remark}
Usually, one requires from a Toeplitz sequence not to be periodic. For convenience, periodic sequences are included in the \cref{top}. We refer the reader, e.g., to~\cite{MR2180227} for more information on Toeplitz sequences.
\end{Remark}

\subsection{Measure-theoretic dynamics: basic notions}
\begin{Def}
A \emph{measure-theoretic dynamical system} is a 4-tuple $(T,X,\mathcal{B},\mu)$, where $(X,\mathcal{B},\mu)$ is a standard probability Borel space and $T$ is an automorphism of~$(X,\mathcal{B},\mu)$. The set of all automorphisms of $(X,\mathcal{B},\mu)$ will be denoted by $Aut(X,\mathcal{B},\mu)$.
\end{Def}

\begin{Def}\label{ergody}
We say that $T\in Aut(X,\mathcal{B},\mu)$ is \emph{ergodic} if, for $A\in\mathcal{B}$, $A=T^{-1}A$ implies $\mu(A)\in \{0,1\}$.
\end{Def}

\begin{Def}
For $T\in Aut(X,\mathcal{B},\mu)$, we define the associated \emph{Koopman operator} $U_T \colon L^2(X,\mathcal{B},\mu)\to L^2(X,\mathcal{B},\mu)$ by setting $U_Tf=f\circ T$.
\end{Def}

\begin{Def}
We say that $\lambda \in \mathbb{S}^1$ is in the \emph{discrete spectrum} of $T\in Aut(X,\mathcal{B},\mu)$ if it is an eigenvalue of $U_T$, i.e., for some $0\neq f\in L^2(X,\mathcal{B},\mu)$, we have
$U_Tf=\lambda f$.
\end{Def}
\begin{Def}
We say that $T\in Aut(X,\mathcal{B},\mu)$ has \emph{purely discrete spectrum} if the eigenfunctions of $U_T$ are linearly dense in $L^2(X,\mathcal{B},\mu)$.
\end{Def}

\begin{Def}[\cite{MR0230877}]
We say that $T\in Aut(X,\mathcal{B},\mu)$ is {\em coalescent} if each endomorphism of $(X,\mathcal{B},\mu)$ commuting with $T$ is invertible.
\end{Def}
\begin{Remark}
All ergodic automorphisms with purely discrete spectrum are coalescent. 
\end{Remark}

\begin{Def}[\cite{MR0213508}]
Let $T\in Aut(X,\mathcal{B},\mu)$, $S\in Aut(Y,\mathcal{C},\nu)$ and let $\rho$ be a $T\times S$-invariant measure on $X\times Y$. We say that $\rho$ is a \emph{joining} of $T$ and $S$ if $\rho|_X=\mu$ and $\rho|_Y=\nu$. In a similar way, joinings of more automorphisms (finitely many and countably many) are defined.
\end{Def}

\begin{Def}
Let $T\in Aut(X,\mathcal{B},\mu)$ and let $C\in \mathcal{B}$ be such that $\mu(C)>0$. Then the function $n_C\colon X\to \N\cup\{\infty\}$ given by
$$
n_C(x)=\min\{n\geq 1 : T^nx \in C\}
$$
is well-defined and finite for $\mu$-a.e.\ $x\in C$. The map $T_C\colon C\to C$ given by $T_Cx=T^{n_C}x$ is called the \emph{induced transformation}. $T_C\in Aut(C,\mathcal{B}_C,\mu_C)$, where $\mathcal{B}_C=\mathcal{B}|_C$ and $\mu_C(A)=\frac{\mu(A)}{\mu(C)}$ for any $A\in \mathcal{B}_C$.
\end{Def}

\subsection{Entropy: basic notions}
There are two basic notions of entropy: \emph{topological entropy} and \emph{measure-theoretic entropy}. We skip the definitions and refer the reader, e.g., to \cite{MR2809170} instead. The topological entropy of $(T,X)$ will be denoted by $h_{top}(T,X)$. The mesure-theoretic entropy of $(T,X,\mathcal{B},\mu)$ will be denoted by $h(T,X,\mu)$. 
\begin{Remark}[Variational principle]
For any topological dynamical system $(T,X)$, we have 
$
h_{top}(T,X)=\sup_{\mu\in\mathcal{P}(T,X)}h(T,X,\mu).
$ 
\end{Remark}
\begin{Def}
If $\mu\in\mathcal{P}(T,X)$ is such that $h(T,X,\mu)=h_{top}(T,X)$, we say that $\mu$ is a \emph{measure of maximal entropy}.
\end{Def}
\begin{Remark}
A measure of maximal entropy may not exist. Subshifts always have at least one measure of maximal entropy.
\end{Remark}
\begin{Def}[\cite{MR0267076}]
$(T,X)$ is said to be \emph{intrinsically ergodic} if it has exactly one measure of maximal entropy.
\end{Def}

\subsection{Topological dynamics: more on minimal subsets}
Let $(T,X)$ be a topological dynamical system.
\begin{Def}
$S\subset \Z$ is called \emph{syndetic} if there exists a finite set $K$ such that $K+S=\Z$.
\end{Def}
\begin{Remark}\label{k2}
There is a well-known characterization of minimality of an orbit closure. Let $x\in X$. Then $(T,\ov{\mathcal{O}_T(x)})$ is minimal if and only if, for any open set $U\ni x$, the set $\{n\in\Z : T^nx\in U\}$ is syndetic. In particular, if $x$ is transitive (i.e.\ its orbit under $T$ is dense in $X$) then $(T,X)$ is minimal if and only if, for any open set $U\subset X$, the set $\{n\in\Z : T^nx\in U\}$ is syndetic.
\end{Remark}

We will be particularly interested in the situation when $(T,X)$ has a unique minimal subset. We first recall well-known results related to the proximal case.
\subsubsection{Proximal case}
\begin{Def}
A pair $(x,y)\in X\times X$ is called \emph{proximal} if $\liminf_{n\to\infty}d(T^nx,T^ny)=0$. We denote the sets of all proximal pairs $(x,y)$ by $\text{Prox}(T)$. $T$ is called \emph{proximal} if $\text{Prox}(T)=X\times X$.
\end{Def}
\begin{Remark}\label{pr:4-}
Note that if $(x,Tx)\in\text{Prox}(T)$ then clearly $T$ has a fixed point. Moreover, $(T,X)$ is proximal if and only if it has a fixed point that is the unique minimal subset of $X$.
\end{Remark}
Recall also the following result:
\begin{Prop}[Auslander - Ellis, see, e.g., \cite{MR2437846}]\label{ausel}
Let $(T,X)$ be a topological dynamical system. Then for any $x\in X$ there exists a minimal point $y\in X$ such that $x$ and $y$ are proximal.
\end{Prop}

\begin{Def}
A pair $(x,y)\in X\times X$ is called  \emph{syndetically proximal} if $\{n\in\Z : d(T^nx,T^ny)<\vep\}$ is syndetic for any $\vep>0$. We denote the set of all syndetically proximal pairs $(x,y)$ by $\text{SyProx}(T)$. $T$ is called \emph{syndetically proximal} if $\text{SyProx}(T)=X\times X$.
\end{Def}
\begin{Remark}\label{subsystem}
Clearly, a subsystem of a (syndetically) proximal system remains (syndetically) proximal.
\end{Remark}
\begin{Remark}\label{equiV}
Both relations, $\text{Prox}$ and $\text{SyProx}$, are reflexive and symmetric. Moreover, $\text{SyProx}$ is always an equivalence relation, whereas $\text{Prox}$ need not be an equivalence relation.
\end{Remark}
\begin{Remark}\label{pr:44}
It is easy to see that if $T$ is syndetically proximal then $T^{\times n}$ is syndetically proximal for each $n\geq 1$.
\end{Remark}

\begin{Prop}[\cite{MR0154269,MR0179775}, see also Theorem 19 in~\cite{MR3205496}]\label{k1}
The following are equivalent:
\begin{itemize}
\item
$\text{Prox}(T)$ is an equivalence relation,
\item
$\text{Prox}(T)=\text{SyProx}(T)$,
\item
the orbit closure of any point $(x,y)\in X\times X$ in the dynamical system $(T\times T,X\times X)$ contains exactly one minimal subset.
\end{itemize}
\end{Prop}
As an immediate consequence of \cref{equiV} and \cref{pr:4-}, we obtain:
\begin{Cor}\label{pr:4}
Suppose that $Tx_0=x_0$ and $\text{SyProx}(T)\cap (\{x_0\}\times X)=\{x_0\}\times X$. Then $\text{Prox}(T)\supset \text{SyProx}(T)=X\times X$, i.e.\ $T$ is syndetically proximal and $\{x_0\}$ is the unique minimal subset of $X$.
\end{Cor}

%%%%%%%%%%%%%%%%%%%%%%%%%%%%%%%%%%%%%%%%%%%%%%%%
%%%%%%%%%%%%%%%%%%%%%%%%%%%%%%%%%%%%%%%%%%%%%%%%
\subsubsection{General case}
%%%%%%%%%%%%%%%%%%%%%%%%%%%%%%%%%%%%%%%%%%%%%%%%
\begin{Prop}\label{LE51new}
Let $(T,X)$ be a topological dynamical system with a transitive point $\eta\in X$. The following are equivalent:
\begin{enumerate}[(a)]
\item\label{HA}
$(T,X)$ has a unique minimal subset $M$.
\item\label{HB}
There exists a closed, $T$-invariant subset $M'\subset X$ such that
for any $x\in M'$, $y\in X$, there exists $(m_n)_{n\geq 1}\subset \Z$ such that $T^{m_n}y\to x$.
\item\label{HC}
There exists $x\in X$ such that for any $y\in X$ there exists $(m_n)_{n\geq 1}\subset \Z$ such that $T^{m_n}y\to x$.
\item\label{HD}
There exists a closed, $T$-invariant subset $M''\subset X$, such that $\{k\in\Z : T^k\eta \in U\}$ is syndetic for any open set $U$ intersecting $M''$.
\item\label{HE}
There exists a sequence of open sets $(U_n)_{n\geq 1}\subset X$ such that:
\begin{itemize}
\item
$\text{diam}(U_n)\to 0\text{ as }n\to\infty$,
\item
$\{k\in\Z : T^k\eta \in U_n\}$ is syndetic.
\end{itemize}
\end{enumerate}
Moreover, if $x\in X$ is as in~\eqref{HC} then $x\in M$, where $M$ is the unique minimal subset of $X$ (in other words, $M$ is equal to the orbit closure of $x$).  Finally, if the above hold then $M'$ and $M''$ with the above properties are also unique and $M=M'=M''$. 
\end{Prop}
\begin{proof}
Suppose that~\eqref{HA} holds and take $x\in M':=M$ and $y\in X$. It follows by~\eqref{HA} that there exists $(m_n)_{n\geq 1}\subset\Z$ and $(x_n)_{n\geq 1}\subset M$  such that $d(T^{m_n}y,x_n)\to 0$ (otherwise, the orbit closure of $y$ would be disjoint from $M$ and would contain another minimal subset). We may assume without loss of generality that $x_n\to x_0\in M$, whence
$d(T^{m_n}y,x_0)\to 0$. Fix $\vep>0$. Let $k_0\in\Z$ be such that $d(T^{{k_0}}x_0,x)<\vep$. Moreover, let $\delta>0$ be sufficiently small, so that $d(z,z')<\delta$ implies $d(T^{{k_0}}z,T^{{k_0}}z')<\vep$ for $z,z'\in X$. Finally, let $\ov{m}\in \Z$ be such that $d(T^{\ov{m}}y,x_0)<\delta$. Then
$$
d(T^{\ov{m}+{k_0}}y,x)\leq d(T^{\ov{m}+{k_0}}y,T^{{k_0}}x_0)+d(T^{{k_0}}x_0,x)<2\vep.
$$
It follows that~\eqref{HB} holds.

Clearly,~\eqref{HB} implies~\eqref{HC}. We will show now that~\eqref{HC} implies~\eqref{HA}. Suppose that $M_1,M_2$ are minimal subsets of $X$. Let $x\in X$ be as in~\eqref{HC} and take $y_i\in M_i$, $i=1,2$. It follows by~\eqref{HC} that $x\in M_1\cap M_2$. This yields $M_1=M_2$.

We will show that~\eqref{HB} implies~\eqref{HD}. Let $U\subset X $ be an open set intersecting $M'':=M'$ and suppose that the orbit of $\eta$ visits $U$ with unbounded gaps. Then there exists $m_n\to \infty$ such that $T^{m_n+k}\eta \not\in U$ for $k\in \{-n,\dots,n\}$. Without loss of generality, we may assume that $T^{m_n}\eta \to y$. Then $T^{m_n+k}\eta\to T^ky \not\in U$ for each $k\in \Z$, i.e.\  the orbit of $y$ avoids $U$. Take $x\in M''\cap U$. It follows that the orbit of $y$ never approaches $x$. This contradicts~\eqref{HB}.

Clearly,~\eqref{HD} implies~\eqref{HE}. We will show now that~\eqref{HE} implies~\eqref{HC}. Suppose that~\eqref{HE} holds. Enlarging the sets $U_n$ if necessary, we may assume that $U_n=B(x_n,1/n)$ for some $x_n\in X$. Moreover, we may assume without loss of generality that $x_n\to x \in X$ as $n\to\infty$. Fix $y\in X$. For $n\geq 1$, let $d_n\geq 1$ be such that the orbit of $\eta$ visits $U_n$ with gaps at most $d_n$. Let $\delta_n>0$ be sufficiently small, so that $d(z,z')<\delta_n$ implies $d(T^mz,T^mz')<1/n$ for $0\leq m\leq d_n-1$. Let $k_n\in\Z$ be such that $d(y,T^{k_n}\eta)<\delta_n$. Finally, let $0\leq m_n\leq d_n-1$ be such that $T^{k_n+m_n}\eta\in U_n$, i.e.\ $d(T^{k_n+m_n}\eta,x_n)<1/n$. Then
\begin{multline*}
d(T^{m_n}y,x)\leq d(T^{m_n}y,T^{m_n+k_n}\eta)+d(T^{k_n+m_n}\eta,x_n)+d(x_n,x)\\
\leq 2/n+d(x_n,x)\to 0.
\end{multline*}
It follows that~\eqref{HC} indeed holds.

The above proof shows that $M\subset M'\subset M''$ where $M'$ and $M''$ are maximal sets with the above properties. Suppose now that $x\not\in M$ is as in~\eqref{HC} and take $y\in M$. Then $\inf_{z\in M}d(z,x)>0$. In particular, we cannot have $d(T^{m_n}y,x)\to 0$. In particular, this yields $M=M'$. We will show now the remaining equality $M=M''$. Let $U_n:=B(x,1/n)$, where $x\in M''$. The proof of implication \eqref{HE} $\Rightarrow$ \eqref{HC} yields that $x$ satisfies \eqref{HC} and we already know that this implies $x\in M$, i.e.\ $M=M''$.
\end{proof}
\begin{Remark}
Notice that the above result includes as a special case the characterization of minimal systems from \cref{k2}. Indeed, if $(T,X)$ is minimal then any open set $U$ intersects $M=X$, whence $\{n\in\Z : T^nx\in U\}$ is syndetic by~\eqref{HD}. On the other hand, if $\{n\in\Z : T^nx\in U\}$ is syndetic for any open set $U$, it follows that $M':=X$ satisfies \eqref{HD}. Therefore the only minimal subset $M$ is also equal to $X$, i.e.\ $(T,X)$ is minimal.
\end{Remark}
%%%%%%%%%%%%%%%%%%%%%%%%%%%%%%%%%%%%%%%%%%%%%%%%
\begin{Remark}\label{ausel1}
It follows by \cref{ausel} that if $(T,X)$ has a unique minimal subset $M$ then for any $x\in X$ there exists $y\in M$ such that $(x,y)\in\text{Prox}(T)$.
\end{Remark}
%%%%%%%%%%%%%%%%%%%%%%%%%%%%%%%%%%%%%%%%%%%%%%%%
\begin{Cor}\label{dlapod}
Let $(T,X)$ be a subshift. Then $(T,X)$ has a unique minimal subset $M$ if and only if there exists an infinite family of pairwise distinct blocks that appear on $\eta$ with bounded gaps.
\end{Cor}
\begin{proof}
This is an immediate consequence of the equivalence of \eqref{HA} and \eqref{HE} in \cref{LE51new}.
\end{proof}

%%%%%%%%%%%%%%%%%%%%%%%%%%%%%%%%%%%%%%%%%%%%%%%%
%%%%%%%%%%%%%%%%%%%%%%%%%%%%%%%%%%%%%%%%%%%%%%%%
If $(T,X)$ is a subshift, sometimes more can be said about the unique minimal subset. Namely, we have the following:
%%%%%%%%%%%%%%%%%%%%%%%%%%%%%%%%%%%%%%%%%%%%%%%%
\begin{Lemma}\label{topl}
Let $\eta\in\{0,1\}^\Z$. Suppose that there exist $B_n\in \{0,1\}^{[\ell_n,r_n]}$ for $n\geq 1$, with $\ell_n \searrow -\infty,r_n\nearrow \infty$, $(m_n)_{n\geq 1}\subset \Z$ and $(d_n)_{n\geq 1}\subset \N$, satisfying, for each $n\geq 1$:
\begin{enumerate}[(a)]
\item\label{WarA}
$d_{n}\divides d_{n+1}$, 
\item\label{WarB}
$d_n \divides m_{n+1}-m_n$,
\item\label{WarC}
$\eta[m_n+kd_n+\ell_n,m_n+kd_n+r_n]=B_n$ for each $k\in\Z$.\footnote{Conditions~\eqref{WarA},~\eqref{WarB} and~\eqref{WarC} imply that $B_{n+1} [\ell_n,r_n]=B_n,\ n\geq 1.$}
\end{enumerate}
Then $\eta$ has a Toeplitz sequence $x$ in its orbit closure $X_\eta$. 
\end{Lemma}
%%%%%%%%%%%%%%%%%%%%
%%%%%%%%%%%%%%%%%%%%
\begin{proof}
Fix $n_0\in \N$ and let $n\geq n_0$. Then, by~\eqref{WarA} and~\eqref{WarB}, we have $d_{n_0} \divides m_n-m_{n_0}$. Therefore, in view of~\eqref{WarC}, for any $k\in \Z$, we have
\begin{multline*}
S^{m_n}\eta[\ell_{n_0}+kd_{n_0},r_{n_0}+kd_{n_0}]\\
=\eta[m_{n_0}+(kd_{n_0}+m_n-m_{n_0})+\ell_{n_0},m_{n_0}+(kd_{n_0}+m_n-m_{n_0})+r_{n_0}]=B_{n_0}.
\end{multline*}
It follows that $x:=\lim_{n\to\infty}S^{m_n}\eta$ is well-defined and Toeplitz.
\end{proof}
\begin{Remark}
Suppose that the assumption of \cref{topl} are satisfied. It follows by \cref{dlapod} that $(S,X_\eta)$ has a unique minimal subset $M$ that is equal to the orbit closure of a Toeplitz sequence.
\end{Remark}

\subsection{Asymptotic densities}
For $A\subset \Z$, we recall several notions of asymptotic density (in fact, these are densities of the positive part of the set $A$, i.e.\ of $A\cap \N$). We have:
\begin{itemize}
\item
$\underline{d}(A):=\liminf_{N\to\infty}\frac1N|A\cap[1,N]| \text{ (\emph{lower density} of $A$)}$,
\item
$\overline{d}(A):=\limsup_{N\to\infty}\frac1N|A\cap[1,N]| \text{ (\emph{upper density} of $A$)}$.
\end{itemize}
If the lower and the upper density of $A$ coincide, their common value $d(A):=\un{d}(A)=\ov{d}(A)$ is called the \emph{density} of~$A$. We also have:
\begin{itemize}
\item
$\underline{\bdelta}(A):=\liminf_{N\to\infty}\frac{1}{\log N}\sum_{a\in A,1\leq a\leq N}\frac{1}{a} \text{ (\emph{lower logarithmic density} of $A$)}$,
\item
$\overline{\bdelta}(A):=\limsup_{N\to\infty}\frac{1}{\log N}\sum_{a\in A, 1\leq a\leq N}\frac{1}{a} \text{ (\emph{upper logarithmic density} of $A$)}$.
\end{itemize}
If the lower and the upper logarithmic density of $A$ coincide, we set $\bdelta(A):=\underline{\bdelta}(A)=\overline{\bdelta}(A)$ (\emph{logarithmic density} of $A$).

The following relations between the above notions are well-known:
\begin{equation}\label{gestos}
\underline{d}(A)\leq \underline{\bdelta}(A) \leq \overline{\bdelta}(A) \leq \overline{d}(A).
\end{equation}

\subsection{Sets of multiples, $\mathscr{B}$-free numbers and their density}\label{se:klasy}
For $\mathscr{B}\subset \N$, let
$
\cm_{\mathscr{B}}:=\bigcup_{b\in\mathscr{B}}b\Z \text{ and }\cf_{\mathscr{B}}:=\Z\setminus \cm_{\mathscr{B}}.
$
Sometimes, additional assumptions are put on $\mathscr{B}$.
\begin{Def}
We say that:
\begin{itemize}
\item
$\mathscr{B}$ is \emph{coprime}, if $\gcd(b,b')=1$ for $b\neq b'$ in $\mathscr{B}$,
\item
$\mathscr{B}$ is \emph{thin} if $\sum_{b\in\mathscr{B}}\nicefrac1{b}<+\infty$,
\item
$\mathscr{B}$ has \emph{light tails} if $\lim_{K\to\infty}\ov{d}\left( \bigcup_{b>K}b\Z\right)=0$,
\item
$\mathscr{B}$ is \emph{taut}~\cite{MR1414678} if for any $b\in \mathscr{B}$, we have
${\boldsymbol{\delta}}(\mathcal{M}_\mathscr{B})>\boldsymbol{\delta}(\mathcal{M}_{\mathscr{B}\setminus\{b\}})$.
\end{itemize}
\end{Def}
\begin{Remark}[see Chapter 0 in~\cite{MR1414678}]\label{stopkaprimitive}
Let $P(\mathscr{B})$ be the intersection of all sets $\mathscr{B}'\subset\N$ such that $\mathcal{M}_{\mathscr{B}}=\mathcal{M}_{\mathscr{B}'}$. Then $\mathcal{M}_{P(\mathscr{B})}=\mathcal{M}_{\mathscr{B}}$. Moreover, $P(\mathscr{B})$ is \emph{primitive} (i.e.\ no element of $P(\mathscr{B})$ divides any other). Therefore, throughout the paper, whenever $\mathscr{B}$ is arbitrary, we will tacitly assume that it is primitive.
\end{Remark}
\begin{Remark}\label{LLtauty2}
Since
$\ov{d}\left( \bigcup_{b>K}b\Z\right)\leq\sum_{b>K}\nicefrac1{b}$,
$$
\mathscr{B} \text{ is thin }\Rightarrow \mathscr{B}\text{ has light tails.}
$$
\end{Remark}
\begin{Def}
Following~\cite{MR1414678}, we say that $\mathscr{B}$ is \emph{Besicovitch} if $d(\cm_{\mathscr{B}})$ exists. Clearly, this is equivalent to the existence of $d(\cf_\mathscr{B})$.
\end{Def}
\begin{Remark}
Clearly, each finite $\mathscr{B}$ is Besicovitch.
\end{Remark}
Recall that $d(\mathcal{M}_\mathscr{B})$ may not exist -- the first counterexample was provided by Besicovitch~\cite{MR1512943}. Recall also the result by Erd\"os:
\begin{Th}[\cite{MR0026088}]\label{erd}
$\mathscr{B}=\{b_k : k\geq 1\}$ is Besicovitch if and only if
$$
\lim_{0<\vep\to 0}\limsup_{n\to \infty}\frac1n\sum_{n^{1-\vep}<b_k\leq n}|[0,n]\cap b_k\Z\cap\cf_{\{b_1,\ldots,b_{k-1}\}}|=0.
$$
\end{Th}
On the other hand, we have the following result of Davenport and Erd\"os:
\begin{Th}[\cite{Davenport1936,MR0043835}]\label{da-er}
For any $\mathscr{B}$, the logarithmic density $\bdelta(\cm_\sB)$ of $\cm_\mathscr{B}$ exists. Moreover,
\begin{equation}\label{i1}
\bdelta(\cm_\mathscr{B})=\un{d}(\cm_{\mathscr{B}})=\lim_{K\to\infty}d(\mathcal{M}_{\{b\in \mathscr{B} : b\leq K\}}).
\end{equation}
\end{Th}
\begin{Remark}\label{uwa9wrz}
Formula~\eqref{i1} follows from the proof of \cref{da-er} from~\cite{MR0043835} (see also~\cite{MR1414678}). Notice that~\eqref{i1} implies that $\mathscr{B}$ is Besicovitch if and only if
$$
\lim_{K\to \infty}\ov{d} (\mathcal{M}_{\{b\in\mathscr{B} : b>K\}} \setminus \mathcal{M}_{\{b\in\mathscr{B} : b\leq K\}})=0.
$$
In particular,
$$
\mathscr{B} \text{ has light tails }\Rightarrow \mathscr{B} \text{ is Besicovitch}.\footnote{This follows also by \cref{erd}.}
$$
\end{Remark}

We will need the following consequence of \cref{da-er}:
\begin{Cor}\label{daerdogolne}
Let $\mathscr{A}=\mathscr{A}_1\cup \mathscr{A}_2\cup\dots$ Then
$$
\bdelta(\mathcal{M}_{\mathscr{A}})=\underline{d}(\mathcal{M}_{\mathscr{A}})=\lim_{K\to \infty}\bdelta(\mathcal{M}_{\mathscr{A}_1\cup\mathscr{A}_2\cup\dots\cup\mathscr{A}_K}).
$$
\end{Cor}

\begin{proof}
Let
$$
\Delta(\sA):=\lim_{K\to \infty}\bdelta(\cm_{\sA_1\cup\dots \cup \sA_K}).
$$
Clearly, $\Delta(\sA)\leq  \bdelta(\cm_\sA)$.
We will show now that $\bdelta(\cm_\sA)\leq \Delta(\sA)$. For $K\geq 1$, let $N_K$ be such that
$$
{\{a\in \sA : a\leq K\}} \subset {\sA_1\cup\dots\cup \sA_{N_K}}.
$$
Using \cref{da-er}, we obtain
$$
\bdelta(\cm_\sA)=\lim_{K\to \infty}\bdelta(\cm_{\{a\in \sA : a\leq K\}}) \leq \lim_{K\to \infty}\bdelta(\cm_{\sA_1\cup \dots \cup \sA_{N_K}})=\Delta(\sA).
$$
This completes the proof.
\end{proof}

\begin{Remark}[Cf.\ \cref{uwa9wrz}]\label{uwa9wrz1}
Let $\mathscr{A}=\mathscr{A}_1\cup \mathscr{A}_2\cup\dots$ and suppose additionally that the density of $\mathscr{A}_1\cup\dots\cup \mathscr{A}_K$ exists, for each $K\geq 1$.
As a consequence of \cref{daerdogolne}, we obtain that $\mathscr{A}$ is Besicovitch if and only if
$$
\lim_{K\to \infty}\ov{d}(\mathcal{M}_\mathscr{A} \setminus \mathcal{M}_{\mathscr{A}_1\cup\dots \cup\mathscr{A}_K})=0.
$$
\end{Remark}

\begin{Def}
Following~\cite{MR1414678}, we say that $\mathscr{B}\subset \N\setminus \{1\}$ is \emph{Behrend} if $\boldsymbol{\delta}(\cm_{\mathscr{B}})=1$.
\end{Def}
\begin{Remark}\label{LLtauty1a}
Clearly, any superset of a Behrend set that does not contain $1$ remains Behrend. Moreover,
$$
\mathscr{B} \text{ is Behrend}\Rightarrow \mathscr{B}\text{ is Besicovitch}.
$$
Note also that by~\cref{da-er}, $\mathscr{B}\subset \N\setminus \{1\}$ is Behrend if and only if $d(\cm_{\mathscr{B}})=1$.
\end{Remark}
\begin{Prop}[\cite{MR1414678}, Corollary 0.14]\label{beh2}
$\mathscr{A}\cup \mathscr{B}$ is Behrend if and only if at least one of $\mathscr{A}$ and $\mathscr{B}$ is Behrend.
\end{Prop}

For $\mathscr{B}$, $a\in\N\setminus\{1\}$ let
$$
\mathscr{B}'(a):=\left\{\frac{b}{\gcd(b,a)}:b\in\mathscr{B}\right\}.
$$

\begin{Prop}[\cite{MR1414678}, Theorem 0.8]\label{beh}
Let $a\not\in\mathcal{M}_\mathscr{B}$. Then
$$
\boldsymbol{\delta}(\mathcal{M}_{\mathscr{B}\cup\{a\}})>\boldsymbol{\delta}(\mathcal{M}_\mathscr{B})
$$
if and only if $\mathscr{B}'(a)$ is not Behrend.
\end{Prop}
\begin{Prop}[\cite{MR1414678}, Corollary 0.19]\label{beh3}
$\mathscr{B}$ is taut if and only if it is primitive and does not contain $c\mathscr{A}$ with $\mathscr{A}\subset \N\setminus\{1\}$ that is Behrend.
\end{Prop}
\begin{Cor}\label{behnowy}
Suppose that $\mathscr{B}$ is taut. If $\boldsymbol{\delta}(\mathcal{M}_{\mathscr{B}\cup\{a\}})=\boldsymbol{\delta}(\mathcal{M}_\mathscr{B})$ then $a\in\mathcal{M}_\mathscr{B}$.
\end{Cor}
\begin{proof}
Suppose that $\boldsymbol{\delta}(\mathcal{M}_{\mathscr{B}\cup\{a\}})=\boldsymbol{\delta}(\mathcal{M}_\mathscr{B})$ and $a\not\in\mathcal{M}_\mathscr{B}$. By \cref{beh}, $\mathscr{B}'(a)$ is Behrend. Since $a$ has finitely many divisors, it follows by \cref{beh2} that at least one of the sets
$$
\mathscr{B}'_d(a):=\left\{\frac{b}{d} : b\in\mathscr{B} \text{ and }\gcd(b,a)=d\right\},
$$
where $d\divides a$, is Behrend. Moreover, $d\cdot \mathscr{B}'_d(a)\subset \mathscr{B}$. Notice that $1\not\in \mathscr{B}'_d(a)$. Indeed, if $1\in \mathscr{B}'_d(a)$ then $d=\gcd(d,a)\in \mathscr{B}$. In particular, $d\divides a$, i.e.\ $a\in\mathcal{M}_\mathscr{B}$, which is not possible by the choice of $a$. It follows by \cref{beh3} that $\mathscr{B}$ cannot be taut. This contradicts the assumptions and the result follows.
\end{proof}

The following is an immediate consequence of \cref{beh3}:
\begin{equation}\label{tautnotbeh}
\mathscr{B} \text{ is taut } \Rightarrow \mathscr{B}\text{ is not Behrend, unless }\mathscr{B}=\{1\}.
\end{equation}
Furthermore, notice that
\begin{equation}\label{LLtauty1}
\mathscr{B}\text{ has light tails (and is primitive) }\Rightarrow \mathscr{B} \text{ is taut}.
\end{equation}
Indeed, if $\mathscr{B}$ is not taut, by \cref{beh3}, we have that $\mathscr{B}\supset c\mathscr{A}$ with $\mathscr{A}$ Behrend. Moreover, given $K\geq 1$, there exists $L=L(K,c)$ such that
$$
c\cdot \{a\in \sA : a>L\}\subset \bigcup_{b>K}b\Z.
$$
But, in view of \cref{beh2}, $ \{a\in \sA : a>L\}$ is Behrend. It follows that $\boldsymbol{\delta}(\bigcup_{b>K}b\Z)\geq 1/c$ for all $K\geq 1$, which means that $\sB$ cannot have light tails. In particular, we obtain
$$
\mathscr{B} \text{ is finite }\Rightarrow \mathscr{B}\text{ is taut}.
$$

\subsection{Canonical odometer associated with $\mathscr{B}$}\label{canonical}

To simplify the notation we will now restrict ourselves to the case when $\mathscr{B}$ is infinite and we will denote the elements of $\mathscr{B}$ by $b_k$, $k\geq 1$ (if $\mathscr{B}$ is finite similar objects can be defined, with obvious changes).

Consider the compact Abelian group $ G _\mathscr{B}:=\prod_{k\geq 1}\Z/b_k\Z$,
with the coordinatewise addition. The product topology on $ G _\mathscr{B}$ is metrizable with a (bounded) metric $d$ given by
\begin{equation}\label{stop}
d( g , g ')=\sum_{k\geq1}\frac1{2^k}\frac{| g _k- g '_k|}{1+| g _k- g '_k|}.
\end{equation}
Let $\PP_{ G _\mathscr{B}}$ be Haar measure of $ G _\mathscr{B}$, i.e.\ $\PP_{ G _\mathscr{B}}=\bigotimes m_{\Z/b_k\Z}$.\footnote{For $c\in\N$, $m_{\Z/c\Z}$ stands for the counting measure on $\Z/c\Z$.}
For $n\in\Z$, let
\begin{equation}
\underline{n}_\mathscr{B}:=(n\bmod b_1,n\bmod b_2,\ldots)\in  G _\mathscr{B}.
\end{equation}
Denote by $ G$ the smallest closed subgroup of $ G _\mathscr{B}$ that contains $\underline{1}_\mathscr{B}$, i.e.
\begin{equation}\label{GIE}
 G:=\overline{\{\underline{n}_\mathscr{B} : n\in\Z\}}\subset  G _\mathscr{B}.
\end{equation}
\begin{Remark}
By its definition, $ G\subset G _\mathscr{B}$ contains a dense cyclic subgroup, i.e.\ $ G$ is monothetic and the homeomorphism
\begin{equation}\label{homeo}
T g = g +\underline1_\mathscr{B}
\end{equation}
yields a uniquely ergodic dynamical system $(T, G)$ (with Haar measure $\PP$ as the only invariant measure).
\end{Remark}
We will now provide another model of $(T,G)$. First, given $1\leq k<\ell$, denote by
$$
\pi_{k,\ell}\colon\Z/\lcm(b_1,\ldots, b_k, \ldots, b_{\ell})\Z\to
\Z/\lcm(b_1,\ldots, b_k)\Z
$$
the natural homomorphism given, for each $r\in\Z/\lcm(b_1,\ldots, b_k, \ldots, b_{\ell})\Z$, by
\begin{equation}\label{dze2}
\pi_{k,\ell}(r)= r \bmod \lcm(b_1,\ldots, b_k).
\end{equation}
Note that whenever $1\leq k<\ell<m$,
\begin{equation}\label{dze3}
\pi_{k,\ell}\circ\pi_{\ell,m}=\pi_{k,m}.
\end{equation}
Also, for each $k\geq1$, we set
$$
\pi_{k}:=\pi_{k,k+1}.
$$
This yields an inductive system
$$
\Z/b_1\Z\stackrel{\pi_1}{\leftarrow}\Z/\lcm(b_1, b_2)\Z\stackrel{\pi_2}{\leftarrow}\ldots\stackrel{\pi_{k-1}}{\leftarrow}\Z/\lcm(b_1,\ldots, b_k)\Z\stackrel{\pi_k}{\leftarrow}\ldots
$$
and we define
\begin{multline}\label{gie}
 G' :=\varprojlim\Z/\lcm(b_1,\ldots, b_k)\Z\\
 =\left\{ g \in\prod_{k\geq 1}\Z/\lcm(b_1,\ldots, b_k)\Z :\pi_k( g _{k+1})= g _k \; \mbox{for each}\;k\geq1\right\},
\end{multline}
where $ g =( g _1, g _2,\ldots)$. Then $G'$ is closed and invariant under the coordinatewise addition. Hence, $G'$ is Abelian, compact and metrizable, cf.~\eqref{stop}. We denote by $\PP'$ Haar measure on $G'$. Note that in view of~(\ref{dze3}), for each $n\geq1$, we have
\begin{equation}\label{dze5}
\underline{n}:=(n\bmod b_1,n\bmod \lcm(b_1,b_2),\ldots)\in G' ,
\end{equation}
in particular, $\underline1\in G' $. On $G'$, we also define a homeomorphism:
\begin{equation}
T' g = g +\underline{1}.
\end{equation}
\begin{Remark}\label{monothe}
Notice that if $(g_1,g_2,\dots)\in G'$ then, since $g_k=g_j \bmod \lcm(b_1,\dots, b_j)$ for $j=1,\dots,k$, we have
$$
(g_k,g_k,\dots)\to (g_1,g_2,\dots) \text{ when }k\to \infty.
$$
It follows that $\{\un{n}: n\in \Z\}$ is dense in $G'$ (and hence $G'$ is monothetic).
\end{Remark}
\begin{Lemma}\label{izomo}
The map $W\colon \{\un{n}_\mathscr{B} : n\in\Z\} \to G'$ given by $W(\un{n}_\mathscr{B})=\un{n}$ extends continuously to $G$ in a unique way. Moreover, it yields a topological isomorphism of the dynamical systems $(T,G)$ and $(T',G')$.
\end{Lemma}
\begin{proof}
Notice first that $W$ is uniformly continuous (and equivariant). Indeed, for any $K\geq 1$, such that if $d(\un{n}_\mathscr{B},\un{m}_\mathscr{B})$ is sufficiently small then $n=m\bmod b_k$ for $1\leq k\leq K$. It follows that $n=m\bmod \lcm(b_1,\dots,b_k)$ for $1\leq k\leq K$, i.e.\ $d(\un{n},\un{m})$ is small, provided that $K$ is large. Therefore, $W$ extends to a continuous map from $G$ to $G'$. Moreover, by \cref{monothe}, $W\colon G\to G'$ is surjective.

It remains to show that $W$ is injective. For this, it suffices to show that the map $\un{n}\mapsto \un{n}_\mathscr{B}$ is also uniformly continuous. Fix $K\geq 1$. If $d(\un{n},\un{m})$ is sufficiently small then then $n=m\bmod \lcm(b_1,\dots,b_k)$ for $1\leq k\leq K$. It follows clearly that, for $1\leq k\leq K$, we have $n=m\bmod b_k$, i.e. $d(\un{n}_\mathscr{B},\un{m}_\mathscr{B})$ is arbitrarily small, provided that $K$ is large. This completes the proof.
\end{proof}

\begin{Def}
We say that $(T,G,\PP)$ is the \emph{canonical odometer} associated to $\mathscr{B}$.
\end{Def}

\begin{Remark}\label{uw215}
It follows by the proof of the above lemma that for $g \in  G$, we have
\begin{equation}\label{formulaW}
W( g )=( g _1\bmod b_1,  g _2\bmod b_2,\dots).
\end{equation}
\end{Remark}

\begin{Example}\label{coinci}
When $\mathscr{B}$ is coprime then
$
\Z/\lcm(b_1,\ldots,b_k)\Z=\Z/(b_1\cdot\ldots\cdot b_k\Z)$ is, by the Chinese Remainder Theorem, canonically isomorphic to $\Z/b_1\Z\times\ldots\times\Z/b_k\Z$ via 
$$
j\mapsto (j\bmod  b_1,\ldots,j\bmod b_k),
$$
so $\pi_k$ corresponds to
$$
{\rm proj}_{k}\colon\Z/b_1\Z\times\ldots\times\Z/b_k\Z\times \Z/b_{k+1}\Z\to
\Z/b_1\Z\times\ldots\times\Z/b_k\Z,
$$
i.e.\ the projection on the $k$ first coordinates. The inverse limit $G'$ given by the system $\{{\rm proj}_k : k\geq1\}$ is naturally identified with the direct product $ G _\mathscr{B}$. Moreover,
$\underline1\in G' $ corresponds to $\underline{1}_\mathscr{B}\in G _\mathscr{B}$. It follows that $G=G_\mathscr{B}$ and thus the canonical odometer associated to $\mathscr{B}$ is the same as in~\cite{Ab-Le-Ru} whenever $\mathscr{B}$ is coprime.
\end{Example}

We will now show that the canonical odometer ``outputs'' $\cf_\mathscr{B}$. Consider the following sets:
\begin{align}
C:=&\{( g _1, g _2,\ldots)\in G: \text{ for all }k\geq 1,\ g _k\not\equiv 0\bmod b_k\},\label{tujestC} \\
C':=&\{( g _1, g _2,\ldots)\in G': \text{ for all }k\geq 1,\ g _k\not\equiv0\bmod b_k\}.
\end{align}
\begin{Remark}\label{dzl1}
By \cref{uw215}, we have $W(C)=C'$. In particular, for each $n\in\Z$, we have
$$
\underline n_\mathscr{B}\in C \iff \underline{n}\in C'  \iff n\in\cf_{\mathscr{B}}.
$$
\end{Remark}
Let $\eta\in\{0,1\}^{\Z}$ be the sequence corresponding to $\raz_{\cf_{\mathscr{B}}}$. Denote by $X_\eta$ the subshift generated by $\eta$, i.e.
$$
X_\eta:=\{x\in\{0,1\}^{\Z}: \mbox{each block appearing on $x$ appears on $\eta$}\}.
$$
In other words, $X_\eta=\overline{\{S^k\eta: k\in\Z\}}$, where $S$ stands for the shift transformation. We call $(S,X_\eta)$ the \emph{$\mathscr{B}$-free subshift}.

Define
$\va\colon G \to \{0,1\}^{\Z}$
by setting $\varphi( g )(n):=\raz_{C}(T^n g )$ and notice that
\begin{equation}\label{dze8equi}
\va( g )(n)=1\iff n\not\equiv -g_k \bmod b_k \text{ for all }k\geq 1.
\end{equation}
Finally, notice that
\begin{equation}\label{fikomut}
\varphi\circ T = S\circ \varphi.
\end{equation}
and
$
\eta=\varphi(0,0,\dots).
$

\subsection{Admissibility}
\begin{Def}[\cite{sarnak-lectures, Ab-Le-Ru}]
We call a sequence $x\in\{0,1\}^{\Z}$ {\em admissible} (or $\mathscr{B}$-{\em admissible}) if $|\text{supp }x \bmod b|<b$ for each $b\in\mathscr{B}$.
We denote by $X_{\mathscr{B}}$ the subshift of admissible sequences (it is easy to check that $X_\mathscr{B}$ is closed and $S$-invariant). We call $(S,X_\mathscr{B})$ the \emph{$\mathscr{B}$-admissible subshift}.
\end{Def}
\begin{Remark}\label{Rk220}
Consider $\va_\mathscr{B}\colon G _\mathscr{B}\to \{0,1\}^\Z$ given, for $ g \in G _\mathscr{B}$, by the same formula as in~\eqref{dze8equi}.  Arguing as in \cite{Ab-Le-Ru}, we easily obtain $\va_{\mathscr{B}}( G _\mathscr{B})\subset X_{\mathscr{B}}$. In particular, since $\eta=\varphi_{\mathscr{B}}(0,0,\dots)$, we have $\eta\in X_\mathscr{B}$, so
\begin{equation}\label{jm18}
X_\eta\subset X_{\mathscr{B}}.
\end{equation}
\end{Remark}
\begin{Def}[cf.\ \cite{MR2317754,MR3007694}]
We say that $X\subset \{0,1\}^\Z$ is \emph{hereditary} if for $x\in X$ and $y\in \{0,1\}^\Z$ with $y\leq x$ (coordinatewise), we have $y\in X$.
\end{Def}
It follows directly from the definition of admissibility that
\begin{equation}\label{jm19}
\mbox{$X_{\mathscr{B}}$ is hereditary}.\end{equation}
Denote by $\widetilde{X}_\eta$ the smallest hereditary subshift containing $X_\eta$. In view of~(\ref{jm18}) and~(\ref{jm19}),
\begin{equation}\label{jm20}
X_\eta\subset \widetilde{X}_\eta\subset X_{\mathscr{B}}.\end{equation}
\begin{Remark}\label{nieprz}
Note that $X_\mathscr{B}$ is always uncountable. Indeed, for $\mathscr{B}$ infinite, it suffices to notice that
$$
A:=\{b_1\cdot\ldots b_k : k\geq 1\} \text{ is $\mathscr{B}$-admissible and infinite}
$$
and apply~\eqref{jm19} (the set $\{x\in \{0,1\}^\Z : x\leq \raz_A\}$ is uncountable). If $\mathscr{B}=\{b_1,\dots,b_k\}$ is finite then 
$$
A:=\{(b_1\cdot\ldots\cdot b_k)^\ell: \ell\geq 1\} \text{ is $\mathscr{B}$-admissible and infinite}
$$
and we again apply~\eqref{jm19}.
\end{Remark}
For $\mathscr{B}$ infinite, coprime and thin we have $X_\eta=X_{\mathscr{B}}$, see \cite{Ab-Le-Ru}. This need not always be the case:
%%%%
\begin{Example}[$X_\eta\subsetneq \widetilde{X}_\eta\subsetneq X_{\mathscr{B}}$]\label{ex:2.4}
Let $\mathscr{B}:=\mathscr{P}$, i.e. $\sB$ is the set of all primes. Then $\cf_\mathscr{B}=\{\pm 1\}$. It follows that
\begin{align*}
X_\eta&=\{S^n \eta : n\in\Z\} \cup\{(\dots,0,0,0,\dots)\},\\
\widetilde{X}_\eta&=\{S^n \eta : n\in\Z\}\cup \{S^n(\dots,0,1,0,\dots)\} \cup\{(\dots,0,0,0,\dots)\}.
\end{align*}
In particular, $X_\eta \subsetneq \widetilde{X}_\eta$ and both these sets are countable. Moreover, by \cref{nieprz}, we have $\widetilde{X}_\eta\subsetneq X_\mathscr{B}$.
\end{Example}
\begin{Remark}
The set $\mathscr{B}$ from Example~\ref{ex:2.4} is Behrend.
\end{Remark}
\begin{Example}[$\widetilde{X}_\eta\subsetneq X_{\mathscr{B}}$, see \cref{Rk21}]\label{X_eta<>X_B}
Suppose that $4,6\in \mathscr{B}$ and $b>12$ for $b\in\mathscr{B}\setminus \{4,6\}$.
Let $y\in \{0,1\}^\Z$ be such that 
$$
y[1,12]=110011100110
$$
and $y(n)=0$ for all $n\in\Z\setminus\{1,2,\ldots,12\}$. It follows that $y\in X_\mathscr{B}$. We claim that $y\not\in\widetilde{X}_\eta$. Suppose that 
\begin{equation}\label{L1}
y[1,12]\leq\eta[k,k+11]\text{ for some $k\in\Z$.}
\end{equation}
Recall that $4\in\mathscr{B}$. Since $y[1]=\eta[k]=y[2]=\eta[k+1]=1$, it follows that $4\divides k+2$ or $4\divides k+3$. Since $y[7]=\eta[k+6]=1$, we cannot have $4\divides k+2$. Hence $4\divides k+3$. On the other hand, we have $6\in\mathscr{B}$. Since $y[i+1]=\eta[k+i]=1$ for $i\in \{0,1,4,5,6\}$ and $k+2$ is odd, we have $6\divides k+3$. It follows that $6\divides k+9$, whence $\eta[k+9]=0$. This, however, contradicts~\eqref{L1}.
\end{Example}
\begin{Remark}\label{Rk21}
In Example~\ref{X_eta<>X_B}, $\mathscr{B}$ can be chosen so that the density of $\mathcal{F}_\mathscr{B}$ exists and $d(\mathcal{F}_\mathscr{B})>0$. We will see in Section~\ref{s6} that, by imposing additional conditions on $\mathscr{B}$ from this example, one can obtain both $X_\eta=\widetilde{X}_\eta\subsetneq X_{\mathscr{B}}$ and $X_\eta\subsetneq\widetilde{X}_\eta\subsetneq X_{\mathscr{B}}$ (and still have $d(\mathcal{F}_\mathscr{B})>0$), see \cref{X_eta<>X_Ba}.
\end{Remark}

\begin{Example}[$X_\eta\subsetneq \widetilde{X}_\eta= X_{\mathscr{B}}$, cf.\ Question~\ref{q2}, page~\pageref{q2}]\label{ex:2.4a}
Let $\mathscr{B}:=\{2\}$. Then $X_{\eta}=\{\eta, S\eta\}\subsetneq \widetilde{X}_\eta$. Take $x\in X_\mathscr{B}$. Then either $\text{supp }x\bmod 2\subset 2\Z$ or $\text{supp }x\bmod 2\subset2\Z+1$. In other words, $\text{supp }x\bmod 2 \subset \text{supp }\eta$ or $\text{supp }x\bmod 2\subset \text{supp }S\eta$, which means that $\widetilde{X}_\eta=X_{\mathscr{B}}$.
\end{Example}

The subshift $(S,\widetilde{X}_\eta)$ has some natural $S$-invariant subsets we will be interested in. To study them, first, for $0\leq s_k\leq b_k$, $k\geq 1$, let
\begin{align*}
Y_{s_1,s_2,\dots}&:=\{x\in \{0,1\}^\Z : |\text{supp }x \bmod b_k|=b_k-s_k \text{ for each }k\geq 1\},\\
Y_{\geq s_1,\geq s_2,\dots}&:=\{x\in \{0,1\}^\Z : |\text{supp }x \bmod b_k|\leq b_k-s_k \text{ for each }k\geq 1\}
\end{align*}
(if $\mathscr{B}$ is finite, we define analogous subsets, with obvious changes). 
\begin{Remark}
For $0\leq s_k\leq b_k$, $k\geq 1$, define auxiliary subsets
\begin{align*}
Y^k_{s_k}&:=\{x\in \{0,1\}^\Z : |\text{supp }x\bmod b_k|=b_k-s_k\},\\
Y^k_{\geq s_k}&:=\{x\in \{0,1\}^\Z: |\text{supp }x\bmod b_k|\leq b_k-s_k\}.
\end{align*}
Then $Y^k_{s_k}=Y^k_{\geq s_k}\setminus Y^k_{\geq s_k+1}$ and $Y^k_{\geq s_k}, Y^k_{\geq s_k+1}$ are closed. Moreover
$$
Y_{s_1,s_2,\dots}=\bigcap_{k\geq 1}Y^k_{s_k},\ Y_{\geq s_1,\geq s_2,\dots}=\bigcap_{k\geq 1}Y^k_{\geq s_k}.
$$
In particular, $Y_{s_1,s_2,\dots}$ is Borel and $Y_{\geq s_1,\geq s_2,\dots}$ is closed, for any choice of $0\leq s_k\leq b_k$, $k\geq 1$.
Additionally, sets $Y_{s_1,s_2,\dots}$ are pairwise disjoint for different choices of $(s_1,s_2,\dots)$ and
$$
\{0,1\}^\Z=\bigcup_{0\leq s_k\leq b_k, k\geq 1}Y_{s_1,s_2,\dots}.
$$
We will write $Y$ for $Y_{1,1,\dots}$. Notice also that $Y_{\geq s_1,\geq s_2,\dots}$ is the smallest hereditary subshift containing $Y_{s_1,s_2,\dots}$. 
\end{Remark}
Following~\cite{Peckner:2012fk}, we define a map $\theta\colon Y\cap \widetilde{X}_\eta\to  {G}_\mathscr{B}$ by
\begin{equation}\label{THETA}
\theta(y)=g \iff (\text{supp }y) \cap (b_k\Z-g_k)=\emptyset \text{ for each }k\geq 1.
%-g_k\not\in\text{supp }y\bmod b_k \text{ for each }k\geq 1.
\end{equation}
Notice that given $y\in Y$ and $k_0\geq 1$, there exists $N\geq 1$ such that
\begin{equation}\label{wksk}
|(\text{supp }y)\cap [-N,N]\bmod b_k|=b_k-1 \text{ for }1\leq k\leq k_0
\end{equation}
\begin{Remark}
Notice that
\begin{equation}\label{jestwG}
\theta(Y\cap \widetilde{X}_\eta)\subset G.
\end{equation}
Indeed, take $y\in Y\cap \widetilde{X}_\eta$. Given $k_0\geq 1$, let $N\geq 1$ be such that~\eqref{wksk} holds and let $M\in\Z$ be such that $y[-N,N] \leq \eta[-N+M,N+M]$. It follows that $\theta(y)=(g_1,g_2,\dots)$, where $g_k\equiv -M\bmod b_k$ for $1\leq k\leq k_0$. This yields~\eqref{jestwG}.
\end{Remark}
\begin{Remark}\label{cgtheta}
Note also that $\theta$ is continuous. Indeed, given $y\in Y$ and $k_0\geq 1$, let $N$ be such that~\eqref{wksk} holds. Then, if $y'\in Y$ is sufficiently close to $y$ then~\eqref{wksk} holds for $y'$ as well. Therefore, if $y_n\to y$ in $Y$ then $\theta(y_n)\to \theta(y)$.
\end{Remark}

\begin{Remark}\label{1.2.6}
Note that:
\begin{itemize}
	\item\label{F1}
	$T\circ \theta=\theta \circ S$,
	\item\label{F2}
	for each $y\in Y\cap \widetilde{X}_\eta$, $y\leq \varphi( \theta (y))$,
	\item for any $\nu\in\mathcal{P}(S,Y\cap \widetilde{X}_\eta)$, $\theta_\ast(\nu)=\PP$
\end{itemize}
(the first two properties follow by a direct calculation, the third one is a consequence of the unique ergodicity of $T$).
\end{Remark}

\subsection{Mirsky measure $\nu_\eta$}
\begin{Def}
The image $\nu_{\eta}:=\va_\ast(\PP)$ of $\PP$ via $\va$ is called the {\em Mirsky measure} of $\mathscr{B}$.
\end{Def}

\begin{Remark}[cf.\ \cref{coinci}]
In the previous works \cite{Ab-Le-Ru,MR3356811}, the Mirsky measure was defined in a different way. In the new notation, the ``old Mirsky measure'' was given by ${\nu}_{\mathscr{B}}:=({\va_\mathscr{B}})_\ast(\PP_{\Omega_\mathscr{B}})$. We
$$
\nu_\sB(\{x\in \{0,1\}^\Z : x(0)=1\})=\prod_{b\in\sB}\left(1-\frac{1}{b}\right)
$$
(we follow word for word the proof of this formula from \cite{Ab-Le-Ru}). This implies that $\nu_\sB\neq\delta_{(\dots,0,0,0,\dots)}$ if and only if $\sB$ is thin. An advantage of $\nu_\eta$ is that $\nu_\eta\neq\delta_{(\dots,0,0,0,\dots)}$ whenever $\sB\subset \N$ is not Behrend (see \cref{kiedybeh}). Moreover, we will see, that $\nu_\eta$ plays a similar role and has similar properties as the ``old Mirsky measure''. This is why we call $\nu_\eta$ the Mirsky measure, not~$\nu_\mathscr{B}$. Notice that if $\sB$ is infinite, coprime and thin, we have $\nu_\eta=\nu_\mathscr{B}$.
\end{Remark}
%%%%%%%%%%%%%%%%%%%%%%%%%%%%%%%%%%%%%%%%%%%%%%%%

\section{Topological dynamics}
\subsection{Unique minimal subset (proof of \cref{TTA})}\label{proofTTA}
In the square-free case, i.e.\ when $\mathscr{B}=\{p^2 : p\in\mathcal{P}\}$, the subshift $(S,X_\eta)$ is proximal~\cite{sarnak-lectures}. In particular, by \cref{pr:4-}, it has a fixed point that yields the only minimal subset of $X_\eta$.\footnote{This fixed point is the sequence $(\dots,0,0,0,\dots)$.} It turns out that in general, even though there are $\mathscr{B}$-free subshifts $(S,X_\eta)$ that are not proximal,\footnote{This happens, e.g., when $\mathscr{B}$ is finite, we will see more examples later (we give necessary and sufficient conditions for proximality in \cref{TTBsekcja}). } the following holds:
\begin{Prop}\label{PR51}
For any $\mathscr{B}\subset \N$, $(S,X_\eta)$ has a unique minimal subset.
\end{Prop}
%%%%%%%%%%%%%%%%%%%%%%%%%%%%%%%%%%%%%%%%%%%%%%%%
\begin{proof}
We apply \cref{dlapod} to $\eta=\raz_{\mathcal{F}_\mathscr{B}}$. Suppose first that on $\eta$ there are arbitrarily long blocks consisting of zeros. Since each zero appears on $\eta$ with some period, it follows that each such block appears on $\eta$ with bounded gaps. Applying \cref{dlapod} to $\eta=\raz_\mathcal{F}$, we conclude that $(S,X_\eta)$ has a unique minimal subset.

Suppose now that the length of blocks consisting of zeros that appear on $\eta$ is bounded. The sequence $(B_n)_{n\geq 1}$ necessary to apply \cref{dlapod} will be defined inductively. Let $B_1$ be the {\bf longest} block of zeros appearing on $\eta$. Suppose that $B_1,\dots, B_n$ are chosen. For $n$ odd, let $B_{n+1}$ be the {\bf shortest} possible block of the form $B_n1\ldots 1$ that appears on $\eta$. For $n$ even, let $B_{n+1}$ be the {\bf longest} possible block of the form $B_n0\ldots 0$ that appears on $\eta$. Now, it suffices to show that each $B_n$, for $n$ even, appears on $\eta$ with bounded gaps. Since each zero appears on $\eta$ with some period, it follows that the pattern of zeros from $B_n$ appears on $\eta$ along some infinite arithmetic progression. Moreover, by the choice of $B_{k+1}$ for $k$ odd, it follows that whenever we see the pattern of zeros from $B_n$ on $\eta$, we actually see a copy of block $B_n$ at the same position on $\eta$. The result follows by \cref{dlapod}.
\end{proof}
%%%%%%%%%%%%%%%%%%%%%%%%%%%%%%%%%%%%%%%%%%%%%%%%
\begin{Prop}\label{PR52}
For any $\mathscr{B}\subset \N$, there exists a Toeplitz sequence $x\in X_\eta$.
\end{Prop}
\begin{proof}
Suppose first that on $\eta$ there are arbitrarily long blocks consisting of zeros. Then the Toeplitz sequence $(\dots,0,0,0,\dots)$ is in $X_\eta$.

Suppose now that the length of blocks consisting of zeros that appear on $\eta$ is bounded. We will use \cref{topl} and $(B_n)_{n\in\N}, (m_n)_{n\in\N}, (d_n)_{n\in\N}$ will be constructed inductively. First, we will choose the longest block of zeros that appears on $\eta$. Then we will extend it to the right and to the left by the shortest possible blocks of ones such that the extended block appears on $\eta$. Next, the obtained block will be extended to the right and then to the left by the longest possible blocks of zeros, so that the block we obtain still appears on $\eta$. This procedure will be repeated to obtain longer and longer blocks.

Let $B_1$ be the longest block of zeros that appears on $\eta$. For convenience, we will treat $B_1$ as an element of $\{0,1\}^{[0,|B_1|-1]}$ (i.e.\ we set $\ell_1:=0$, $r_1:=|B_1|-1$). Then, since $\eta=\raz_{\mathcal{F}_\mathscr{B}}$, there exists $d_1\in\N$ such that $B_1$ appears on $\eta$ periodically, with period $d_1$, i.e., for some $m_1\in\Z$, we have
$$
\eta[m_1+kd_1+\ell_1,m_1+kd_1+r_1]=B_1 \text{ for each }k\in\Z.
$$
Suppose now that $B_n\in \{0,1\}^{[\ell_n,r_n]}$, $m_n\in\Z$, $d_n\in\N$ for $1\leq n\leq 4n_0+1$ are chosen so that~\eqref{WarA} and~\eqref{WarB} from \cref{topl} hold for $1\leq n\leq 4n_0$ and~\eqref{WarC} from \cref{topl} holds for $1\leq n\leq 4n_0+1$. We will now define $B_n\in\{0,1\}^{[\ell_{n},r_{n}]}$, $m_n\in\Z$, $d_n\in\N$ for $4n_0+2\leq n\leq 4n_0+5$. 

Let $B_{4n_0+2}\in \{0,1\}^{[\ell_{4n_0+2},r_{4n_0+2}]}$, where $\ell_{4n_0+2}=\ell_{4n_0+1}$ (and $r_{4n_0+2}=\ell_{4n_0+2}+|B_{4n_0+2}|-1$), be the {\bf shortest} block of the form $B_{4n_0+1}1\ldots 1$ that appears on $\eta$ and begins at position $m_{4n_0+1}+\ell_{4n_0+1}+k_0d_{4n_0+1}$ for some $k_0\in\Z$, i.e.\
$$
\eta[m_{4n_0+2}+\ell_{4n_0+2},m_{4n_0+2}+r_{4n_0+2}]=B_{4n_0+2},
$$
where $m_{4n_0+2}=m_{4n_0+1}+k_0d_{4n_0+1}$. Then, clearly, $d_{4n_0+1}\divides m_{4n_0+2}-m_{4n_0+1}$. Moreover, by the definition of $B_{4n_0+2}$, we have
$$
\eta[m_{4n_0+2}+\ell_{4n_0+2}+kd_{4n_0+1},m_{4n_0+2}+r_{4n_0+2}+kd_{4n_0+1}]=B_{4n_0+2}
$$
for each $k\in\Z$, i.e.\ we may set $d_{4n_0+2}:=d_{4n_0+1}$. This way, we have extended our block $B_{4n_0+1}$ to the right by a block of ones.

The block $B_{4n_0+3}$ is defined in a similar way as $B_{4n_0+2}$, but now we extend $B_{4n_0+2}$ to the left. Let $B_{4n_0+3}\in \{0,1\}^{[\ell_{4n_0+3},r_{4n_0+3}]}$, where $r_{4n_0+3}=r_{4n_0+2}$ (and $\ell_{4n_0+3}=r_{4n_0+3}-|B_{4n_0+3}|+1$), be the {\bf shortest} block of the form $1\ldots 1B_{4n_0+2}$ that appears on $\eta$ and ends at position $m_{4n_0+2}+r_{4n_0+2}+k_0d_{4n_0+2}$ for some $k_0\in\Z$, i.e.\
$$
\eta[m_{4n_0+3}+\ell_{4n_0+3},m_{4n_0+3}+r_{4n_0+3}]=B_{4n_0+3},
$$
where $m_{4n_0+3}=m_{4n_0+2}+k_0d_{4n_0+2}$. Then, clearly, $d_{4n_0+2}\divides m_{4n_0+3}-m_{4n_0+2}$. Moreover, by the definition of $B_{4n_0+3}$, we have
$$
\eta[m_{4n_0+3}+\ell_{4n_0+3}+kd_{4n_0+3}, m_{4n_0+3}+r_{4n_0+3}+kd_{4n_0+3}]=B_{4n_0+3}
$$
for each $k\in\Z$, i.e.\ we may set $d_{4n_0+3}:=d_{4n_0+2}$. This way, we have extended our block $B_{4n_0+2}$ to the left by a block of ones.

Let $B_{4n_0+4}\in \{0,1\}^{[\ell_{4n_0+4},r_{4n_0+4}]}$, where $\ell_{4n_0+4}=\ell_{4n_0+3}$ (and $r_{4n_0+4}=\ell_{4n_0+4}+|B_{4n_0+4}|-1$), be the {\bf longest} block of the form $B_{4n_0+3}0\ldots 0$ that appears on $\eta$ and begins at position $m_{4n_0+3}+\ell_{4n_0+3}+k_0d_{4n_0+3}$ for some $k_0\in\Z$, i.e.\
$$
\eta[m_{4n_0+4}+\ell_{4n_0+4},m_{4n_0+4}+r_{4n_0+4}]=B_{4n_0+4},
$$
where $m_{4n_0+4}=m_{4n_0+3}+k_0d_{4n_0+3}$. Then, clearly, $d_{4n_0+3}\divides m_{4n_0+4}-m_{4n_0+3}$. Moreover, since each zero on $\eta$ appears with some period, there exists $d'_{4n_0+4}$ such that the pattern of zeros from $B_{4n_0+4}$ repeats on $\eta$ periodically, with period $d'_{4n_0+4}$. Thus, by taking $d_{4n_0+4}:=\lcm(d'_{4n_0+4},d_{4n_0+3})$, we obtain
$$
\eta[m_{4n_0+4}+\ell_{4n_0+4}+kd_{4n_0+4},m_{4n_0+4}+r_{4n_0+4}+kd_{4n_0+4}]= B_{4n_0+4}
$$
for each $k\in\Z$.

Finally, let $B_{4n_0+5}\in \{0,1\}^{[\ell_{4n_0+5},r_{4n_0+5}]}$, where $r_{4n_0+5}=r_{4n_0+4}$ (and $\ell_{4n_0+5}=r_{4n_0+5}-|B_{4n_0+5}|+1$), be the {\bf longest} block of the form $0\ldots 0B_{4n_0+4}$ that appears on $\eta$ and ends at position $m_{4n_0+4}+r_{4n_0+4}+k_0d_{4n_0+4}$ for some $k_0\in\Z$, i.e.\
$$
\eta[m_{4n_0+5}+\ell_{4n_0+5},m_{4n_0+5}+r_{4n_0+5}]=B_{4n_0+5},
$$
where $m_{4n_0+5}=m_{4n_0+4}+k_0d_{4n_0+4}$. Then, clearly $d_{4n_0+4}\divides m_{4n_0+5}-m_{4n_0+4}$. Moreover, since each zero on $\eta$ appears with some period, there exists $d'_{4n_0+5}$ such that the pattern of zeros from $B_{4n_0+5}$ repeats on $\eta$ periodically, with period $d'_{4n_0+5}$. Thus, by taking $d_{4n_0+5}:=\lcm(d'_{4n_0+5},d_{4n_0+4})$, we obtain
$$
\eta[m_{4n_0+5}+\ell_{4n_0+5}+kd_{4n_0+5},m_{4n_0+5}+r_{4n_0+5}+kd_{4n_0+5}]= B_{4n_0+5}
$$
for each $k\in\Z$.
\end{proof}
\cref{TTA} is an immediate consequence of \cref{PR51} and \cref{PR52}. Moreover, \cref{TTAwn} follows from \cref{TTA} and \cref{ausel1}. By \cref{TTAwn1}, $(S,X_\eta)$ is minimal if and only if it is Toeplitz. In fact, even $\eta$ may even happen to be a Toeplitz sequence:
\begin{Example}\label{ex54}
Let $\mathscr{B}:=\{b_i2^i : i\geq 1\}$, where $b_i\geq 2$ for $i\geq 1$. We will show that $\eta$ is a Toeplitz sequence. Indeed, for each $n\in\Z$ such that $\eta(n)=0$, there is $k_n\geq 1$ such that $\eta(n+jk_n)=0$ for all $j\in\Z$. Let now $n\in\Z$ be such that $\eta(n)=1$, i.e.
\begin{equation}\label{MM1}
n\not\equiv 0\bmod b_i{2^i} \text{ for each }i\geq 1.
\end{equation}
Let $m$ be odd, such that $n=m2^a$. We claim that
\begin{equation}\label{okres}
\eta(n+jb_1\ldots b_a2^{a+1})=1 \text{ for all }j\in\Z.
\end{equation}
Suppose not, so that for some $i_0$, we have
\begin{equation}\label{e2}
n+j_0b_1\ldots b_a2^{a+1}=K_0b_{i_0}2^{i_0}\text{ for some }j_0,K_0\in\Z.
\end{equation}
Then $i_0\leq a$; if not, by~\eqref{e2},  $2^{a+1}\divides n$ which is impossible. But now, again by~\eqref{e2}, $b_{i_0}2^{i_0} \divides n$ which contradicts~\eqref{MM1}.
\end{Example}
\begin{Remark}
Notice that, by \cref{B-skonczony}, it is easy to find $\sB$ such that $\eta$ is a Toeplitz sequence that is not periodic.
\end{Remark}
\begin{Remark}
Note that the Toeplitz sequence from \cref{ex54} is regular,\footnote{For the definition of a regular Toeplitz sequence, we refer the reader, e.g., to~\cite{MR2180227}.} so, in particular, $(S,X_\eta)$ is minimal and uniquely ergodic. To show the regularity of $\eta$, consider $d_n:=b_1\cdot\ldots\cdot b_n 2^{n+1}$. Consider two cases: $s\in\cf_\sB$, $s\in\cm_\sB$:
\begin{itemize}
\item
If $s\in\mathcal{M}_\mathscr{B}$ then $b_i2^i \divides s$ for some $i\geq 1$. If $i\leq n$ then $s+d_n\Z\subset \mathcal{M}_\mathscr{B}$. Otherwise, we have $2^{n+1}\divides s$. 
\item
If $s\in \cf_\sB$, let $m$ be odd, such that $s=m\cdot 2^a$. Then, by~\eqref{okres}, $s+b_1\cdot\ldots\cdot b_a2^{a+1}\Z\subset \cf_{\mathscr{B}}$. If $a\leq n$ then clearly $s+d_n\Z\subset \cf_\mathscr{B}$. Otherwise, we have $2^{n+1} \divides s$. 
\end{itemize}
It follows that if $s\in\Z$ satisfies
$$
(s+d_n\Z)\cap\cm_\mathscr{B} \neq\emptyset \text{ and }(s+d_n\Z)\cap\cf_\mathscr{B} \neq\emptyset,
$$
then $2^{n+1}\divides s$. The proportion of such $s$ in each integer interval of length $d_n$ equals $(b_1\cdot\ldots\cdot b_n)^{-1}$ and tends to zero as $n\to\infty$.
\end{Remark}

\subsection{Proximality}
We will now study the proximality of $(S,X_\eta)$. We will first show that for $\mathscr{B}$ pairwise coprime and infinite, $(S,\widetilde{X}_\eta)$ is proximal. This implies, by \cref{subsystem}, the proximality of $(S,X_\eta)$. By the same token, if $\widetilde{X}_{\eta'}\subset \widetilde{X}_{\eta}$ then $(S,\widetilde{X}_{\eta'})$ and $(S,{X}_{\eta'})$ are both proximal. Our aim (see \cref{AUR}) is to show that this is the only possible way to obtain a proximal $\mathscr{B}'$-free system $(S,X_{\eta'})$.

\subsubsection{Coprime case}
\begin{Prop}\label{pr:5}
If $\mathscr{B}\subset \N$ is infinite and coprime then $(S,\widetilde{X}_\eta)$ is syndetically proximal. In particular, $(S,X_\eta)$ is syndetically proximal.
\end{Prop}
\begin{proof}
By \cref{pr:4}, it suffices to show that
for any $x\in \widetilde{X}_\eta$ and $\vep>0$ the set
\begin{equation}\label{eq:2}
\{n\in\Z\colon d(S^nx,(\dots,0,0,0,\dots))<\vep\} \text{ is syndetic}.
\end{equation}
Fix $x\in X_\eta$. For $n\in \N$ and $k\geq 1$ there exists $m=m_{n,k}\in\Z$ such that
$$
x[n,\ldots,n+b_1\cdot\ldots \cdot b_k+k-1]\leq\eta[m,\ldots,m+b_1\cdot\ldots \cdot b_k+k-1].
$$
By the Chinese Remainder Theorem, there exists a unique $0\leq i_0\leq b_1\cdot\ldots\cdot b_k-1$ ($i_0=i_0(m,n)$) such that
$$
m+i_0+j\equiv 0\bmod b_{j+1}\text{ for }0\leq j\leq k-1,
$$
i.e.\ $x(n+i_0+j)\leq\eta(m_{n,k}+i_0+j)=0$ for $0\leq j\leq k-1$. This yields~\eqref{eq:2} and completes the proof.
\end{proof}
As an immediate consequence of \cref{pr:5} and \cref{pr:44}, we obtain the following:
\begin{Cor}\label{Krotriv}
For $\mathscr{B}\subset \N$ infinite and coprime, the maximal equicontinuous factor of $(S^{\times N},X_{\eta}^{\times N})$ is trivial for each $N\geq 1$.
\end{Cor}

\subsubsection{General case (proof of \cref{TTB})}\label{TTBsekcja}
\begin{Def}
We say that $\mathscr{B}\subset\N$ satisfies condition~\eqref{Au1}, whenever
\begin{equation}\label{Au1}
\text{there exists infinite pairwise coprime $\mathscr{B}'\subset\mathscr{B}$}.\tag{Au}
\end{equation}
\end{Def}
\begin{Def}
We say that $\mathscr{B}\subset \N$ satisfies condition~\eqref{Au}, whenever
\begin{multline}\label{Au}
\text{for any $k\in\N$ there exist }b^{(k)}_1,\ldots,b^{(k)}_k\in \mathscr{B}\text{ such that}\\
\gcd(b_i^{(k)},b_j^{(k)})\divides (j-i) \text{ for all }1\leq i<j\leq k.\tag{T$_{\text{prox}}$}
\end{multline}
\end{Def}
\begin{Th}\label{AUR}
Let $\sB\subset \N$. The following conditions are equivalent:
\begin{enumerate}[(a)]
\item\label{condA}
$(S,X_\mathscr{B})$ is proximal,
\item\label{condB}
$(S,\widetilde{X}_\eta)$ is proximal,
\item\label{condC}
$(S,X_\eta)$ is proximal,
\item\label{condD}
$(\ldots,0,0,0,\ldots)\in X_\eta$,\label{Af}
\item\label{condE}
$\mathscr{B}$ satisfies~\eqref{Au},
\item\label{zawi}
for any choice of $q_1,\dots,q_m>1$, $m\geq 1$, we have $\mathscr{B}\not\subset \bigcup_{i=1}^{m}\Z q_i$,
\item\label{condF}
$\mathscr{B}$ satisfies~\eqref{Au1},
\item $\cf_{\mathscr{B}}$ does not contain an infinite arithmetic progression.\label{condG}
\end{enumerate}
\end{Th}
Clearly, \cref{TTB} is an immediate consequence of \cref{AUR} and \cref{pr:4-}. Before we prove \cref{AUR}, we concentrate on its consequences.
\begin{Remark}\label{Auzawiera}
Clearly, if~\eqref{Au1} holds then $\eta \leq \eta'$, whence $X_\eta\subset \widetilde{X}_{\eta'}$.
\end{Remark}
By \cref{Auzawiera} and \cref{AUR}, we have the following:
\begin{Cor}\label{easy2}
If $(S,X_\eta)$ is proximal then $X_{\eta}\subset \widetilde{X}_{\eta'}$ with $\mathscr{B}'$ coprime.
\end{Cor}
\begin{Remark}\label{hp}
Recall (see~\cite{Ab-Le-Ru}) that if $\mathscr{B}$ is coprime and thin then $X_\eta=X_\mathscr{B}$. In particular, $X_\eta$ is hereditary.
\end{Remark}
By the implication \eqref{condD} $\Rightarrow$ \eqref{condC} in \cref{AUR}, we obtain the following:
\begin{Cor}\label{hertoprox}
If $X_\eta$ is hereditary then $(S,X_\eta)$ is proximal.
\end{Cor}

\begin{Question}[Cf.\ Example~\ref{ex:2.4a}]\label{q2}
Is it possible that $X_\eta\subsetneq \widetilde{X}_\eta=X_\mathscr{B}$ with $X_\eta$ proximal?
\end{Question}

The proof of \cref{AUR} will be divided into several observations.
\begin{Remark}
Since $X_\eta\subset \widetilde{X}_\eta\subset X_\mathscr{B}$, by \cref{subsystem}, we have $\eqref{condA} \Rightarrow \eqref{condB} \Rightarrow \eqref{condC}$.
\end{Remark}
\begin{Lemma}
We have \eqref{condC} $\Rightarrow$ \eqref{condD}.
\end{Lemma}
\begin{proof}
If $(S,X_\eta)$ is proximal then, by \cref{pr:4-}, it has a fixed point, i.e.\ either $(\dots,0,0,0,\dots)\in X_\eta$ or $(\dots,1,1,1,\dots)\in X_\eta$. The latter of the two is impossible, since each zero on $\eta$ appears on $\eta$ with bounded gaps and the claim follows.
\end{proof}
\begin{Lemma}
We have \eqref{condD} $\Rightarrow$ \eqref{condE}.
\end{Lemma}
\begin{proof}
If $(\ldots,0,0,0,\ldots)\in X_\eta$ then there are arbitrarily long blocks of consecutive zeros on $\eta$. In other words, given $k\geq1$, we can solve the systems of congruences:
$$
i_0+i-1 \equiv 0\bmod b_{s_i},\ 1\leq i\leq k.
$$
Suppose that $d\divides (b_{s_i},b_{s_j})$. Then $d\divides i_0+i-1$ and $d\divides i_0+j-1$, whence $d\divides (j-i)$. This completes the proof.
\end{proof}
\begin{Lemma}
We have \eqref{condE} $\Rightarrow$ \eqref{zawi}.
\end{Lemma}
\begin{proof}
Suppose that \eqref{condE} holds. Without loss of generality, we can assume that $\{q_1,\dots,q_m\}$ is coprime (indeed, we can always find a coprime set $\{q_1',\dots,q_n'\}$ such that $\bigcup_{i=1}^mq_i\Z\subset \bigcup_{i=1}^{n}q_i'\Z$). Let $k\ge q_1\ldots q_m$ and choose  $b^{(k)}_1,\ldots,b^{(k)}_k\in \mathscr{B}$ satisfying condition~\eqref{Au}. For $i=1,\ldots,m$, let
$$
M_i:=\{1\leq \ell\leq k:  b^{(k)}_\ell\in q_i\Z \}.
$$
Then, by~\eqref{Au},   $q_i \divides (\ell+\ell')$ for any $\ell,\ell'\in M_i$, whence 
\begin{equation}\label{VV1}
M_i\subset q_i\Z +r_i \text{ for some }r_i.
\end{equation}
For $i=1,\ldots,m$, choose a natural number $r'_i$ such that $q_i\ndivides (r_i-r_i')$. By the Chinese Remainder Theorem there exists a natural number $j\le q_1\ldots q_m$ (note that $j\leq k$) such that $j\equiv r'_i\bmod q_i$ for $i=1,\ldots,m$. By \eqref{VV1}, it follows that $j\notin M_i$ for any $i=1,\ldots,m$. It follows that $b^{(k)}_j\notin q_1\Z \cup\ldots\cup q_m\Z$.
\end{proof}
\begin{Lemma}
We have \eqref{zawi} $\Rightarrow$ \eqref{condF}.
\end{Lemma}
\begin{proof}
We will proceed inductively. Fix $c_1\in \mathscr{B}$. Suppose that for $k\ge 1$ we have found pairwise coprime subset $\{c_1,\ldots,c_k\}\subset \mathscr{B}$. Let $\{q_1,\ldots,q_m\}$ be the set of all prime divisors of $c_1,\ldots, c_k$. Then any $c_{k+1}\in \mathscr{B}\setminus (q_1\Z \cup\ldots\cup q_m\Z)$ is coprime with each of $c_1,\ldots,c_{k}$.
\end{proof}
\begin{Remark}
If~\eqref{condF} holds then, by \cref{Auzawiera}, we have $X_\eta\subset \widetilde{X}_{\eta'}$. By \cref{pr:5}, $\widetilde{X}_{\eta'}$ is proximal. Hence, by \cref{subsystem}, we obtain $\eqref{condF} \Rightarrow \eqref{condA}$.
\end{Remark}
\begin{Remark}
Condition \eqref{condD} implies that $\cm_{\mathscr{B}}$ contains intervals of integers of arbitrary length. Hence \eqref{condD} $\Rightarrow$ \eqref{condG}.
\end{Remark}
\begin{Lemma}
We have \eqref{condG} $\Rightarrow$ \eqref{zawi}.
\end{Lemma}
\begin{proof}
Suppose that \eqref{zawi} does not hold and let $q_1,\dots,q_k$, $k\geq 1$, be such that $\mathscr{B}\subset\bigcup_{i=1}^k \Z q_i$. Let $M:=q_1\cdot\ldots\cdot q_k$. We claim that
$$
b\ndivides \ell M+1 \text{ for every }b\in\mathscr{B},\text{ i.e.\ }\ell M+1\in \mathcal{F}_\mathscr{B}.
$$
Indeed, given $b\in\mathscr{B}$, there exists $q_i$ ($1\leq i\leq k$) such that $q_i\divides b$. If $b\divides \ell M+1$ then $q_i \divides \ell M+1$. This is however impossible since $q_i\divides M$.
\end{proof}
The proof of \cref{AUR} is complete in view of the above remarks and lemmas.
%%%%%%%%%%%%%%%%%%%%%%%%%%%%%%%%%%%%%%%%%%%%%%%%%\section{Proximality}\label{s5}
%%%%%%%%%%%%%%%%%%%%%%%%%%%%%%%%%%%%%%%%%%%%%%%%%%%%%%%%%%%%%%%%%%%%%%%%%%%
We will give now one more characterization of proximal $(S,X_\eta)$, in terms of the maximal equicontinuous factor (cf.\ \cref{Krotriv}):
\begin{Th}\label{krone}
$(S,X_\eta)$ is proximal if and only if its maximal equicontinuous factor is trivial.
\end{Th}
For the proof, we will need the following lemma:
\begin{Lemma}\label{dd'}
Let $d\geq 1$ and let $A\subset \{0,1,\dots,d-1\}$. Suppose that for any $k\geq 1$ there exist $n_k\in\Z$ and $0\leq r_k\leq d-1$ such that
\begin{equation}\label{fpap}
A+md+r_k\subset \cf_\mathscr{B} \text{ for }n_k\leq m\leq n_k+k.
\end{equation}
Then, for any $0\leq r\leq d-1$ such that there are infinitely many $k\geq 1$ with $r_k=r$, we have
\begin{equation}\label{ddclaim}
A+\Z d+r\subset \cf_\mathscr{B}.
\end{equation}
\end{Lemma}
\begin{proof}
Let $0\leq r\leq d-1$ be such that there are infinitely many $k\geq 1$ satisfying~\eqref{fpap} with $r_k=r$, i.e.
\begin{equation}\label{pom:1}
A+md+r\subset \cf_\mathscr{B} \text{ for }n_k\leq m\leq n_k+k.
\end{equation}
Suppose that~\eqref{ddclaim} fails. Then, for some $a\in A$ and $k\in \Z$, we have $a+kd+r\in\cm_\mathscr{B}$. In other words, for some $b\in\mathscr{B}$, we have $b\divides a+kd+r$. It follows that for any $\ell\in\Z$
$$
b\divides a+(k+\ell b)d+r.
$$
This, however, contradicts~\eqref{pom:1}.
\end{proof}
\begin{proof}[Proof of \cref{krone}]
Since proximality implies that the maximal equicontinuous factor is trivial, we only need to show the converse implication. Suppose that $(S,X_\eta)$ is not proximal. Let $d\geq 1$ be the smallest number such that $\cf_\mathscr{B}$ contains an infinite arithmetic progression with difference $d$ (such $d$ exists by \cref{AUR} \eqref{condG}). Let $F\subset \{0,\dots,d-1\}$ be the maximal set such that
\begin{equation}\label{zal}
F+\Z d \subset \cf_\mathscr{B}
\end{equation}
($F\neq\emptyset$ by the definition of $d$). We claim that for any $y\in X_\eta$, there exists a unique $0\leq r<d$ such that
\begin{equation}\label{withr}
y(a+md+r)=1 \text{ for all }a\in F \text{ and }m\in\Z.
\end{equation}
Since $y\in X_\eta$, it follows by~\eqref{zal} that such $r$ exists and we only need to show uniqueness. Suppose that~\eqref{withr} holds for $r=r_1,r_2$, where $d\ndivides (r_1-r_2)$, i.e., we have
$$
y(a+md)=1 \text{ for all }a\in (F+r_1)\cup (F+r_2) \text{ and }m\in\Z.
$$
Since $y\in X_\eta$, each block from $y$ appears on $\eta$ and it follows that the assumptions of Lemma~\ref{dd'} hold for $A:=(F+r_1)\cup (F+r_2) \bmod d$. Therefore, using additionally~\eqref{zal},
$$
[F\cup (F+r_1+s)\cup (F+r_2+s)]+\Z d\subset \cf_\mathscr{B}\text{ for some }s.
$$
Note that by the minimality of $d$, we have $F+i \neq F \bmod d$ for $0<i<d$. Therefore, 
$$
F\subsetneq F\cup (F+r_1+s)\cup (F+r_2+s).
$$
This contradicts the maximality of $F$ and thus indeed implies the uniqueness of~$r$. It follows that
$$
X_\eta=\bigcup_{i=0}^{d-1}X_\eta^{(i)},\ X_\eta^{(i)}=\{y\in X_\eta : \eqref{withr}\text{ holds for }r=i \}
$$
is a decomposition of $X_\eta$ into $d$ pairwise disjoint sets. Clearly, each $X_\eta^{(i)}$ is closed and $S X_\eta^{(i)}=X_\eta^{(i-1)}$, where $X_\eta^{(-1)}=X_\eta^{(d-1)}$. It follows that $(S,X_\eta)$ has the (minimal) rotation on $d$ points as a topological factor, which completes the proof.
\end{proof}

The following natural question arises:
\begin{Question}\label{Q323}
Given $\sB\subset \N$, what is the maximal equicontinuous factor of $(S,X_\eta)$?
\end{Question}
We provide below the answer to \cref{Q323} in the simplest case of finite sets $\mathscr{B}$, where $(S,X_\eta)$ turns out to be equicontinuous. Moreover, we will show that if $X_\eta=X_\eta\cap Y$ then $(T,G)$ defined as in \cref{canonical}, is the maximal equicontinuous factor of $(S,X_\eta)$.

We will need the following well-known fact:
\begin{equation}\label{liniowe}
  \parbox{0.8\linewidth}{Let $m,a,b\in\N$. The equation $ax\equiv b\bmod m$ has a solution in $x\in\Z$ if and only if $\gcd(m,a)\divides b$.}
\end{equation}

\begin{Prop}\label{B-skonczony}
Let $\mathscr{B}\subset \N$. Then $\mathscr{B}$ is finite if and only if $\eta$ is periodic, with the minimal period $m=\lcm(\mathscr{B})$.\footnote{Recall that we assume that $\mathscr{B}$ is primitive.}
\end{Prop}
\begin{proof}
If $\mathscr{B}$ is finite then $\eta$ is periodic with period $\lcm \mathscr{B}$. Suppose now that $\eta$ is periodic and denote its period by $m$. Let $1\leq r_1 < r_2 < \ldots < r_s\leq m$ be such that $(\text{supp }\eta) \cap [1,m]=\{1,\dots,m\}\setminus \{r_1,\dots,r_s\}$. Then
\[
\bigcup_{b\in \mathscr{B}} b\Z=\bigcup_{\ell=1}^s (m\Z+r_\ell).
\]
For $1\leqslant \ell \leqslant s$, let $d_\ell:=\gcd(m,r_\ell)$. By the definition of $d_\ell$,
\begin{equation}\label{VV2}
d_\ell\Z\supset m\Z+r_\ell
\end{equation}
Then, by~\eqref{liniowe}, there exists $k_l\in\Z$ such that $r_\ell k_\ell\equiv d_\ell \mod m$. Since $\eta(r_\ell)=0$, we have $\eta(r_\ell k_\ell)=0$, which, by periodicity, yields $\eta(d_\ell)=0$. This and~\eqref{VV2} imply
\begin{equation}\label{dwojkapon}
\bigcup\limits_{b\in \mathscr{B}} b\Z=\bigcup\limits_{\ell=1}^s d_{\ell} \Z.
\end{equation}
Fix $b\in\mathscr{B}$. It follows from~\eqref{dwojkapon} that $d_\ell \divides b$ for some $1\leq \ell\leq s$. On the other hand, there exists $b'\in\mathscr{B}$ such that $b' \divides d_\ell$. By the primitivity of $\mathscr{B}$, we have $b'\divides b$, whence $b=b'$ and $d_\ell=b$. We conclude that $\mathscr{B}\subset \{d_\ell : 1\leq \ell\leq s\}$, i.e.\ $\mathscr{B}$ is finite. Moreover, since $d_\ell\divides m$ for $1\leq \ell\leq s$, we obtain $b\divides m$ for each $b\in\mathscr{B}$. This yields $\lcm(\mathscr{B})\divides m$.
\end{proof}
As an immediate consequence of \cref{B-skonczony}, we have:
\begin{Cor}
If $\mathscr{B}$ is finite then $(S,X_\eta)$ is finite whence equicontinuous.
\end{Cor}

\begin{Prop}\label{popopop}
Suppose that $X_\eta=X_\eta\cap Y$. Then $(T,G)$ is the maximal equicontinuous factor of $(S,X_\eta)$. In particular, if we additionally assume that $\sB$ is infinite then the maximal equicontinuous factor of $(S,X_\eta)$ is infinite.
\end{Prop}
\begin{proof}
Notice first that, by \cref{cgtheta}, $\theta \colon X_\eta\to G$ is well-defined and continuous. Thus, $(T,G)$ is an equicontinuous factor of $(S,X_\eta)$ and we only need to show its maximality. 
Notice that the (discrete) spectrum of the maximal equicontinuous factor of $(S,X_\eta)$ is alwyas included in the discrete part of the spectrum of $(S,X_\eta,\nu)$ for any $\nu\in\mathcal{P}(S,X_\eta)$. Therefore, to prove the maximality of $(T,G)$, it suffices to find $\nu$ such that the discrete part of the spectrum of $(S,X_\eta,\nu)$ agrees with the (discrete) spectrum of $(T,G,\PP)$. We have
$$
(T,G,\PP)\xrightarrow{\varphi} (S,X_\eta,\nu_\eta) \xrightarrow{\theta} (T,G,\PP).
$$
It follows by the coalescence of $(T,G,\PP)$ that $\varphi$ yields an isomorphism of $(T,G,\PP)$ and $(S,X_\eta,\nu_\eta)$. In particular, the (discrete) spectrum of $(T,G,\PP)$ is the same as the (discrete) spectrum of $(S,X_\eta,\nu_\eta)$ and the claim follows.
\end{proof}

\begin{Example}
Let $\mathscr{B}$ be as in \cref{ex54}. Then $\sum_{i\geq 1}\frac{1}{2^ib_i}\leq \sum_{i\geq 1}\frac{1}{2^i}$ is thin and it follows by \eqref{LLtauty1}, \cref{LLtauty1a} and by \cref{gdzieeta} that $\eta\in Y$. Moreover, by the minimality of $(S,X_\eta)$, for each $0\leq s_k\leq b_k$, $k\geq 1$, we have that either $X_\eta \cap Y^k_{\geq s_k}=X_\eta$ or $X_\eta \cap Y^k_{\geq s_k}=\emptyset$. Since $\eta\in Y$, it follows that $X_\eta \cap Y^k_{\geq s_k}=\emptyset$ whenever $s_k\geq 2$. Since $X_\eta=X_\eta \cap (\bigcup_{1\leq s_k\leq b_k}Y^k_{s_k})$ for each $k\geq 1$, it follows that $X_\eta=X_\eta\cap Y$. By \cref{popopop}, the associated canonical odometer $(T,G)$ is the maximal equicontinuous factor of $(S,X_\eta)$.
\end{Example}

%%%%%%%%%%%%%%%%%%%%%%%%%%%%%%%%%%%%%%%%%%%%%%%%%%%%%%%%%%%%%%%%%%%%%%%%%%%
%%%%%%%%%%%%%%%%%%%%%%%%%%%%%%%%%%%%%%%%%%%%%%%%%%%%%%%%%%%%%%%%%%%%%%%%%%%
%%%%%%%%%%%%%%%%%%%%%%%%%%%%%%%%%%%%%%%%%%%%%%%%%%%%%%%%%%%%%%%%%%%%%%%%%%%

\subsection{Transitivity}

\subsubsection{Transitivity of $(S,\widetilde{X}_\eta)$ and $(S,X_\mathscr{B})$}

\begin{Prop}\label{tanzyt}
For any $\mathscr{B}\subset \N$ such that the support of $\eta$ is infinite, the following conditions are equivalent:
\begin{enumerate}[(a)]
\item\label{tanzytA}
$(S,\widetilde{X}_\eta)$ is transitive.
\item\label{tanzytB}
$(S,\widetilde{X}_\eta)$ does not have open wandering sets of positive diameter.
\item\label{tanzytC}
For any block $B$ that appears on $\eta$ there exists a block $B'\geq B$ (coordinatewise) that appears on $\eta$ infinitely often.
\end{enumerate}
\end{Prop}
The implication  \eqref{tanzytA} $\Rightarrow$ \eqref{tanzytB} from \cref{tanzyt} is a consequence of the following general lemma:
\begin{Lemma}\label{tanzytlemma}
Let $(T,X)$ be a topological dynamical system with a transitive point $x\in X$. Then  $(T,X)$ has no open wandering sets of positive diameter.
\end{Lemma}
\begin{proof}
Let $U$ be an open wandering set for $(T,X)$. Then the orbit of $x$ visits $U$ exactly once. It follows that $U$ must be a singleton.
\end{proof}
\begin{proof}[Proof of \cref{tanzyt}]
In view of \cref{tanzytlemma}, it remains to show \eqref{tanzytB} $\Rightarrow$ \eqref{tanzytC} $\Rightarrow$ \eqref{tanzytA}. We will prove first \eqref{tanzytB} $\Rightarrow$ \eqref{tanzytC}. Suppose that \eqref{tanzytC} does not hold. Let $B$ be a block on $\eta$ such that all blocks $B'\geq B$ appear on $\eta$ (at most) finitely many times. Let
\begin{align*}
K&:=\min\{k\in \Z : \eta[k,k+|B|-1]\geq B\},\\
L&:=\max\{k+|B|-1 : \eta[k,k+|B|-1]\geq B\}
\end{align*}
(in particular, blocks $B'\geq B$ do not appear on $\eta$ outside $\eta[K,L]$). We claim that, for any $x\in\widetilde{X}_\eta$, the block $C:=\eta[K,L]$ appears on $x$ at most once. Suppose that, for some $x\in\widetilde{X}_\eta$, $C$ appears on $x$ twice. It follows that a block of the form $C' D C''$, where $C',C''\geq C$, appears on $\eta$ and this is impossible by the choice of $C$. Thus, the cylinder set 
$$
\mathcal{C}:=\{x\in\widetilde{X}_\eta : x[K,L]=C\}
$$
corresponding to $C$ is an open wandering set. Clearly, we have $\eta \in \mathcal{C}$. Moreover, since the support of $\eta$ is infinite, we also have $x\in\mathcal{C}$ for $x$ given by
$x(n)=\eta(n)$ for $n\in[K,L]$, $x(n)=0$ otherwise. It follows that $|\mathcal{C}|\geq 2$, i.e.\ the diameter of $\mathcal{C}$ is positive and we conclude that \eqref{tanzytB} fails.

We will now prove \eqref{tanzytC} $\Rightarrow$ \eqref{tanzytA}. By \cref{tranzyrem}, given blocks $B',C'$ that appear on $\eta$ and $B\leq B', C\leq C'$, it suffices to show that there exists $x\in \widetilde{X}_\eta$ such that both $B$ and $C$ appear on $x$. It follows by \eqref{tanzytC} that there exists $B''\geq B'$ that appears on $\eta$ infinitely often. Therefore for some block $D$, a block of the form $C' D B''$ or a block of the form $B'' D C'$ appears on $\eta$. Hence, $x:=(\dots,0,0, B ,\underbrace{0,\dots,0}_{|D|}, C,0,0,\dots)\in \widetilde{X}_\eta$ and the result follows.
\end{proof}
As an immediate consequence of \cref{tanzyt}, we obtain the following:
\begin{Cor}\label{Coroloro}
Let $\mathscr{B}\subset \N$ be such that $\eta$ is recurrent. Then $(S,\widetilde{X}_\eta)$ is transitive.
\end{Cor}
In particular, by \cref{Coroloro},  \cref{quasi-gen} and \cref{OOF}, we have the following:
\begin{Cor}
The subshift $(S,\widetilde{X}_\eta)$ is transitive whenever $\mathscr{B}$ has light tails.
\end{Cor}
\begin{Remark}
Let $\sB$ be as in \cref{ex:2.4}, i.e.\ $\sB=\mathscr{P}$. Then $(S,\widetilde{X}_\eta)$ fails to be transitive.
\end{Remark}

Clearly, if $X_\eta=X_\mathscr{B}$ then $(S,X_\mathscr{B})$ is transitive.\footnote{Recall that $X_\eta=X_{\sB}$ holds for $\sB$ satisfying~\eqref{settingerdosa}.} We will now give an example, where $(S,X_\mathscr{B})$ fails to be transitive.
\begin{Example}
Let $\mathscr{B}$ be as in \cref{X_eta<>X_B}, i.e.\ $4,6\in\mathscr{B}$ and $b>12$ for $b\in\mathscr{B}\setminus \{4,6\}$. Let
\begin{align*}
A_1&:=110011100110,\\
A_2&:=011101010111=\eta[0,11].
\end{align*}
Suppose that both $A_1,A_2$ appear on $x\in\{0,1\}^\Z$. We will show that $x\not\in X_\mathscr{B}$. Indeed, we have
\begin{align*}
\Z/4\Z \setminus (\text{supp }A_1 \bmod 4)&=\{3\},\\
\Z/4\Z \setminus (\text{supp }A_2 \bmod 4)&=\{0\}.
\end{align*}
Let $k,\ell\in\Z$ be such that $x[k,k+11]=A_1$ and $x[\ell,\ell+11]=A_2$. It follows that if $x$ is $\{4\}$-admissible then $4\divides k+3+\ell$. In a similar way, if $x$ is $\{6\}$-admissible then $6\divides k+2+\ell$. Since one of the numbers $k+3+\ell$ and $k+2+\ell$ is odd, we conclude that $x$ is not $\{4,6\}$-admissible, so all the more, it is not $\mathscr{B}$-admissible.
\end{Example}

\subsubsection{$(S\times S,X_{\eta}\times X_{\eta})$ is not transitive}
Our main goal in this section is to show that $(S\times S,X_{\eta}\times X_{\eta})$ is not transitive. As a consequence, we will have the following whenever $(S,X_\eta)$ is proximal:
\begin{itemize}
\item
$(S,X_\eta)$ is transitive with trivial maximal equicontinuous factor,
\item
$(S\times S,X_{\eta}\times X_{\eta})$ has trivial equicontinuous factor, but it is not transitive.
\end{itemize}
Analogous phenomenon is impossible in ergodic theory. Our main tool is the following result:
\begin{Prop}\label{sarnak10}
$(S,X_{\eta})$ has a non-trivial topological joining with $(T, G )$.
\end{Prop}
\begin{proof}
Let
$$
N:=\overline{\mathcal{O}_{T\times S}(\un{0},\eta)},
$$
where $\underline{0}=(0,0,\dots)$, i.e.\ $N$ is the closure of the graph of $\varphi$ along the orbit of $\underline{0}$ (indeed, we have $S^n\eta=S^n\varphi(\underline{0})=\varphi(T^n\underline{0})$). Since the orbit of $\underline{0}$ under $T$ is dense in $G$ and the orbit of $\eta$ under $S$ is dense in $X_\eta$, it follows that $N$ has full projection on both coordinates. Moreover, $N$ is closed and $T\times S$-invariant. It remains to show that $N\neq G\times X_\eta$. Take $(\dots,0,0,0,\dots)\neq x\in X_\eta$. We claim that $\{g\in G : (g,x)\in N\}\neq G$. Indeed, let $k_0\in\Z$ be such that $x(k_0)=1$ and suppose that $(T^{n_i}\times S^{n_i})(\underline{0},\eta)\to (g,x)$. Then $S^{n_i}\eta\to x$, whence, for $i$ sufficiently large, $\eta(k_0+n_i)=S^{n_i}\eta(k_0)=x(k_0)=1$. It follows that $n_i+k_0 \in \mathcal{F}_\mathscr{B}$, i.e.\ $n_i+k_0 \neq 0\bmod b_k$ for each $k\geq 1$. On the other hand, we have $T^{n_i}\underline{0}\to g$, i.e.\ $(n_i,n_i,\dots)\to (g_1,g_2,\dots)$. Thus, $g_k \neq -k_0 \bmod b_k$ for each $k\geq1$. Hence, $\{g\in G : (g,x)\in N\}\neq G$ for $x\neq (\dots,0,0,0,\dots)$, which completes the proof.
\end{proof}
\begin{Remark}
Suppose that $\mathscr{B}$ is taut and $1\not\in\mathscr{B}$. By \cref{gdzieeta}, $\eta\in Y$, i.e.\ for each $k\geq 1$, we have $\mathcal{F}_\mathscr{B} \bmod b_k = (\Z/b_k\Z)\setminus \{0\}$. It follows by the above proof that $\{g\in G : (g,\eta)\in N\}=\{\underline{0}\}$. In a similar way, if $x\in Y$ then $\{g\in G : (g,x)\in N\}$ is a singleton, in particular, for each $n\in\Z$, the set $\{g\in G : (g,S^n\eta)\in N\}$ is a singleton.
\end{Remark}

\begin{Cor}\label{sarnak11} $(S\times S,X_{\eta}\times X_{\eta})$ is not transitive.
\end{Cor}
\begin{proof}
In view of \cref{sarnak10}, we can use the theorem about disjointness of topologically weakly mixing systems with (minimal) equicontinuous systems (see Thm.\  II.3 in \cite{MR0213508}).
\end{proof}

%%%%%%%%%%%%%%%%%%%%%%%%%%%%%%%%%%%%%%%%%%%%%%%%

\section{Tautness}\label{s3}
\subsection{$\eta$ is quasi-generic for $\nu_\eta$ (proof of \cref{OOE})}
\begin{Th}\label{quasi-gen}
Given $\mathscr{B}\subset \N$, let $(N_k)$ be such that
$$
\underline{d}(\mathcal{M}_\mathscr{B})=\lim_{k\to \infty}\frac{1}{N_k}|[1,N_k]\cap\mathcal{M}_\mathscr{B}|.
$$
Then $\eta$ is quasi-generic for $\nu_\eta$ along $(N_k)$. In particular, if $\mathscr{B}$ is Besicovitch then $\eta$ is generic for $\nu_\eta$.
\end{Th}
\begin{proof}
To simplify the notation, we will only deal with the (most involved) case when $\mathscr{B}$ is infinite. According to \cite{Ab-Le-Ru}, by a pure measure theory argument, we only need to prove that
$$
\frac{1}{N_k}\sum_{n\leq N_k}\raz_{\varphi^{-1}(A)}(T^n\underline{0})\to \PP(\varphi^{-1}(A))
$$
for each
$A=\{x\in \{0,1\}^\Z : x(j_s)=0, s=1,\dots, r\},\ j_1<\dots<j_r, r\geq 1$.
Recall that
$$
C=\{( g _1, g _2,\ldots)\in{ G }: g _k\not\equiv 0\bmod b_k \text{ for }k\geq 1\}
$$
and, for $K\geq 1$, define
$$
C_K:=\{( g _1, g _2,\ldots)\in{ G }: g _k\not\equiv 0\bmod b_k \text{ for }1\leq k\leq K\}.
$$
Then each $C_K$ is clopen and $C_K \searrow C$ when $K\to \infty$.
We have
$$
\varphi^{-1}(A)=\bigcap_{s=1}^{r} T^{-j_s}C^c,
$$
whence
\begin{equation}\label{KK1}
\bigcap_{s=1}^{r} T^{-j_s}C_K^c\subset \varphi^{-1}(A) \subset \bigcap_{s=1}^{r} T^{-j_s}C_K^c \cup \bigcup_{s=1}^{r}T^{-j_s}(C^c\setminus C_K^c).
\end{equation}
Moreover, since $\raz_{\bigcap_{s=1}^rT^{-j_s}C_K^c}$ is continuous, by the unique ergodicity of $T$ in \cref{monoerg}, we have
\begin{equation}\label{KK2}
\frac{1}{N_k}\sum_{n\leq N_k}\raz_{\bigcap_{s=1}^{r}T^{-j_s}C_K^c}(T^n\underline{0})\to \PP(\bigcap_{s=1}^{r}T^{-j_s}C_K^c)
\end{equation}
and, given $\vep>0$, for $K$ sufficiently large, we have
\begin{equation}\label{KK3}
\PP(\bigcap_{s=1}^{r}T^{-j_s}C_K^c) \geq \PP(\bigcap_{s=1}^{r}T^{-j_s}C^c)-\vep.
\end{equation}
Notice that
$$
T^n\underline{0} \in C^c\setminus C_K^c \iff n\in \mathcal{M}_{\mathscr{B}}\setminus \mathcal{M}_{\{b_1,\dots, b_K\}}.
$$
By \cref{da-er}, if $K$ is large enough then
$$
d(\mathcal{M}_{\{b_1,\dots, b_K\}})\geq \underline{d}(\mathcal{M}_\mathscr{B})-\vep.
$$
Therefore, and by the choice of $(N_k)$,
\begin{equation}\label{KK4}
\limsup_{k\to \infty}\frac{1}{N_k}\sum_{n\leq N_k}\raz_{\bigcup_{s=1}^{r}(C^c\setminus C_K^c)}(T^n\underline{0})\leq \vep.
\end{equation}
Putting together~\eqref{KK1},~\eqref{KK2},~\eqref{KK3} and~\eqref{KK4} completes the proof.
\end{proof}
\begin{Remark}\label{kiedybeh}
Notice that by \cref{quasi-gen}, we have
$$
\nu_\eta(\{x\in\{0,1\}^\Z : x(0)=1\})=\lim_{k\to\infty}\frac{1}{N_k}|\{1\leq n\leq N_k : \eta(n)=1\}|=\ov{d}(\cf_\sB).
$$
It follows immediately that $\sB$ is Behrend if and only if $\nu_\eta=\delta_{(\dots,0,0,0,\dots)}$.
\end{Remark}

By \cref{quasi-gen}, if a block does not appear on $\eta$ then the Mirsky measure of the corresponding cylinder set is zero. As a consequence, we obtain:
\begin{Cor}\label{dzc1}
$\nu_\eta(X_\eta)=1$.
\end{Cor}
In view of \cref{quasi-gen} and \cref{dzc1}, \cref{OOE} has been proved.

\begin{Remark}\label{tautnotdelta}
As an immediate consequence of \cref{quasi-gen}, we have 
$$
\ov{d}(\cf_\mathscr{B})>0 \iff \nu_\eta\neq\delta_{(\dots,0,0,0,\dots)}.
$$
In particular, it follows by~\eqref{tautnotbeh} that $\nu_\eta\neq\delta_{(\dots,0,0,0,\dots)}$ whenever $\mathscr{B}\neq \{1\}$ is taut.
\end{Remark}

\subsection{Tautness and Mirsky measures (\cref{TTC} -- first steps)}\label{tautyuberall}
In this section our main goal is to prove the following:

\begin{Th}\label{tautywszystko}
For each $\mathscr{B}\subset \N$, there exists a taut set $\mathscr{B}'\subset \N$, such that $\cf_{\mathscr{B}'}\subset \cf_{\sB}$ and $\nu_\eta=\nu_{\eta'}$.\footnote{We will see later that, in fact, the equality $\nu_\eta=\nu_{\eta'}$ determines $\mathscr{B}'$, cf.\ \cref{zawiewnio1}.}
\end{Th}
In course of the construction of $\mathscr{B}'$ and to prove that $\mathscr{B}'$ satisfies the required properties, we will use the following general lemmas (they are easy consequences of \cref{beh2} and \cref{beh3}):
\begin{Lemma}\label{ZZ1}
Suppose that $\mathscr{B}\subset \N$ is primitive. Then $\mathscr{B}$ is taut if and only if $\sB$ there exists a cofinite subset of $\mathscr{B}$ that is taut.
\end{Lemma}
\begin{proof}
Let $\sB\subset \N$ be primitive. It suffices to show that if $\mathscr{B}\setminus \{b\}$ is taut for some $b\in\sB$ then $\sB$ is taut. Suppose that $\sB$ fails to be taut. By \cref{beh3}, there exist $c\in\N$ and a Behrend set $\sA$ such that $c\sA\subset \sB$. Then $c\sA' \subset \sB\setminus \{b\}$, where $\sA'=\sA\setminus \{b/c\}\subset$ and $\sA'$ is Behrend by \cref{beh2}. Applying again \cref{beh3}, we conclude that $\mathscr{B}\setminus \{b\}$ also fails to be taut.
\end{proof}
\begin{Lemma}\label{Mnajmn}
Suppose that $\mathscr{B}\subset \N$ is primitive. If $\sB$ is not taut then, for some $c\in\N$, the set 
\begin{equation}\label{jakA}
\sA_c:=\left\{ \frac{b}{c} : b\in\mathscr{B} \text{ and }c\divides b\right\}
\end{equation}
is Behrend. 
\end{Lemma}
\begin{proof}
Clearly, for any $c\in \N$, we have $c \sA_c \subset \sB$, where $\sA_c$ (possibly empty) is as in \eqref{jakA}.
By \cref{beh3}, we have
$$
C:=\{c\in\N : c\sA'_c\subset \sB \text{ for some Behrend set }\sA'_c\}\neq\emptyset
$$
and, for any $c\in C$, we have $\sA'_c\subset \sA_c$, whence $\sA_c$ is Behrend. This completes the proof.
\end{proof}
\begin{Lemma}\label{rozbicie}
Let $\sB_1,\sB_2\subset \N$ be disjoint and such that $\sB:=\sB_1\cup \sB_2$ is primitive. Then $\sB$ is taut if and only if both $\sB_1$ and $\sB_2$ are taut.
\end{Lemma}
\begin{proof}
If $\sB_i$ is not taut for some $i\in\{1,2\}$ then, by \cref{beh3}, there exist $c\in\N$ and a Behrend set $\sA$ such that $c\sA\subset \sB_i\subset \sB$. Applying again \cref{beh3}, we deduce that $\sB$ also fails to be taut.  On the other hand, if $\sB$ is not taut then, by \cref{beh3}, there exist $c\in\N$ and a Behrend set $\sA$ such that $c\sA\subset \sB$. Let
$$
\sA_i:=\left\{\frac{b}{c} : b\in \sB_i\right\}, i=1,2.
$$
Clearly, $\sA=\sA_1\cup \sA_2$. Moreover, by \cref{beh2}, $\sA_i$ is Behrend for some $i\in\{1,2\}$. We obtain $c\sA_i\subset \sB_i$ for this $i$ and, by \cref{beh3}, we conclude that $\sB_i$ fails to be taut.
\end{proof}

\paragraph{Construction.}\label{cons}
We may assume without loss of generality that $\sB$ is primitive (cf.\ \cref{stopkaprimitive}).
\paragraph{Step 0.}
If $1\in \mathscr{B}$, we set $\mathscr{B}':=\{1\}$. 
%%%%%%%%%%%%%%%%%%%%%%%%%%%
\paragraph{Step 1.}
Suppose now that $1\not\in\mathscr{B}$ and suppose that $\mathscr{B}$ is not taut. Let $c_1\in\N$ be the smallest natural number such that
$$
\mathscr{A}^1:=\left\{\frac{b}{c_1} : b\in\mathscr{B} \text{ and }c_1 \divides b\right\}
$$
is Behrend (such $c_1$ exists by \cref{Mnajmn}). By the definition of $\sA^1$, we have $\sB\setminus c_1\sA^1=\sB\setminus c_1\Z$. Let
\begin{equation}\label{STEP1}
\mathscr{B}^1:=(\mathscr{B}\setminus c_1\Z) \cup \{c_1\}=(\mathscr{B}\setminus c_1\mathscr{A}^1)\cup\{c_1\}.
\end{equation}
We claim that $\mathscr{B}^1$ is primitive. Indeed, if this is not the case then, by the primitivity of $\sB$, for some $b\in \mathscr{B}\setminus c_1\Z$, we have $b\divides c_1$ or $c_1\divides b$. The latter is impossible for $b \not \in c_1\Z$, whence $b \divides c_1$. This implies $b \divides c_1 a_1\in\sB$ for any $a_1\in\sA^1$. By the primitivity of $\sB$, it follows that $b=c_1a_1$ for infinitely many $a_1$, which is impossible and we obtain that $\sB^1$ is indeed primitive. If $\mathscr{B}^1$ is taut, we stop the procedure here and set $\mathscr{B}':=\mathscr{B}^1$. Otherwise, we continue.
%%%%%%%%%%%%%%%%%%%%%%%%%%%
\paragraph{Step 2.}
If $\mathscr{B}^1$ is not taut then, by \cref{ZZ1}, $\sB \setminus c_1\Z$ is not taut. Let $c_2\in\N$ be the smallest number such that
$$
\mathscr{A}^2:=\left\{ \frac{b}{c_2}: b\in\mathscr{B}\setminus c_1\Z \text{ and }c_2 \divides b \right\}
$$
is Behrend (such $c_2$ exists by \cref{Mnajmn}). Note that (by the definition of $c_1$ and $c_2$)
\begin{equation}\label{PPP1}
c_2> c_1 \text{ and }c_2\not\in c_1\Z
\end{equation}
(if $c_1 \divides c_2$ then $c_2 \Z \cap (\sB\setminus c_1\Z)=\emptyset$). Moreover, $\sB \setminus (c_1\Z\cup c_2\Z)=\sB \setminus (c_1\sA^1 \cup c_2\sA^2)$. Let
$$
\sB^2:=(\sB \setminus (c_1\Z\cup c_2\Z))\cup \{c_1,c_2\}=(\mathscr{B} \setminus (c_1\sA^1 \cup c_2 \sA^2))\cup \{c_1,c_2\}.
$$
We claim that $\sB^2$ is primitive. Indeed, if this is not the case then, by the primitivity of $\sB^1$ and by~\eqref{PPP1}, for some $b\in \sB\setminus (c_1\Z\cup c_2\Z)$, we have $b \divides c_2$ or $c_2\divides b$. The latter is impossible for $b\not\in c_2\Z$, whence $b\divides c_2$. This implies $b \divides c_2 a_2$ for any $a_2 \in \sA^2$. By the primitivity of $\sB$, it follows that $b =c_2a_2$ for infinitely many $a_2$, which is impossible and we obtain that $\sB^2$ is indeed primitive. If $\mathscr{B}^2$ is taut, we stop here and set $\mathscr{B}':=\mathscr{B}^2$. Otherwise we continue our construction in a similar way. 
%%%%%%%%%%%%%%%%%%%%%%%%%%%
\paragraph{Step $n$.}
Suppose that from the previous step we have 
\begin{align*}
\sB^{n-1}&=(\sB \setminus (c_1\Z \cup\dots\cup c_{n-1}\Z) )\cup \{c_1,\dots,c_{n-1}\}\\
&=(\sB \setminus (c_1\sA^1 \cup\dots\cup c_{n-1}\sA^{n-1}) )\cup \{c_1,\dots,c_{n-1}\}
\end{align*}
that is primitive but not taut. Then, by \cref{ZZ1}, $\sB \setminus (c_1\Z\cup \dots \cup c_{n-1}\Z)$ is not taut. Let $c_n\in\N$ be the smallest number such that
$$
\sA^n:=\left\{ \frac{b}{c_n}: b\in\mathscr{B} \setminus (c_1\Z\cup \dots\cup c_{n-1}\Z) \text{ and }c_n \divides b \right\}
$$
is Behrend (such $c_n$ exists by \cref{Mnajmn}). Note that (by the definition of $c_1,\dots, c_n$)
\begin{equation}\label{cn}
c_n>c_{n-1} \text{ and }c_n\not\in c_1\Z\cup \dots \cup c_{n-1}\Z.
\end{equation}
Moreover,  
\begin{equation}\label{MMA}
\sB \setminus (c_1\Z\cup\dots \cup c_n\Z)=\sB \setminus (c_1\sA^1\cup\dots \cup c_n\sA^n).
\end{equation}
Let
\begin{align}
\begin{split}\label{BEEN}
\sB^n:&=\sB \setminus (c_1\Z\cup\dots\cup c_n\Z) \cup \{c_1,\dots,c_n\}\\
&=\sB \setminus (c_1\sA^1\cup\dots\cup c_n\sA^n) \cup \{c_1,\dots,c_n\}.
\end{split}
\end{align}
Again, $\sB^n$ is primitive.  If $\mathscr{B}^n$ is taut, we stop the procedure and set $\mathscr{B}':=\mathscr{B}^n$. 
\paragraph{Step $\infty$.}
If $\mathscr{B}^n$ is not taut for all $n\geq 1$, we set
\begin{equation}\label{KKindu6}
\mathscr{B}':=( \mathscr{B}\setminus \bigcup_{n\geq 1}c_n\Z ) \cup \{c_n: n\geq 1\}=( \mathscr{B}\setminus \bigcup_{n\geq 1}c_n\mathscr{A}^n) \cup \{c_n: n\geq 1\},
\end{equation}
where the above equality follows from~\eqref{MMA}. Note that for any $b,b'\in\sB'$ there exists $n\geq 1$ with $b,b'\in \sB^n$. Therefore, by the primitivity of $\sB^n$, $n\geq 1$, also $\sB'$ is primitive.

From now on, for the sake of readability, we will restrict ourselves to the case when $\mathscr{B}'$ is defined by~\eqref{KKindu6}.\footnote{This is the most involved case. When  $\mathscr{B}=\mathscr{B}^n$ for some $n\geq 1$, the proof goes along the same lines, with some simplifications.}
%%%%%%%%%%%%

%%%%%%%%%%%%%%%%%%%%%%%%%%%%%%%%%%%%

\begin{Remark}\label{comniej}
It follows by \eqref{MMA} that 
$$
\sB=(\sB\setminus \bigcup_{n\geq 1}c_n\Z) \cup \bigcup_{n\geq 1}c_n\sA^n.
$$
Therefore, $\mathcal{M}_\mathscr{B}\subset \mathcal{M}_{\mathscr{B}'}$. Moreover, $\eta'\leq \eta$ and $\widetilde{X}_{\eta'}\subset \widetilde{X}_\eta$.
\end{Remark}
%%%%%%%%%%%%%%%%%%%%%%%%%%%%%%%%%%%%

\begin{Lemma}\label{l:3.6}
$\mathscr{B}'$ is taut.
\end{Lemma}
\begin{proof}
Recall that $\sB'$ is primitive. In view of \cref{rozbicie}, it suffices to show that $\sB\setminus \bigcup_{n\geq 1}c_n\Z$ and $\{c_n :n\geq 1\}$ are taut. Suppose that $\sB\setminus \bigcup_{n\geq 1}c_n\Z$ fails to be taut. Then, by \cref{beh3}, for some $c\in\N$ and a Behrend set $\sA$, we have
$$
c\sA \subset \sB \setminus \bigcup_{n\geq 1}c_n\Z.
$$
Therefore, for any $n\geq 1$,
$$
c\sA \subset \sB \setminus (c_1\Z \cup \dots\cup c_n\Z).
$$
By the definition of $c_{n+1}$, we obtain $c\geq c_{n+1}$. Since $n\geq 1$ is arbitrary and the sequence $(c_n)_{n\geq 1}$ is strictly increasing, this yields a contradiction.

Suppose now that $\mathscr{C}:=\{c_n : n\geq 1\}$ fails to be taut. Then, for some $n_0\geq 1$, we have
$$
\bdelta(\cm_{\mathscr{C}})=\bdelta(\cm_{\mathscr{C}\setminus \{c_{n_0}\}}).
$$
Note that by~\eqref{cn}, we have 
$$
c_{n_0}\not\in \bigcup_{n\neq n_0}c_n\Z.
$$
Therefore, by \cref{beh},
$$
\left\{ \frac{c_n}{\gcd(c_n,c_{n_0})} :n\neq n_0 \right\}  \text{ is Behrend}.
$$
We have
$$
\left\{ \frac{c_n}{\gcd(c_n,c_{n_0})} :n\neq n_0 \right\}=\bigcup_{d_{n_0} \divides c_{n_0}}\left\{ \frac{c_n}{d_{n_0}} : n\neq n_0, \gcd(c_n,c_{n_0})=d_{n_0} \right\}.
$$
It follows by \cref{beh2} that at least one of the sets in the union above is Behrend. Denote this set by $\sA(d_{n_0})$ and, for $m> n_0$, define
$$
\sA_m:=\left\{ \frac{c_n}{d_{n_0}} : n\geq m, \gcd(c_n,c_{n_0})=d_{n_0} \right\}.
$$
Since each $\sA_m$ differs from $\sA_{d_{n_0}}$ by at most finitely many elements, it follows by \cref{beh2} that $\sA_m$ is Behrend for $m > n_0$. Let
$$
\sA'_m:=\bigcup_{\substack{n\geq m\\ \gcd(c_n,c_{n_0})=d_{n_0}}}\frac{c_n}{d_{n_0}}\sA^n.
$$
Using \cref{da-er}, \cref{daerdogolne} and the fact that $\sA^n$ is Behrend, we obtain
\begin{multline*}
\bdelta(\cm_{\sA'_m})=\lim_{K\to \infty}\bdelta\Big(\cm_{\bigcup_{ \substack{m\leq n\leq K \\ \gcd(c_n,c_{n_0})=d_{n_0}} } \frac{c_n}{d_{n_0}}\sA^n }\Big)\\
=\lim_{K\to \infty}\bdelta\Big(\cm_{\left\{\frac{c_n}{\gcd(c_n,d_{n_0})} : m\leq n\leq K, \gcd(c_n,c_{n_0})=d_{n_0}\right\}}\Big)=\bdelta(\cm_{\sA_m})=1,
\end{multline*}
since $\sA_m$ is Behrend (the sets $\sA^n$ are the same as in the construction of $\sB'$). By the definition of $\sA'_m$ and $\sA^n$, $n\geq m>n_0$, it follows that
$$
d_{n_0}\sA'_m \subset \bigcup_{n\geq m}c_n\sA^n\subset \sB \setminus \bigcup_{n<m}c_n\Z.
$$
Moreover, by the definition of $c_m$, it follows that $d_{n_0}\geq c_m$, which is impossible as $m\geq n_0$ is arbitrary. This completes the proof.
\end{proof}

%%%%%%%%%%%%%%%%%%%%%%%%%%%%%%%%%%%%%%%%%%%%%%%%%%%%%%%%
\begin{Lemma}\label{KKgest}
We have $\nu_\eta=\nu_{\eta'}$.
\end{Lemma}
\begin{proof}
We will show first that $\un{d}(\cm_\sB)=\un{d}(\cm_{\sB'})$. Let
\begin{align*}
&\mathscr{A}_1:=\mathscr{B}\setminus \bigcup_{n\geq 1}c_n\Z,\ \mathscr{A}_k:=c_{k-1}\mathscr{A}^{k-1}  \text{ for }k\geq 2,\\
&\mathscr{A}'_1:=\mathscr{B}\setminus \bigcup_{n\geq 1}c_n\Z,\ \mathscr{A}'_k:=\{c_{k-1}\}  \text{ for }k\geq 2.
\end{align*}
Then
\begin{equation}\label{postac}
\sB=\bigcup_{n\geq 1}\sA_n \text{ and }\sB'=\bigcup_{n\geq 1}\sA_n'.
\end{equation}
Since each of the sets $\sA^k$, $k\geq 1$, is Behrend, we have
\begin{equation}\label{postac1}
\bdelta(\cm_{\mathscr{A}_1\cup\dots\cup\mathscr{A}_K})=\bdelta(\cm_{\mathscr{A}_1'\cup\dots\cup\mathscr{A}_K'}) \text{ for each }K\geq 1.
\end{equation}
It follows by \eqref{postac}, \eqref{postac1} and by \cref{daerdogolne} that
\begin{multline*}
\un{d}(\cm_\mathscr{B})=\bdelta(\cm_\mathscr{B})=\lim_{K\to \infty}\bdelta(\cm_{\mathscr{A}_1\cup\dots\cup\mathscr{A}_K})\\
=\lim_{K\to\infty}\bdelta(\cm_{\mathscr{A}_1'\cup\dots\cup\mathscr{A}_K'})=\bdelta(\cm_{\mathscr{B}'})=\un{d}(\cm_{\mathscr{B}'}).
\end{multline*}
Moreover, since $\mathcal{M}_\mathscr{B}\subset\mathcal{M}_{\mathscr{B}'}$, it follows that whenever $(N_k)_{k\geq 1}$ satisfies
$$
\lim_{k\to\infty}\frac1{N_k}|\cm_{\mathscr{B}'}\cap[1,N_k]|=\un{d}(\cm_{\mathscr{B}'}),
$$
then
$$
\lim_{k\to\infty}\frac1{N_k}|\cm_{\mathscr{B}}\cap[1,N_k]|=\un{d}(\cm_{\mathscr{B}}).
$$
Since $\eta$ and $\eta'$ differ, along $(N_k)_{k\geq 1}$, on a subset of zero density, it follows by~\cref{quasi-gen} that $\eta$ and $\eta'$ are generic along $(N_k)_{k\geq 1}$ for the same measure, i.e.\ $\nu_\eta=\nu_{\eta'}$.
\end{proof}

%%%%%%%%%%%%%%%%%%%%%%%%
\cref{tautywszystko} follows by \cref{KKgest} and \cref{l:3.6}.

%%%%%%%%%%%%%%%%%%%%%%%%%%%%%%%%%%%%%%%%%%%%%%%%%%%%%%%%%%%%
\subsection{Different classes of $\mathscr{B}$-free numbers}\label{s4}
%%%%%%%%%%%%%%%%%%%%%%%%%%%%%%%%%%%%%%%%%%%%%%%%%%%%%%%%%%%%

%\subsection{Mutual relations}
In Section~\ref{se:klasy}, we defined several classes of $\mathscr{B}$-free numbers and described some basic relations between them. In particular, we showed that
$$
\mathscr{B} \text{ is thin }\Rightarrow \mathscr{B}\text{ has light tails}
$$
and
$$
\mathscr{B}\text{ has light tails (and is primitive) }\Rightarrow \mathscr{B} \text{ is taut}.
$$
We will continue now this discussion. In particular, we will show that the implications converse to the above do not hold. The relations between various classes of $\mathscr{B}$-free numbers for primitive $\mathscr{B}\subset \N$ are summarized in this diagram (all depicted regions are non-empty):
\begin{figure}[h!]
  \centering
    \includegraphics[width=0.3\textwidth]{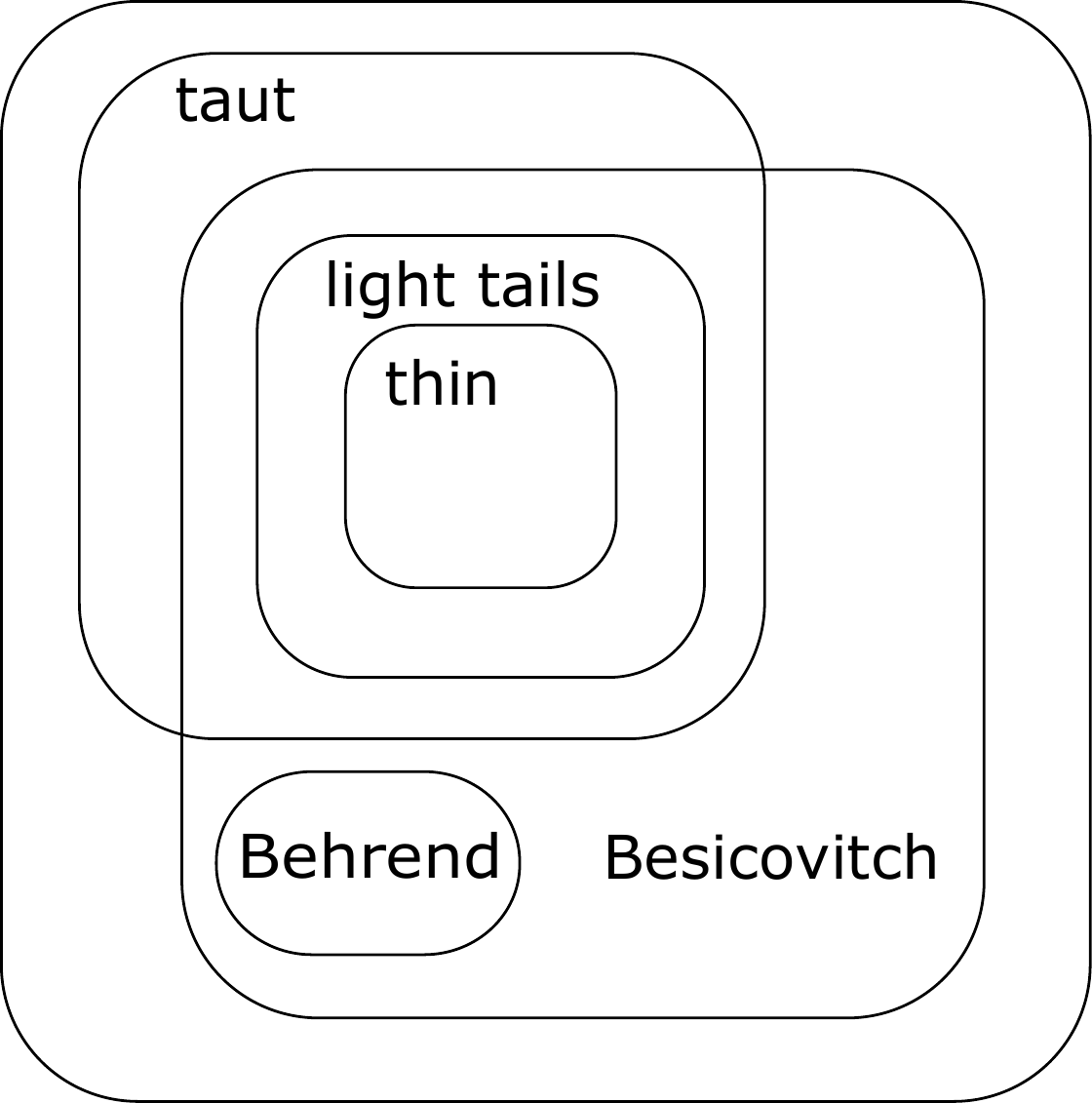}
\end{figure}

\begin{Remark}\label{LLOGOLNE}
Let $\sB,\sB'\subset \N$ be such that:
\begin{itemize}
\item
for each $b'\in\sB'$ there exists $b\in\sB$ such that $b \divides b'$,
\item
for each $b\in\sB$ there exists $b'\in\sB'$ such that $b \divides b'$.
\end{itemize}
Then, clearly, $\cf_\sB \subset \cf_{\sB'}$. Suppose additionally that $\sB$ has light tails and for each $b\in\sB$ the set $\{b'\in\sB' : b\divides b'\}$ is finite. Then, given $K\geq 1$, there exists $N_K\geq 1$ such that 
$$
\text{if  $b\in\sB,b'\in\sB'$, $b\divides b'$ and $b' > N_K$ then $b >K$.}
$$
It follows that
$$
\bigcup_{b'>N_K}b'\Z \subset\bigcup_{b>K}b\Z.
$$
Therefore, if $\sB$ has light tails then also $\sB'$ has light tails. In particular, this applies when $\sB$ is thin (see \cref{prz1sk} below).
\end{Remark}

\begin{Example}[$\mathscr{B}$ has light tails $\centernot\Rightarrow$ $\mathscr{B}$ is thin]\label{prz1sk}
Let $(q_n)_{n\geq1}$ be a thin sequence of primes, i.e., $\sum_{n\ge 1}\frac{1}{q_n}<+\infty$. We arrange the remaining primes into countably many finite pairwise disjoint sets of the form
$\{p_{n,1},p_{n,2},\ldots,p_{n,k_n}\}$ such that
$$
\frac{1}{p_{n,1}}+\frac{1}{p_{n,2}}+\ldots+\frac{1}{p_{n,k_n}}\ge q_n
$$
for any $n$.
Let $\mathscr{B}:=\{q_np_{n,j} : n\in\N, j=1,\ldots, k_n\}$. By \cref{LLOGOLNE}, $\sB$ has light tails. We will show now that $\sB$ is not thin. Indeed,
$$
\sum_{b\in \mathscr{B}}\frac 1b=\sum_{n\ge 1}\left(\frac{1}{q_np_{n,1}}+\frac{1}{q_np_{n,2}}+
\ldots+\frac{1}{q_np_{n,k_n}}\right)\ge \sum_{n\ge 1}1=+\infty.
$$
\end{Example}

\begin{Remark} 
Notice that $\mathscr{B}$ from \cref{prz1sk} is not coprime ($q_np_{n,1}$ and $q_np_{n,2}$ are clearly not coprime).
This is not surprising -- if $\mathscr{B}$ is coprime then it has light tails if and only if it is thin (indeed, in the coprime case the density of $\cf_{\mathscr{B}}$ exists and it is equal $\prod_{k\geq 1}(1-\frac{1}{b_k})$, see, e.g.,\ \cite{MR1414678}). Note however that $\mathscr{B}$ above is primitive.\end{Remark}

\begin{Remark}\label{exaten}
Let $\mathscr{B}$ be as in \cref{prz1sk}. It follows by \cref{tautnotdelta} that $\nu_\eta\neq \delta_{(\dots,0,0,0,\dots)}$. 
\end{Remark}

\begin{Prop}\label{nobe}
$\mathscr{B} \text{ is taut } \centernot\Rightarrow \mathscr{B}\text{ is Besicovitch}.$
\end{Prop}
In the proof, we will use the following lemma:
\begin{Lemma}\label{l:3.7}
Let $\sB\subset \N$ and let $\sB'$ be as in the proof of \cref{tautywszystko}. Then $\sB$ is Besicovitch whenever $\sB'$ is Besicovitch.
\end{Lemma}
\begin{proof}
Recall that in the notation from the proof of \cref{tautywszystko}, we have
$$
\sB=(\sB\setminus \bigcup_{n\geq 1}c_n\sA^n)\cup \bigcup_{n\geq 1}c_n\sA^n
$$
and
$$
\sB'=(\sB\setminus \bigcup_{n\geq 1}c_n\sA^n)\cup \{c_n :n\geq 1\}.
$$
It follows by \cref{da-er}, by the fact that the sets $\sA^n$ for $n\geq 1$ are Behrend and by \cref{daerdogolne} that we have
\begin{multline*}
\un{d}(\cm_\sB)=\lim_{K\to \infty}\bdelta(\cm_{(\sB\setminus \bigcup_{n\geq 1}c_n\sA^n) \cup \bigcup_{n\leq K}c_n\sA_n})\\
=\lim_{K\to \infty}\bdelta(\cm_{(\sB \setminus\bigcup_{n\geq 1}c_n\sA^n)\cup \{c_n :n\leq K\}})=\un{d}(\cm_{\sB'}).
\end{multline*}
Therefore, is $\sB'$ is Besicovitch, we obtain $\un{d}(\cm_\sB)=d(\cm_{\sB'})$. On the other hand, by \cref{tautywszystko}, we have $\cm_{\sB}\subset \cm_{\sB'}$ and it follows that $\ov{d}(\cm_\sB)\leq d(\cm_{\sB'})$. We obtain $\ov{d}(\cm_\sB)\leq \un{d}(\cm_\sB)$ and conclude that also $\sB$ must be Besicovitch.
\end{proof}

\begin{proof}[Proof of \cref{nobe}]
Consider $\mathscr{B}$ that fails to be Besicovitch. By Lemma~\ref{l:3.7}, the associated set $\mathscr{B}'$ defined as in~\eqref{KKindu6} also fails to be Besicovitch. Moreover, in view of Lemma~\ref{l:3.6}, $\mathscr{B}'$ is taut.
\end{proof}

Since, as noted in \cref{se:klasy}, each $\mathscr{B}$ with light tails is automatically Besicovitch, we have the following immediate consequence of \cref{nobe}:
\begin{Cor}
$\mathscr{B}$ is taut $\centernot\Rightarrow$ $\mathscr{B}$ has light tails.
\end{Cor}
The rest of this section is devoted to the proof of the following more subtle result:

\begin{Th}\label{SK-1}
$\mathscr{B}$ is taut and Besicovitch $\centernot\Rightarrow \mathscr{B}$ has light tails.
\end{Th}
To prove \cref{SK-1}, we will need three lemmas.
\begin{Lemma}\label{rozklad}
Let $\mathscr{R}$ be a union of finitely many arithmetic progressions with steps $d_1,...,d_r$. Then $\mathscr{R}$ is a union of finitely many pairwise disjoint  arithmetic progressions of steps $\lcm(d_1,\dots, d_r)$.
\end{Lemma}
\begin{proof}
Let $\mathscr{R}=\bigcup_{i=1}^r (d_i\Z+a_i)$. Notice that
\begin{equation}\label{UII}
\mathscr{R}=\bigcup_{\{i_1,\dots,i_s\}\in I}\bigcap_{i=1}^s (d_{i_s}\Z +a_{i_s}),
\end{equation}
where $\{i_1,\dots,i_s\}\in I$ if and only if $\bigcap_{i=1}^s (d_{i_s}\Z +a_{i_s})\neq\emptyset$ and $(d_j\Z+a_j)\cap\bigcap_{i=1}^s (d_{i_s}\Z +a_{i_s})=\emptyset$ for any $j\not\in \{i_1,\dots,i_s\}$. Moreover, the elements of the union in~\eqref{UII} are pairwise disjoint. Finally, notice that if $a\in \bigcap_{i=1}^s (d_{i_s}\Z +a_{i_s})$ then, by \cref{prosciutto},
\begin{align}
\begin{split}\label{UII1}
\bigcap_{i=1}^s (d_{i_s}\Z +a_{i_s})&=\lcm(d_{i_1},\dots,d_{i_s})\Z +a\\
&=\bigcup_{\ell=0}^{L-1}(\lcm(d_1,\dots,d_r)\Z+ \ell\lcm(d_{i_1},\dots,d_{i_s}) +a),
\end{split}
\end{align}
where $L=\lcm(d_1,\dots,d_r)/\lcm(d_{i_1},\dots,d_{i_s})$ and the elements of the union~\eqref{UII1} are pairwise disjoint.
\end{proof}
\begin{Lemma}\label{independent}
Assume that $B,C\subset \N$ are thin, with $\gcd(b,c)=1$ for any $b\in B$, $c\in C$.
Let $BC:=\{bc:b\in B, c\in C\}$. Then
\begin{equation}\label{teza:ka1}
d(\mathcal{M}_{BC})=d(\mathcal{M}_{B}\cap\mathcal{M}_{C})=d(\mathcal{M}_{B})d(\mathcal{M}_{C}).
\end{equation}
\end{Lemma}
\begin{proof}
Since $\lcm(b,c)=bc$  for any $b\in B$ and  $c\in C$, it follows that
$$
 \mathcal{M}_{BC}=\mathcal{M}_{B}\cap\mathcal{M}_{C}.
$$
It remains to show the right hand side equality in~\eqref{teza:ka1} and it is enough to show its validity for finite sets $B$, $C$ (since $BC$ is thin, it is Besicovtich and we can use \cref{da-er} to pass to a limit).

Let $B=\{b_1, \dots ,b_n\}$, $C=\{c_1, \dots ,c_m\}$ and set
$$
b':=\lcm(b_1,\dots,b_n),\ c':=\lcm(c_1,\dots,c_m).
$$
Then, by \cref{rozklad},
$$
\mathcal{M}_{B}=\bigcup\limits_{r\in R}(b'\Z+r),\;\mathcal{M}_{C}=\bigcup\limits_{s\in S}(c'\Z+s)
$$
for some finite sets $R,S\subset \N$. Note that
\begin{equation}\label{1.1}
d(\mathcal{M}_{B})=\frac{|R|}{b'},\; d(\mathcal{M}_{C})=\frac{|S|}{c'}.
\end{equation}
Since $\gcd(b',c')=1$, we get
\begin{equation}\label{1.2}
d((b'\Z+r)\cap (c'\Z+s))=\frac{1}{b'c'}
\end{equation}
for any $r\in R$, $s\in S$. Hence, by  \eqref{1.2} and \eqref{1.1}, we obtain
$$
d(\mathcal{M}_{B}\cap\mathcal{M}_{C})
=d(\bigcup\limits_{(r,s)\in R\times S}(b'\Z+r)\cap (c'\Z+s))
=\frac{|R\times S|}{b'c'}=d(\mathcal{M}_{B})d(\mathcal{M}_{C})
$$
and the result follows.
\end{proof}

\begin{Lemma}\label{paczki}
Let $P\subset\N$ be pairwise coprime with $\sum_{p\in P}\nicefrac{1}{p}=+\infty$. For any $0<\beta<1$ there exists a  finite (resp.\ infinite and thin) set $P'\subset P$  such that
$$
\beta< d(\mathcal{M}_{ P'})< 1.
$$
\end{Lemma}
\begin{proof}
For $n\geq 1$, let $P_n:=\{p\in P : p\leq n\}$. By \cref{da-er}, we have
$$
\lim_{n\to\infty}d(\cm_{P_n})=d(\cm_P)=1.
$$
Therefore, for $n\geq 1$ large enough, we have $\beta<d(\cm_{P_n})<1$ and we can take $P':=P_n$ to obtain a finite set satisfying the assertion. To obtain an infinite set $P'$, let the sequence $(p_m)_{m\geq 1}\subset P$ be such that $d(\cm_{P_n})+\sum_{m\geq 1}\nicefrac{1}{p_m}<1$ and take $P':=P_n\cup\{p_m : m\geq 1\}$.
\end{proof}

\paragraph{Construction.}\label{cons1}
Fix $0<\gamma<1$ and choose a sequence $(\gamma_k)_{k\geq 1}\subset (0,1)$ such that  $\prod_{k\geq 1}\gamma_k=\gamma$ (for instance, $\gamma_k=\gamma^{1/2^k}$). Applying Lemma~\ref{paczki}, we construct a collection $\{B_k,C_k:k\in\N\}$ of pairwise disjoint thin sets of primes such that
\begin{equation}\label{dotb}
\gamma_k<d(\mathcal{M}_{B_k})<1\text{ for }k\geq 1
\end{equation}
 and
\begin{equation}\label{dotc}
1-\frac{1}{k}<d(\mathcal{M}_{C_k})\text{ for }k\geq 1.
\end{equation}
Let
\begin{equation}\label{konstrukcja}
\mathscr{B}:=B_1C_1\cup B_1B_2C_2\cup \dots \cup B_1 \dots B_nC_n \cup \dots
\end{equation}
Notice that $B_1C_1\cup B_1B_2C_2\cup \dots \cup B_1 \dots B_nC_n $ is thin for any $n\in\N$.

\begin{proof}[Proof of \cref{SK-1}]
Let $\mathscr{B}$ be defined as in~\eqref{konstrukcja}. We claim the following:
\begin{enumerate}[(a)]
\item\label{LLA}
$\mathscr{B}$ is Besicovitch,
\item\label{LLB}
$\mathscr{B}$ does not have light tails,
\item\label{LLC}
$\mathscr{B}$ is taut.
\end{enumerate}

We will first prove~\eqref{LLA}. For $k\geq m$, we have
$$
\mathcal{M}_{B_1 \dots B_kC_k}\subset \mathcal{M}_{B_1 \dots B_k}\subset \mathcal{M}_{B_1 \dots B_m}.
$$
Thus,
\begin{equation}\label{a1}
\overline{d}(\mathcal{M}_{\bigcup_{k\geq m+1}B_1 \dots B_kC_k}\setminus \mathcal{M}_{B_1 \dots B_mC_m})\le d(\mathcal{M}_{B_1 \dots B_m}\setminus \mathcal{M}_{B_1 \dots B_mC_m}).
\end{equation}
By Lemma \ref{independent} and by (\ref{dotc}), we get
$$
d(\mathcal{M}_{B_1 \dots B_mC_m})=d(\mathcal{M}_{B_1 \dots B_m})d(\mathcal{M}_{C_m})\ge d(\mathcal{M}_{B_1 \dots B_m})(1-\frac{1}{m}),
$$
whence
\begin{equation}\label{a2}
d(\mathcal{M}_{B_1 \dots B_m}\setminus \mathcal{M}_{B_1 \dots B_mC_m})\le \frac{1}{m}d(\mathcal{M}_{B_1 \dots B_m})\le \frac{1}{m}.
\end{equation}
Using (\ref{a1}) and (\ref{a2}), we obtain
$$
\overline{d}(\mathcal{M}_{\bigcup_{i=m+1}^{\infty}B_1 \dots B_iC_i}\setminus \mathcal{M}_{B_1 \dots B_mC_m})\le \frac{1}{m}.
$$
In view of \cref{uwa9wrz1}, this implies that $\mathscr{B}$ is Besicovitch.

We will now show~\eqref{LLB}. By Lemma \ref{independent}, (\ref{dotb}) and (\ref{dotc}), we have
$$
d(\mathcal{M}_{B_1 \dots B_mC_m})\ge \gamma_1 \dots \gamma_m(1-\frac{1}{m})\rightarrow \gamma>0
$$
as $m\rightarrow +\infty$, which yields~\eqref{LLB}.

It remains to prove~\eqref{LLC}. Suppose that $\mathscr{B}$ is not taut. Since $\mathscr{B}$ is primitive, it follows by \cref{beh3} that for some $c\in\N$ and a Behrend set $\mathscr{A}\subset \N\setminus \{1\}$, we have $c\mathscr{A}\subset \mathscr{B}$. Let $m\in\N$ be such that $c$ is coprime to all elements of $B_{m+1}$ (such $m$ exists since $B_n$, $n\in\N$, are pairwise disjoint sets of primes). Let
$$
\mathscr{B}_1:=B_1C_1\cup\dots\cup B_1B_2\dots B_mC_m \text{ and }\mathscr{B}_2:=\bigcup_{n>m}B_1B_2\dots B_nC_n.
$$
Then clearly, $\mathscr{B}=\mathscr{B}_1\cup \mathscr{B}_2$. Moreover, let
$$
\mathscr{A}_1:=\left\{\frac{b}{c}: b\in c\mathscr{A}\cap \mathscr{B}_1\right\}\text{ and } \mathscr{A}_2:=\left\{\frac{b}{c} : b\in c\mathscr{A}\cap \mathscr{B}_2\right\}.
$$
Then clearly, $\mathscr{A}=\mathscr{A}_1\cup\mathscr{A}_2$. Since $\mathscr{B}_1$ is thin, it follows by \cref{LLtauty2} and by \eqref{LLtauty1} that $\mathscr{B}_1$ is taut. Therefore, since $c\mathscr{A}_1\subset \mathscr{B}_1$, it follows by \cref{beh3} that $\mathscr{A}_1$ is not Behrend. Since $\mathscr{A}$ is Behrend, we obtain by \cref{beh2}  that $\mathscr{A}_2$ must be Behrend. Moreover, we have $c\mathscr{A}_2\subset \mathscr{B}_2$. Take $a\in\mathscr{A}_2$. Since $c$ is coprime to each element of $B_{m+1}$, it follows that $a\in \cm_{B_{m+1}}$. Hence, $\cm_{\mathscr{A}_2}\subset \cm_{B_{m+1}}$, which is impossible since $d(\cm_{\mathscr{A}_2})=1$, whereas $d(\cm_{B_{m+1}})<1$ since $B_{m+1}$ is thin. We conclude that $\mathscr{B}$ is taut, which completes the proof.
\end{proof}
%%%%%%%%%%%%%%%%%%%%%%%%%%%%%%%%%%%%%%%%%%%%%%%%%%%%%%%%%%%%

\subsection{Tautness and combinatorics (\cref{OOJ} -- first steps)}\label{Jfirst}
%%%%%%%%%%%%%%%%%%%%%%%%%%%%%%%%%%%%%%%%%%%%%%%%%%%%%%%%%%%%
Since $\eta\in X_\mathscr{B}$, a natural question arises how many residue classes are missing on $\text{supp }\eta \bmod b_k$, $k\geq 1$. We will answer this question in the class of taut sets $\mathscr{B}$.
Recall first the following result:
\begin{Th}[Dirichlet]\label{Dirichlet}
Let $a,r\in\N$. If $\gcd(a,r)=1$ then $a\Z+r$ contains infinitely many primes. Moreover, $\sum_{p\in (a\Z+r)\cap\mathscr{P}}1/p=+\infty$.
\end{Th}
Since each set containing a pairwise coprime set with divergent sum of reciprocals is automatically Behrend, we obtain the following:
\begin{Cor}\label{jebe}
Let $a,r\in\N$. If $\gcd(a,r)=1$ then the set $(a\Z+r)\cap \mathscr{P}$ is Behrend.
\end{Cor}

\begin{Prop}\label{druggi2}
Assume that $\mathscr{B}\subset \N$ is taut, $a\in\N$ and $1\leq r\leq a$. If
\begin{equation}\label{wszy}
a\Z+ r \subset \bigcup_{b\in \mathscr{B}}b\Z
\end{equation}
then there exists $b\in\mathscr{B}$ such that $b\divides \gcd (a,r)$. In particular, if $a\in\mathscr{B}$ then $r=a$.
\end{Prop}

\begin{proof}
Suppose that $a\in\N$ and $1\leq r\leq a$ are such that~\eqref{wszy} holds.
Let $d:=\gcd(a,r)$, $a':=a/d$, $r':=r/d$, i.e.\ we have
$$
d\cdot (a'\Z +r')\subset \bigcup_{b\in\mathscr{B}}b\Z .
$$
Applying \cref{jebe} to $a'$ and $r'$, we obtain $d(\cm_{a'\Z+r'})=1$, whence $\bdelta(\cm_\sB)=\bdelta(\cm_{\sB\cup\{d\}})$. If $d\in\cm_\sB$, then there exists $b\in\sB$ such that $b\divides d$, whence $b\divides \gcd (a,r)$. Suppose now that $d\not\in\cm_\sB$. Then, by \cref{beh}, we have that $\mathscr{B}'(d)$ is Behrend. Moreover, $1\not\in\mathscr{B}'(d)$. Since $\gcd(b,d)$, $b\in\mathscr{B}$, takes only finitely many values, we can represent $\mathscr{B}'(d)$ as a finite union:
$$
\mathscr{B}'(d)=\bigcup_{c\divides d} \left\{ b/c : b\in\mathscr{B}, \gcd(b,d)=c \right\}.
$$
Therefore, in view of \cref{beh2}, for some $c\divides d$, the sequence
$$
\mathscr{A}:=\left\{b/c : b\in\mathscr{B}, \gcd(b,d)=c \right\}
$$
is Behrend. Hence $\mathscr{B}\supset c\mathscr{A}$, with $\mathscr{A}\subset \N\setminus\{1\}$ that is Behrend. This however, in view of \cref{beh3}, contradicts the assumption that $\mathscr{B}$ is taut.

Suppose now that $a\in\mathscr{B}$ and~\eqref{wszy} holds. By the first part of the proof, we have $b \divides \gcd(a,r)$ for some $b\in\sB$. It follows that $b\divides a$ and, since $a,b\in\sB$, by the primitivity of $\sB$, we obtain $a=b$. Therefore, using the relation $b\divides \gcd(a,r)$, we conclude that $b\divides r$ and, since $1\leq r\leq b$, this yields $r=b$.
\end{proof}
\begin{Remark}\label{druggi22}
Let $N\geq 1$. Note that the assertion of \cref{druggi2} remains true if we replace condition~\eqref{wszy} with
$$
(a\Z+ r) \cap [N,\infty) \subset \bigcup_{b\in \mathscr{B}}b\Z.
$$
Indeed, by~\cref{jebe} and \cref{beh2}, $(a\Z+ r) \cap [N,\infty)\cap \mathscr{P}$ remains Behrend and we repeat the rest of the proof of \cref{druggi2}.
\end{Remark}

\begin{Cor}\label{gdzieeta}
Assume that $\sB\subset \N$ is taut. Then, for each $b\in\sB$ and $1\leq r\leq b-1$, there exists infinitely many $m\in \cf_\sB$ such that $m\equiv r\bmod b$. In particular, $\eta\in Y$.
\end{Cor}
\begin{proof}
Fix, $N\geq 1$, $b\in\sB$ and consider $b\Z+r$ for $1\leq r\leq b-1$. By \cref{druggi2} and \cref{druggi22}, $(b\Z+r)\cap [N,\infty)\not\subset \cm_\sB$, i.e.,
$$
(\cf_\sB \cap [N,\infty)) \bmod b = \{1,\dots, r-1\},
$$
and the result follows.
\end{proof}
\begin{Remark}
Note that if $\eta\in Y$ then $\sB$ is primitive. Indeed, if $\sB$ is not primitive then, for some $b,b'\in\sB$, we have $b\divides b'$. If
$|\text{supp }\eta \bmod b'|=b'-1$ then  $|\text{supp }\eta \bmod b|=b$. The latter is impossible as $\eta\in X_\sB$ and it follows that $\eta\not\in Y$.
\end{Remark}
The following example shows that the converse of \cref{gdzieeta} does not hold:
\begin{Example}\label{ex414}
Consider $\{(p_i,r_i):i\geq 1 \}=\{(p,r):p\in \mathscr{P}, 0<r<p\}$. Every progression $p_i\Z+r_i$ contains infinitely many primes; given $i\geq 1 $ let, for $n\geq 1$, 
$$
q_i^n\in (p_i\Z+r_i)\cap\mathscr{P} \text{ be such that }q_i^n>2^n\cdot i^2.
$$
We set $\mathscr{B}:=\mathscr{P}\setminus\{q_i^n:i,n\geq 1 \}$. Since $\sum_{i,n\geq 1 }\frac{1}{q_i^n}<\infty$, it follows that $\sB$ is Behrend, so, in particular, $\sB$ is not taut.

Let $b\in\mathscr{B}$ and $0<r<b$ and let $i\geq 1$ be such that $(b,r)=(p_i,r_i)$. Then, for each $n\geq 1$,  $q_i^n \equiv r \bmod b$ by the choice of $q_i^n$. Moreover, $q_i^n\in\mathcal{F}_\mathscr{B}$ since it is a prime not belonging to $\mathscr{B}$.
\end{Example}

In~\cite{MR3356811} it has been proved that for $\mathscr{B},\mathscr{B}'\subset\N$ coprime and thin the following holds:
\begin{itemize}
\item $X_\mathscr{B}\subset X_{\mathscr{B}'}$ $\iff$ for each $b'\in\mathscr{B}'$ there exists $b\in\mathscr{B}$ with $b \divides b'$,
\item
$X_\mathscr{B}=X_{\mathscr{B}'}$ $\iff$ $\mathscr{B}=\mathscr{B}'$.
\end{itemize}
We will now extend these results to the case of taut sets.

\begin{Cor}\label{zawiewnio}
Let $\sB,\sB'\subset \N$ and suppose that $\mathscr{B}$ is taut. Then the following conditions are equivalent:
\begin{enumerate}[(a)]
\item $X_\mathscr{B} \subset X_{\mathscr{B}'}$,\label{J:A}
\item for each $b'\in\mathscr{B}'$ there exists $b\in\mathscr{B}$ with $b\divides b'$,\label{J:B}
\item $\eta\leq \eta'$,\label{J:C}
\item $\widetilde{X}_\eta \subset \widetilde{X}_{\eta'}$,\label{J:D}
\item $\eta\in \widetilde{X}_{\eta'}$,\label{J:E}
\item $\eta\in X_{\mathscr{B}'}$.\label{J:F}
\end{enumerate}
\end{Cor}
\begin{proof}
Clearly, we have \eqref{J:B} $\Rightarrow$ \eqref{J:C} $\Rightarrow$ \eqref{J:D} $\Rightarrow$ \eqref{J:E} $\Rightarrow$ \eqref{J:F} and \eqref{J:A} $\Rightarrow$ \eqref{J:F}. Therefore, to complete the proof it suffices to show \eqref{J:B} $\Rightarrow$ \eqref{J:A} and \eqref{J:F} $\Rightarrow$ \eqref{J:B}. 

Suppose that \eqref{J:B} holds and let $A\subset \N$ be $\sB$-admissible. Take $b'\in\sB$ and let $b\in\sB$ be such that $b\divides b'$. It follows by the $\{b\}$-admissibility of $A$ that for some $0\leq r\leq b-1$, we have $(b\Z+r)\cap A=\emptyset$, so all the more, we have $(b'\Z+r)\cap A=\emptyset$, i.e., $A$ is $\{b'\}$-admissible and \eqref{J:A} follows.

Suppose that \eqref{J:F} holds. Then, for each $b'\in\mathscr{B}'$ there exists $1\leq r'\leq b'$ such that $r'\not\in\cf_\sB \bmod b'$, i.e.,
$$
b'\Z+r' \subset \bigcup_{b\in\mathscr{B}}b\Z
$$
It follows by \cref{druggi2} that there exists $b\in\sB$ such that $b \divides \gcd(b',r')$, so, in particular, $b\divides b'$, i.e. \eqref{J:B} holds.
\end{proof}
\begin{Cor}\label{zawiewnio1}
Suppose that $\mathscr{B},\mathscr{B}'$ are taut. Then the following conditions are equivalent:
\begin{enumerate}[(a)]
\item $X_\mathscr{B} = X_{\mathscr{B}'}$,\label{J:A1}
\item $\mathscr{B}=\mathscr{B}'$,\label{J:B1}
\item $\eta=\eta'$,\label{J:C1}
\item $\widetilde{X}_\eta = \widetilde{X}_{\eta'}$,\label{J:D1}
\item $\eta\in \widetilde{X}_{\eta'}$ and $\eta'\in \widetilde{X}_{\eta}$,\label{J:E1}
\item $\eta\in X_{\mathscr{B}'}$ and $\eta'\in X_{\mathscr{B}}$,\label{J:F1}
\item $X_\eta=X_{\eta'}$.\label{J:G1}
\end{enumerate}
\end{Cor}
\begin{proof}
We have immediately \eqref{J:B} $\Rightarrow$ \eqref{J:C} $\Rightarrow$ \eqref{J:D} $\Rightarrow$ \eqref{J:E} $\Rightarrow$ \eqref{J:F}, \eqref{J:B} $\Rightarrow$ \eqref{J:A} $\Rightarrow$ \eqref{J:F} and \eqref{J:C1} $\Rightarrow$ \eqref{J:G1} $\Rightarrow$ \eqref{J:D1}. We will show now the remaining implication \eqref{J:F1} $\Rightarrow$ \eqref{J:B1}. By the corresponding implication in \cref{zawiewnio}, for any $b\in\mathscr{B}$
there exist $b'\in\mathscr{B}'$ and $b''\in\mathscr{B}$ such that $b'' \divides b' \divides  b$. Since $\mathscr{B}$ is taut, it is, in particular, primitive which yields $b=b'=b''$, i.e.\ $\mathscr{B}\subset \mathscr{B}'$. Reversing the roles of $\mathscr{B}$ and $\mathscr{B}'$, we obtain $\mathscr{B}=\mathscr{B}'$.
\end{proof}

%%%%%%%%%%%%%%%%%%%%%%%%%%%%%%%%%%%%%%%%%%%%%%%%%%%%%%%%%%%%%%%%%%%%%%%%%%%
%%%%%%%%%%%%%%%%%%%%%%%%%%%%%%%%%%%%%%%%%%%%%%%%%%%%%%%%%%%%%%%%%%%%%%%%%%%
%%%%%%%%%%%%%%%%%%%%%%%%%%%%%%%%%%%%%%%%%%%%%%%%%%%%%%%%%%%%%%%%%%%%%%%%%%%

%%%%%%%%%%%%%%%%%%%%%%%%%%%%%%%%%%%%%%%%%%%%%%%%%%%%%%%%%%%%
%%%%%%%%%%%%%%%%%%%%
\section{Heredity (proofs of \cref{TTD} and \cref{OOF})}\label{s6}
By \cref{hertoprox}, $(S,X_\eta)$ is proximal whenever $X_\eta$ is hereditary. The converse to that does not hold, cf.\ Example~\ref{ex:2.4} (proximality follows from \cref{AUR}). In this section, we will show however that the proximality of $(S,X_\eta)$ and the heredity of $X_\eta$ are equivalent when $\mathscr{B}$ has light tails.

\begin{Def}
We say that $A\subset \N$ is \emph{$\eta$-admissible} whenever
\begin{equation}\label{eta-adm}
\{k+1,\ldots,k+n\}\cap\mathcal{F}_{\mathscr{B}}=A+k
\end{equation}
for some $k,n\in\N$ (in other words, $\text{supp }\eta[k+1,k+n]=A+k$).
\end{Def}
\begin{Def}
We say that $A$ satisfies condition~\eqref{Ther} whenever
\begin{equation}\label{Ther}
\parbox{0.75\textwidth}{ 
there exists $\{n_b \in\Z \colon b\in\sB\}$ such that $A\cap (b\Z+n_b)=\emptyset$ and $\gcd(b,b')\divides n_b-n_{b'}$ for any $b,b'\in\sB$.
}\tag{T$_{\text{her}}$}
\end{equation}
\end{Def}

Our main goal in this section is to prove the following:
\begin{Th}\label{charakt}
Assume that $\mathscr{B}\subset \N$ has light tails and satisfies~\eqref{Au1}. Let $n\in\N$ and  $A\subset\{1,\ldots,n\}$. The following conditions are equivalent:
\begin{enumerate}[(a)]
\item\label{UIA}
$A$ satisfies~\eqref{Ther},
\item\label{UIB}
$A$ is $\eta$-admissible.
\end{enumerate}
In particular, $X_\eta$ is hereditary, i.e. $X_\eta=\widetilde{X}_\eta$. 
\end{Th}
\begin{Remark}\label{IU1}
Clearly, if $A'\subset A\subset \Z$ and $A$ satisfies \eqref{Ther} then also $A'$ satisfies \eqref{Ther}. Thus, \cref{TTD}, i.e., the assertion that $X_\eta$ is hereditary in \cref{charakt}, follows immediately by the equivalence of \eqref{UIA} and \eqref{UIB}.
\end{Remark}

As an immediate consequence of \cref{AUR} and \cref{charakt}, we have:
\begin{Cor}\label{wniosekdzis}
Assume that $\mathscr{B}\subset \N$ has light tails. Then $X_\eta$ is hereditary if and only if $(S,X_\eta)$ is proximal.
\end{Cor}

\begin{Example}[cf.\ \cref{X_eta<>X_B}]\label{X_eta<>X_Ba}
Let $\mathscr{B}\subset \N$ be as in Example~\ref{X_eta<>X_B}. If additionally $\sB$ has light tails and satifies~\eqref{Au1}, then, by \cref{charakt}, $X_\eta=\widetilde{X}_\eta$. E.g.\ one can take $\mathscr{B}=\{4,6\}\cup\{p^2 : p\in\mathscr{P},\ p>12 \}$.

On the other hand, if~\eqref{Au1} fails then, by \cref{AUR}, $X_\eta$ fails to be proximal. Hence, by \cref{hertoprox}, $X_\eta$ also fails to be hereditary. E.g.\ one can take $\mathscr{B}=\{4,6\}\cup\{5p^2 : p\in\PP,\ p>12 \}$.
\end{Example}

We leave the following question open:
\begin{Question}
Are the heredity of $X_\eta$ and proximality of $(S,X_\eta)$ the same whenever $\mathscr{B}$ is taut?
\end{Question}

\begin{Remark}
Notice that $\mathscr{B}$ from the construction on page~\pageref{cons1} satisfies condition~\eqref{Au1} whenever $B_1$, $C_1$ are infinite, i.e.\ $(S,X_\eta)$ is proximal. We do not know whether in this example $X_\eta=\widetilde{X}_\eta$.
\end{Remark}

For the proof of \cref{charakt}, we will need several auxiliary results.
\begin{Lemma}\label{onedir}
Let $n\in\N$ and suppose that $A\subset \{1,\dots,n\}$ is $\eta$-admissible. Then $A$ satisfies~\eqref{Ther}.
\end{Lemma}
\begin{proof}
Suppose that
$$
\{k+1,\ldots,k+n\}\cap\mathcal{F}_{\mathscr{B}}=A+k.
$$
for some $k$. For $b\in\mathscr{B}$, let $n_b:=-k$. Since for any $i\in A$, $i+k\in \cf_\sB$, we have $i+k\notin b\Z $ for any $b\in\mathscr{B}$. This means that $i\notin b\Z  -k=b\Z +n_b$. It follows immediately that $A$ satisfies~\eqref{Ther}.
\end{proof}

Lemma~\ref{onedir} gives the implication (b) $\Rightarrow$ (a) in the assertion of \cref{charakt}. Now, we will cover the converse implication.
For $n\geq 1$, let
$$
\mathscr{B}^{(n)}:=\{b\in \mathscr{B} : p\le n \text{ for any }p\in{\rm Spec}(b)\},
$$
where $\text{Spec}(b)$ stands for the set of all prime divisors of $b$.\footnote{For $A\subset \N$ the set $\text{Spec}(A)$ is defined as the union of $\text{Spec(a)}$, $a\in A$.} Our main tools are the following two results:
\begin{Prop}\label{new1}
Assume that $\mathscr{B}\subset \N$ satisfies~\eqref{Au1} and $\mathscr{B}^{(n)}\subset \mathscr{A}\subset \mathscr{B}$.
Suppose that
\begin{equation}\label{dziedz:as}
\{k+1,\ldots,k+n\}\cap {\cal M}_{\mathscr{A}}=\{k+i_1,k+i_2,\ldots,k+i_r\}
\end{equation}
for some $1\le i_1,\ldots,i_r\le n$, $r<n$.\footnote{If $r=0$, we interpret the right hand side of \eqref{dziedz:as} as the empty set.}
Then, for arbitrary $i_0\in\{1,\ldots,n\}$, there exist $\mathscr{B}^{(n)}\subset \mathscr{A}'\subset \mathscr{B}$ and $k'\in\Z$ such that
$$
\{k'+1,\ldots,k'+n\}\cap {\cal M}_{\mathscr{A}'}=\{k'+i_0,k'+i_1,\ldots,k'+i_r\}.
$$
\end{Prop}
\begin{Prop}\label{new1a}
Assume that $\mathscr{B}\subset \N$ has light tails and $\mathscr{B}^{(n)}\subset \mathscr{A}\subset \mathscr{B}$.
Suppose that
\begin{equation}\label{dziedz:as1}
\{k+1,\ldots,k+n\}\cap {\cal M}_{\mathscr{A}}=\{k+i_0,k+i_1,\ldots,k+i_r\}
\end{equation}
for some $1\le i_0,\ldots,i_r\le n$, $r<n$. Then the density of $k'\in\N$ such that
$$
\{k'+1,\ldots,k'+n\}\cap {\cal M}_{\mathscr{B}}=\{k'+i_0,k'+i_1,\ldots,k'+i_r\}
$$
is positive.\footnote{For the purposes of this section it would be sufficient to know that such $k'$ exists. We will use this result in its full form later.\label{sto:here}}
\end{Prop}

Before we give the proofs of \cref{new1} and \cref{new1a}, we will show how these two results yield the implication (a) $\Rightarrow$ (b) in \cref{charakt}. Notice first that an inductive procedure applied to \cref{new1}, together with \cref{new1a}, implies immediately the following:
\begin{Cor}\label{tenwnio}
Assume that $\mathscr{B}$ has light tails and satisfies~\eqref{Au1}. Assume that $\mathscr{B}^{(n)}\subset \mathscr{A}\subset \mathscr{B}$. Suppose that
\begin{equation}\label{UI3}
\{k+1,\ldots,k+n\}\cap {\cal M}_{\mathscr{A}}=k+C
\end{equation}
for some $C\subset\{1,.\ldots,n\}$.
Then, for arbitrary set $C'$ such that $C\subset C'\subset \{1,\ldots,n\}$, the density of the set of $k'\in\Z$ such that
$$
\{k'+1,\ldots,k'+n\}\cap {\cal M}_{\mathscr{B}}=k'+C'
$$
is positive.
\end{Cor}

We will present now some auxiliary results.

\begin{Lemma}\label{finite}
Let $A\subset \N$ be primitive, with $\text{Spec}(A)$ finite. Then $A$ is also finite.
\end{Lemma}

\begin{proof}
The proof will use induction on $|\text{Spec}(A)|$. Clearly, if $|\text{Spec}(A)|=1$ then also $|A|=1$. Suppose that the assertion holds for any set $A$ with $|\text{Spec}(A)|\leq n-1$. Let now $A$ be primitive with $|\text{Spec}(A)|=n$, i.e.\
$$
\text{Spec}(A)=\{p_1,\dots,p_n\}\subset\mathscr{P}.
$$
For $k\geq 0$, let
\begin{align*}
A^{(k)}:=&\{a\in A : k=\max\{\ell\geq 0: (p_1\cdot\ldots \cdot p_n)^\ell \divides a \}\},\\
B^{(k)}:=&\{a/(p_1\cdot\ldots\cdot p_n)^k :a\in A^{(k)}  \}.
\end{align*}
For $1\leq i\leq n$, let
$$
B^{(k)}_i:=\{b\in B^{(k)} : p_i\ndivides b\}.
$$
By the induction hypothesis, each of the sets $B^{(k)}_i$ is finite. Therefore $B^{(k)}$ is finite because $B^{(k)}=\bigcup_{1\leq i\leq n}B_i^{(k)}$. It follows immediately that also
\begin{equation}\label{finiteC}
A^{(k)}\text{ is finite.}
\end{equation}
Suppose that
$|\{k\geq 0 : A^{(k)}\neq \emptyset\}|=\infty$. Choose $a=p_1^{\alpha_1}\cdot\ldots\cdot p_n^{\alpha_n}\in A$. Let $k_0 > \max\{\alpha_i : 1\leq i\leq n\}$ be such that $A^{(k_0)}\neq \emptyset$ and take $a'\in A^{(k_0)}$. Then $a\divides a'$, however $a\neq a'$, which yields a contradiction, i.e.\ we have
\begin{equation}\label{finiteD}
|\{k\geq 0 : A^{(k)}\neq \emptyset\}|<\infty.
\end{equation}
Since $A=\bigcup_{k\geq 0}A^{(k)}$, using~\eqref{finiteC} and~\eqref{finiteD}, we obtain $|A|<\infty$, and the result follows.
\end{proof}

\begin{Lemma}[see, e.g.,~\cite{MR0096565}]\label{kongruencje}
Let $b_1,\ldots,b_k\in\N$, $n_1\ldots,n_k\in\Z$. The system of congruences
$$
m\equiv n_i\bmod b_i,\ 1\leq i\leq k
$$
has a solution $m\in\N$ if and only if $\gcd(b_i,b_j)\divides (n_i-n_j)$ for any $i,j=1,...,k$.
\end{Lemma}

\begin{proof}[Proof of \cref{charakt}]
In view of Lemma~\ref{onedir}, we have (b) $\Rightarrow$ (a). We will now show (a) $\Rightarrow$ (b).
Assume that $A\subset\{1,\ldots,n\}$ satisfies condition~\eqref{Ther}, with $\{n_b : b\in\mathscr{B}\}$ as in the definition. Since $\mathscr{B}$ is primitive, it follows from Lemma~\ref{finite} that $\mathscr{B}^{(n)}$ is finite. Therefore, by Lemma~\ref{kongruencje}, there exists $m\in \N$ such that
$$
m\equiv -n_{b} \bmod b, \; b\in\mathscr{B}^{(n)}.
$$
It follows that
\begin{multline*}
\{m+1,\ldots,m+n\}\cap \mathcal{M}_{\mathscr{B}^{(n)}}\\
= (\{1,\ldots,n\}\cap \bigcup_{b\in\mathscr{B}^{(n)}}(b\Z +n_b))+m\subset (\{1,\ldots,n\}\setminus A)+m.
\end{multline*}
Applying \cref{tenwnio} to $\sA=\sB^{(n)}$, $k=m$, $C=\{1,\dots,n\}\cap \bigcup_{b\in\sB^{(n)}} (b\Z+n_b)$ and $C'=\{1,\ldots,n\}\setminus A$, we conclude that there exists $m'$ such that
$$
\{m'+1,\ldots,m'+n\}\cap \mathcal{M}_{\mathscr{B}}=
(\{1,\ldots,n\}\setminus A)+m'.
$$
Equivalently,
$$
\{m'+1,\ldots,m'+n\}\cap \mathcal{F}_{\mathscr{B}}= A+m',
$$
which yields (a) $\Rightarrow$ (b). In view of \cref{IU1}, this completes the proof.
\end{proof}

What remains to be proved is \cref{new1} and \cref{new1a}.
\begin{proof}[Proof of \cref{new1}]
For $u=1,\ldots,r$, let $j_u$ be such that $b_{j_u}\in\mathscr{A}$ and 
\begin{equation}\label{st0}
b_{j_u}\divides k+i_u.
\end{equation}
Let $B:=\mathscr{B}^{(n)}\cup\{b_{j_1},\ldots,b_{j_r}\}$. Then
\begin{equation}\label{st2}
\mbox{any $b\in \mathscr{B}\setminus B$ has a prime divisor $p>n$}
\end{equation}
and, by Lemma~\ref{finite},
\begin{equation}\label{st1} \mbox{$B$ is finite}.
\end{equation}
Let $\beta_1:=\lcm B$. Using~\eqref{st0} and the assumption \eqref{dziedz:as}, we obtain
\begin{multline*}
\{i_1,\dots,i_r\}\subset
(\{k+\beta_1\ell +1,\dots,k+\beta_1\ell +n\}\cap\mathcal{M}_{B})-(k+\beta_1\ell )\\
=\{k+1,\dots,k+n\}\cap \mathcal{M}_{B}-k
\subset(\{k+1, \dots,k+n\}\cap\mathcal{M}_{\mathscr{A}})-k=\{i_1, \dots,i_r\},
\end{multline*}
i.e.\ for any $\ell\in \Z$ we have
\begin{equation}\label{st3}
(\{k+\beta_1\ell+1,\ldots,k+\beta_1 \ell+n\}\cap\mathcal{M}_{B})-(k+\beta_1\ell)=\{i_1,\ldots,i_r\}.
\end{equation}
Using~\eqref{Au1}, we can find $j_0$ such that $\gcd(b_{j_0},\beta_1)=1$. It follows that there are $\ell_0\in\Z$ and $s\in\Z$ such that
$$
\beta_1\ell_0 - s b_{j_0} =-i_0-k.
$$
Hence, for $k':=k+\beta_1 \ell_0$, we have $b_{j_0}\divides k'+i_0$. Since $b_{j_0}\notin B$, we have $b_{j_0}>n$. It follows that
\begin{equation}\label{rowonly}
b_{j_0}\ndivides k'+i \text{ for any }1\leq i\neq i_0\leq n
\end{equation}
(indeed, if $b_{j_0}\divides k^{\prime}+i$, then $n<b_{j_0}\divides (i-i_0)$). Let $\beta:=\beta_1b_{j_0}$. It follows from~\eqref{rowonly} and~\eqref{st3} (with $l:=l_0+mb_{j_0}$) that
\begin{multline} \label{st5a}
(\{k'+\beta m+1,\ldots,k'+\beta m+n\}\cap\mathcal{M}_{B\cup\{b_{j_0}\}})-(k'+\beta m)\\
=\{i_0,i_1,\ldots,i_r\}
\end{multline}
for any $m\in\N$. Hence, it suffices to take $\mathscr{A}'=B\cup \{b_{j_0}\}$.
\end{proof}

The proof of \cref{new1a} will be proceeded by several lemmas.

\begin{Lemma}\label{prosciutto}
Let $\mathscr{R}$ be the intersection of finitely many arithmetic progressions with steps $d_1,\dots,d_r$. Then either $\mathscr{R}=\emptyset$ or $\mathscr{R}$ is equal to an arithmetic progression of step $\lcm(d_1,\dots,d_r)$.
\end{Lemma}
\begin{proof}
It suffices to notice that if $a\in \mathscr{R}$ then $\mathscr{R}=\lcm(d_1,\dots,d_r)\Z+a$.
\end{proof}

\begin{Lemma}\label{pas1}
Let $\beta, r,n \in \N$, and assume that $p>n$ is a prime that does not divide $\beta$. Assume that $\mathscr{R}$ is a union of finitely many arithmetic progressions with steps not divisible by $p$.
Then
\begin{equation}\label{np}
d\left(\left(\beta\Z+r\right)\cap\left(\bigcup_{i=1}^n\left( p\Z  -i\right)\right)\cap \mathscr{R}\right)=\frac{n}{p} d((\beta\Z+r)\cap \mathscr{R})
\end{equation}
and
\begin{equation}\label{1-np}
d\left(\left(\beta\Z+r\right)\setminus\left(\bigcup_{i=1}^n\left( p\Z  -i\right)\cup \mathscr{R}\right)\right)=\left(1-\frac{n}{p}\right) d\left(\left(\beta\Z+r\right)\setminus \mathscr{R}\right)
\end{equation}
\end{Lemma}

\begin{proof}
By Lemma~\ref{rozklad}, in order to prove~\eqref{np}, it suffices to prove it for $\mathscr{R}=b\Z +j$, where $p\ndivides b$. Moreover, since the progressions $ p\Z  -i$ are pairwise disjoint for $1\leq i\leq n$, what we need to show is
\begin{equation}\label{EEE1}
d((\beta\Z+r)\cap( p\Z  -i)\cap (b\Z +j))=\frac{1}{p} d((\beta\Z+r)\cap (b\Z +j))
\end{equation}
for each $1\leq i\leq n$. Clearly, the above equality holds if $(\beta\Z+r)\cap (b\Z +j)=\emptyset$. Otherwise, let $\beta':=\lcm(\beta,b)$ and take $a\in (\beta\Z+r)\cap (b\Z +j)$. Then, by \cref{prosciutto}, $(\beta\Z+r)\cap (b\Z +j)=\beta'\Z+a$ and \eqref{EEE1} is equivalent to
\begin{equation}\label{EEE2}
d((\beta'\Z+a)\cap (p\Z-i))=\frac{1}{p}d(\beta'\Z+a).
\end{equation}
Since $\gcd(\beta',p)=1$, it follows that $(\beta'\Z+a)\cap (p\Z-i)\neq \emptyset$ and \eqref{EEE2} is a straightforward consequence of \cref{prosciutto}.

In order to prove \eqref{1-np}, note that
\begin{align*}
& d\left(\left(\beta\Z+r\right)\setminus\left(\bigcup_{i=1}^n\left( p\Z  -i\right)\cup \mathscr{R}\right)\right)\\
 &=d\left(\beta\Z+r\right)-d\left(\left(\beta\Z+r\right)\cap \mathscr{R}\right) - d\left(\left(\beta\Z+r\right)\cap\left(\bigcup_{i=1}^n\left( p\Z  -i\right)\right)\right)\\
 &\quad+d\left(\left(\beta\Z+r\right)\cap\left(\bigcup_{i=1}^n\left( p\Z  -i\right)\cap \mathscr{R}\right)\right)\\
&=  d\left(\beta\Z+r\right)-d\left(\left(\beta\Z+r\right)\cap \mathscr{R}\right)- \frac{n}{p}d\left(\beta\Z+r\right)+\frac{n}{p}d\left(\left(\beta\Z+r\right)\cap \mathscr{R}\right)\\
&= \left(1-\frac{n}{p}\right)d((\beta\Z+r)\setminus \mathscr{R}),
\end{align*}
where the second equality follows from \eqref{np}.
\end{proof}

\begin{Lemma} \label{pas2} Let $\beta, r, n,c_1,...,c_m\in \N$.
Assume that $p>n$ is a prime,  $p$ divides $c_1,...,c_k$ and $p$ does not divide $c_{k+1},...,c_m$ nor $\beta$.
Then
\begin{multline}\label{nierpas2}
d\left(\left(\beta\Z+r\right)\cap\bigcap_{i=1}^n\left(\mathcal{F}_{\{c_1,...,c_m\}}-i\right)\right)\\
\ge \left(1-\frac{n}{p}\right)d\left(\left(\beta\Z+r\right)\cap\bigcap_{i=1}^n\left(\mathcal{F}_{\{c_{k+1},...,c_m\}}-i\right)\right).
\end{multline}
\end{Lemma}

\begin{proof}
Notice first that
\begin{equation}\label{prostee}
(A-i)^c=A^c-i\text{ for any }A\subset \Z, i\in\Z.
\end{equation}
Therefore,
\begin{equation}\label{FFF1}
\left(\beta\Z+r\right)\cap\bigcap_{i=1}^n\left(\mathcal{F}_{\{c_1,...,c_m\}}-i\right)=\left(\beta\Z+r\right)\setminus \left(\bigcup_{i=1}^n\left(\mathcal{M}_{\{c_1,...,c_m\}}-i\right)\right).
\end{equation}
Since
$$
\mathcal{M}_{\{c_1,...,c_m\}}\subset \mathcal{M}_{\{p, c_{k+1},...,c_m\}}= p\Z \cup \mathcal{M}_{\{c_{k+1},...,c_m\}},
$$
using~\eqref{FFF1}, we obtain
$$
\left(\beta\Z+r\right)\cap\bigcap_{i=1}^n\left(\mathcal{F}_{\{c_1,...,c_m\}}-i\right)\supset (\beta\Z+r) \setminus \left(\bigcup_{i=1}^n\left( p\Z  -i\right)\cup \bigcup_{i=1}^n\left(\mathcal{M}_{\{c_{k+1},...,c_m\}}-i\right)\right).
$$
To complete the proof, we apply Lemma~\ref{pas1} to $\mathscr{R}=\bigcup_{i=1}^n(\mathcal{M}_{\{c_{k+1},...,c_m\}}-i)$ and use again~\eqref{prostee}.
\end{proof}
\begin{Remark}\label{empy}
In the above lemma, we admit the situation when $k=m$ (we interpret $\{c_{k+1},\dots,c_m\}$ as the empty set and we have $\cf_{\emptyset}=\Z$ and $\cm_\emptyset=\emptyset$).
\end{Remark}

\begin{Lemma}\label{wniosekop}
Let $\beta, r, n\in\N$. Suppose that $\{c_m:m\geq 1\}\subset \N$ is Besicovitch.
Assume that $p>n$ is a prime,  $p$ divides $c_1$ but does not divide $\beta$.
Then the densities of $\left(\beta\Z+r\right)\cap\bigcap_{i=1}^n\left(\mathcal{F}_{\{c_m:m\ge 1\}}-i\right)$ and $\left(\beta\Z+r\right)\cap\bigcap_{i=1}^n\left(\mathcal{F}_{\{c_m:m\ge 2\}}-i\right)$ exist and
\begin{multline*}
d\left(\left(\beta\Z+r\right)\cap\bigcap_{i=1}^n\left(\mathcal{F}_{\{c_m:m\ge 1\}}-i\right)\right)\\
\ge  \left(1-\frac{n}{p}\right)d\left(\left(\beta\Z+r\right)\cap\bigcap_{i=1}^n\left(\mathcal{F}_{\{c_m:m\ge 2\}}-i\right)\right).
\end{multline*}
\end{Lemma}
\begin{proof}
Fix $M\in \N$ and assume that $c_{l_1},...,c_{l_t}$ are the elements of the set $\{c_1,...,c_M\}$ which are not divisible by $p$ ($t$ can be equal to 0, cf. \cref{empy}).
By Lemma \ref{pas2}, it follows that
\begin{multline*}
d\left(\left(\beta\Z+r\right)\cap\bigcap_{i=1}^n\left(\mathcal{F}_{\{c_1,...,c_M\}}-i\right)\right)\\
\ge \left(1-\frac{n}{p}\right)d\left(\left(\beta\Z+r\right)\cap\bigcap_{i=1}^n\left(\mathcal{F}_{\{c_{l_1},...,c_{l_t}\}}-i\right)\right).
\end{multline*}
On the other hand, $\mathcal{F}_{\{c_2,...,c_M\}}\subset \mathcal{F}_{\{c_{l_1},...,c_{l_t}\}}$. Thus, we obtain
\begin{multline*}
d\left(\left(\beta\Z+r\right)\cap\bigcap_{i=1}^n\left(\mathcal{F}_{\{c_1,...,c_M\}}-i\right)\right)\\
\ge \left(1-\frac{n}{p}\right)d\left(\left(\beta\Z+r\right)\cap\bigcap_{i=1}^n\left(\mathcal{F}_{\{c_{2},...,c_{M}\}}-i\right)\right).
\end{multline*}
In view of \cref{da-er}, we can pass to the limit with $M\rightarrow\infty$ and the assertion follows.
\end{proof}

\begin{Lemma}\label{wolnypas}
Suppose that $\mathscr{B}$ has light tails. Assume that $\beta,r,n\in\N$ and $b_{k_1},b_{k_2},...\in\mathscr{B}$ are such that each  $b_{k_j}$ has a prime divisor greater than $n$ and not dividing $\beta$.
Then the density of
$$
\left(\beta\Z+r\right)\cap \bigcap_{i=1}^n\left(\mathcal{F}_{\{b_{k_j}:j\ge 1\}}-i\right)
$$
exists and is positive.
\end{Lemma}
\begin{proof}
Observe that by \cref{wniosekop}, for any $m\ge 1$, we have
\begin{multline*}
d\left(\left(\beta\Z+r\right)\cap \bigcap_{i=1}^n\left(\mathcal{F}_{\{b_{k_m},b_{k_{m+1}},...\}}-i\right)\right)\\
\ge\left(1-\frac{n}{p}\right) d\left(\left(\beta\Z+r\right)\cap \bigcap_{i=1}^n\left(\mathcal{F}_{\{b_{k_{m+1}},...\}}-i\right)\right)
\end{multline*}
where $p>n$ is a prime divisor of $b_{k_m}$.
It follows that
\begin{multline*}
d\left(\left(\beta\Z+r\right)\cap \bigcap_{i=1}^n\left(\mathcal{F}_{\{b_{k_1},b_{k_{2}},...\}}-i\right)\right)\\
\ge \rho\left(m\right) d\left(\left(\beta\Z+r\right)\cap \bigcap_{i=1}^n\left(\mathcal{F}_{\{b_{k_m},b_{k_{m+1}},...\}}-i\right)\right),
\end{multline*}
where $\rho\left(m\right)>0$ depends only on $m$. Since $\mathscr{B}$ has light tails, for $m$ large enough so that $d\left(\mathcal{M}_{\{b_{k_m},b_{k_{m+1}},...\}}\right)<\frac{1}{n\beta}$, we have
$$
d\left(\left(\beta\Z+r\right)\cap \bigcap_{i=1}^n\left(\mathcal{F}_{\{b_{k_{m+1}},b_{k_{m+2}},...\}}-i\right)\right)>0
$$
and the assertion follows.
\end{proof}

\begin{proof}[Proof of \cref{new1a}]
For $u=1,\ldots,r$, let $j_u$ be such that $b_{j_u}\in\mathscr{A}$ and 
\begin{equation}\label{st01}
b_{j_u}\divides k+i_u.
\end{equation}
Without loss of generality, we may assume that $\mathscr{A}=\{b_{j_u} : 0\leq u\leq r\}\cup B^{(n)}$. Then, by \cref{finite}, $\sA$ is finite and we set $\beta:=\gcd(\sA)$. It follows by~\eqref{dziedz:as1} that
\begin{equation}\label{st5A}
(\{k+\beta m +1,\dots,k+\beta m+n\}\cap \mathcal{M}_{\mathscr{A}})-(k+\beta m)=\{i_0,\dots,i_r\}
\end{equation}
for any $m\in\N$. Let
$$
B:=\{b\in\mathscr{B}\setminus \mathscr{A} : \text{ all prime divisors of }b \text{ greater than $n$ divide }\beta\}
$$
($B$ may be empty) and notice that we have $B$ is finite. Indeed, if $p$ is a prime divisor of $b\in B$ then either $p\leq n$ or $p>n$ and divides $\beta$. Hence $|\text{Spec}(B)|<\infty$ and we can use Lemma~\ref{finite}. Since $\sB^{(n)}\subset \sA$, we have $B\subset \sB\setminus \sB^{(n)}$ and it follows that
\begin{equation}\label{st2aa}
\mbox{any $b\in B$ has a prime divisor $p>n$}.
\end{equation}
Let $b\in B$ and take a prime $p\divides b$, $p>n$ (such $p$ exists by~\eqref{st2aa}). By the definition of $B$, we have $p\divides \beta$, whence $p\divides  b_{j_u}$ for some $0\leq u\leq r$. It follows that if $b\divides k+\beta m+i$ for some $1\leq i\leq n$ then $i\in\{i_0,\ldots,i_r\}$ (otherwise, using~\eqref{st01}, we obtain $p\divides i_u-i$, which is impossible). Thus, by~\eqref{st5A}, we obtain
\begin{equation}\label{st6}
(\{k+\beta m+1,\ldots,k+\beta m+n\}\cap\mathcal{M}_{\mathscr{A}\cup B})-(k+\beta m)
=\{i_0,i_1,\ldots,i_r\}
\end{equation}
for any $m\in\N$. Let
$$
(\mathscr{B}\setminus \mathscr{A}) \setminus B=:B'=\{b_{k_1},b_{k_2},\dots\},
$$
i.e.\ each $b_{k_j}$ has a prime divisor greater than $n$, not dividing $\beta$. By \cref{wolnypas}, the density of the set
\begin{equation}\label{dodatniage}
(\Z \beta +k)\cap (\bigcap_{i=1}^n\mathcal{F}_{\{b_{k_j}:j\ge 1\}} -i)
\end{equation}
exists and is positive. Therefore, for $m\in\N$ from some positive density set, we have   $\beta m+k +i\in \mathcal{F}_{\{b_{k_j}:j\ge 1\}}$ for any $i=1,\ldots,n$.
Using~\eqref{st6}, it follows that for each such $m\in\N$, we have
\begin{multline*}
(\{k+\beta m+1,\ldots,k+\beta m+n\}\cap\mathcal{M}_{\mathscr{B}})-(k+\beta m)\\
=(\{k+\beta m+1,\ldots,k+\beta m+n\}\cap\mathcal{M}_{\sA \cup B})-(k+\beta m)=\{i_0,\dots,i_r\},
\end{multline*}
as required.
\end{proof}

\cref{OOF} is an immediate consequence of  \cref{quasi-gen} and of \cref{new1a} (applied to $\mathscr{A}:=\mathscr{B}$).

\section{Entropy}\label{s7}
\subsection{Entropy of $\widetilde{X}_\eta$ and $X_\mathscr{B}$ (proof of \cref{OOI})}\label{seOOI}
In this section our main goal is to prove \cref{OOI}.To fix attention, we will restrict ourselves to the case when $\mathscr{B}$ is infinite. The proof will be very similar to the proof of Theorem 5.3 in~\cite{Ab-Le-Ru}. However, since we dropped the assumptions~\eqref{settingerdosa}, we cannot use the Chinese Remainder Theorem directly and we will need an additional ingredient:
\begin{Lemma}[Rogers, see~\cite{MR687978}, page 242]\label{lerog}
For any $b_k$, $k\geq 1$, any $r_k\in \Z/b_k\Z$ and $K\geq 1$, we have
\begin{equation}\label{rog}
d\Big(\bigcup_{k\leq K}(b_k\Z+r_k) \Big)\geq d\left(\mathcal{M}_{\{b_1,\dots,b_K\}} \right).
\end{equation}
\end{Lemma}
\begin{Remark}\label{rerog}
Clearly, for any $n\in\N$,
$$
d\Big(\bigcup_{k\leq K}(b_k\Z+r_k) \Big)=\frac{1}{n \cdot b_1\cdot\ldots\cdot b_K}\Big| [1,n\cdot b_1\cdot\ldots\cdot b_K]\cap \Big(\bigcup_{k\leq K}(b_k\Z+r_k)\Big)  \Big|.
$$
\end{Remark}
\begin{proof}[Proof of \cref{OOI}]
In view of \cref{da-er}, the result will follow once we show
$$
h_{top}(S,\widetilde{X}_\eta)= h_{top}(S,X_\mathscr{B})= \overline{d}(\mathcal{F}_\mathscr{B}).
$$
For $n\in\N$ let
$$
\gamma(n):=|\{B\in \{0,1\}^n : B\text{ is }\mathscr{B}\text{-admissible}\}|
$$
and, for $K\geq 1$,
$$
\gamma_K(n):=|\{B\in \{0,1\}^n : B\text{ is }\{b_1,\dots,b_K\}\text{-admissible}\}|
$$
Clearly,
$$
\gamma(n)\leq \gamma_K(n) \text{ for any $K\geq 1$}.
$$
Moreover, any $\{b_1,\dots, b_K\}$-admissible $n\cdot b_1\cdot\ldots\cdot b_K$-block $B\in\{0,1\}^{[1,n\cdot b_1\cdot\ldots\cdot b_K]}$ can be obtained in the following way:
\begin{enumerate}[(a)]
\item choose $(r_1,\dots, r_K)\in \prod_{k\leq K}\Z/b_k\Z$ and set $B(j):=0$ for $1\leq j\leq n\cdot b_1\cdot\ldots\cdot b_K$ satisfying $j\equiv r_k \bmod b_k$ for some $1\leq k\leq K$,
\item complete the word by choosing arbitrarily $B(j)\in \{0,1\}$ for all other $1\leq j\leq n\cdot b_1\cdot\ldots\cdot b_K$.
\end{enumerate}
(Clearly, $(\text{supp }B) \cap (b_i\Z+r_i)=\emptyset$.) Notice that once $(r_1,\dots, r_K)\in \prod_{k\leq K}b_k\Z$ is fixed, the freedom in Step~(b) gives
$$
2^{n\cdot b_1\cdot\ldots\cdot b_K\left(1-d\left(\bigcup_{k\leq K}b_k\Z+r_k \right)\right)}
$$
pairwise distinct $\{b_1,\dots, b_K\}$-admissible $n\cdot b_1\cdot\ldots\cdot b_K$-blocks (cf.\ \cref{rerog}). Moreover, in view of Lemma~\ref{lerog}, this number does not exceed
\begin{equation}\label{eq:hh1}
2^{n\cdot b_1\cdot\ldots\cdot b_K(1-d_K)},
\end{equation}
where $d_K={d(\mathcal{M}_{\{b_1,\dots,b_K\}})}$

We will show that $h_{top}(S,X_\mathscr{B})\leq \overline{d}(\mathcal{F}_\mathscr{B})$. Fix $\vep>0$. In view of \cref{da-er}, if $K$ is large enough then $d_K\geq 1-\overline{d}(\mathcal{F}_\mathscr{B})-\vep$. Fix such $K$. It follows by Lemma~\ref{lerog}, \cref{rerog} and the discussion preceeding~\eqref{eq:hh1} that
$$
\gamma_{K}(n\cdot b_1\cdot\ldots\cdot b_K) \leq \prod_{k\leq K}b_k \cdot 2^{n\cdot b_1\cdot\ldots\cdot b_K \cdot  (1-d_K)},
$$
whenever $n=n(K,\vep)$ is sufficiently large. Thus (since the number of possible choices in Step~(a) equals $b_1\cdot \ldots\cdot b_K$), for such $n$, we obtain
$$
\gamma_{K}(n\cdot b_1\cdot\ldots\cdot b_K) \leq \prod_{k\leq K}b_k \cdot 2^{n\cdot b_1\cdot\ldots\cdot b_K \cdot  (\overline{d}(\mathcal{F}_\mathscr{B})+\vep)}.
$$
Therefore,
$$
h_{top}(S,X_\mathscr{B})=\lim_{n\to\infty} \frac{1}{n}\log \gamma(n)\leq \lim_{n\to\infty} \frac{1}{n}\log \gamma_K(n)\leq \overline{d}(\mathcal{F}_\mathscr{B}).
$$

We will now show that $h_{top}(S,\widetilde{X}_\eta)\geq \ov{d}(\mathcal{F}_\mathscr{B})$. For $n\geq 1$, denote by $p(n)$ the number of $n$-blocks occurring on $\widetilde{X}_\eta$. Let $(N_k)$ be such that
$$
\lim_{k\to\infty}\frac{1}{N_k}|[0,N_k]\cap \cf_\mathscr{B}|=\ov{d}(\mathcal{F}_\mathscr{B})
$$
(such a sequence exists by Theorem \ref{quasi-gen}).
Since
$$
p(N_k)\geq 2^{|[0,N_k]\cap \cf_\mathscr{B}|},
$$
it follows that
$$
h_{top}(S,\widetilde{X}_\eta)=\lim_{k\to\infty}\frac{1}{N_k}\log p(N_k) \geq\ov{d}(\mathcal{F}_\mathscr{B}).
$$
This concludes the proof.
\end{proof}

\begin{Remark}
Recall that a hereditary system has zero entropy if and only if $\delta_{(\dots,0,0,0,\dots)}$ is the unique invariant measure (for the proof, see \cite{MR3007694}). Therefore, since both, $\widetilde{X}_\eta$ and $X_\sB$, are hereditary, it follows by \cref{OOI} that the following conditions are equivalent:
\begin{itemize}
\item
$\mathcal{P}(S,X_\sB)=\{\delta_{(\dots,0,0,0,\dots)}\}$,
\item
$\mathcal{P}(S,\widetilde{X}_\eta)=\{\delta_{(\dots,0,0,0,\dots)}\}$,
\item 
$\bdelta(\cf_{\mathscr{B}})=0$.
\end{itemize}
In particular, this applies to $(S,X_\mathscr{P})$ (cf.\ \eqref{bf2}), even though $X_\mathscr{P}$ is uncountable, cf.\ \cref{nieprz}.
\end{Remark}\noindent

\subsection{Entropy of some invariant subsets of $\widetilde{X}_\eta$}
In this section we will prove the following:
\begin{Prop}\label{moreentropy}
If $\mathscr{B}$ is taut then
$$
h_{top}(S,Y_{\geq s_1,\geq s_2,\dots}\cap \widetilde{X}_\eta)<h_{top}(S,\widetilde{X}_\eta),
$$
whenever $s_k>1$ for some $k\geq 1$.
\end{Prop}
For this, we will need some tools.

\begin{Lemma}[cf.\ Lemma 1.17 in \cite{MR1414678} and \cref{da-er}]\label{hol:istnieje}
Let $\mathscr{B}\subset \N$. For any $q\in\N$ and $0\leq r\leq q-1$ the logarithmic density of $\mathcal{M}_\mathscr{B} \cup (q\Z+r)$ exists and
$$
\bdelta(\mathcal{M}_\mathscr{B} \cup (q\Z+r))=\un{d}(\mathcal{M}_\mathscr{B} \cup (q\Z+r))
=\lim_{k\to \infty}d(\mathcal{M}_{\{b_1,\dots, b_k\}}\cup (q\Z+r)).
$$
\end{Lemma}
\begin{proof}
Since
$$
\mathcal{M}_\mathscr{B} \cup (q\Z+r)=(q\Z+r)\cup \bigcup_{0\leq s\neq r\leq q-1}\mathcal{M}_\mathscr{B} \cap (q\Z+s),
$$
it suffices to prove that the logarithmic density of $\mathcal{M}_\mathscr{B} \cap (q\Z+s)$ exists and
\begin{equation}\label{i2}
\bdelta(\mathcal{M}_\mathscr{B} \cap (q\Z+s))=\un{d}(\mathcal{M}_\mathscr{B} \cap (q\Z+s))
=\lim_{k\to \infty}d(\mathcal{M}_{\{b_1,\dots, b_k\}}\cap (q\Z+s))
\end{equation}
for each $0\leq s\leq q-1$. Indeed, if~\eqref{i2} holds, we have
\begin{multline*}
\bdelta(\mathcal{M}_\mathscr{B}\cup (q\Z+r))\geq \un{d}(\mathcal{M}_\mathscr{B}\cup (q\Z+r))\\
\begin{aligned}
&\geq d(q\Z+r)+\bigcup_{0\leq s\neq r\leq q-1}\un{d}(\mathcal{M}_\mathscr{B} \cap (q\Z+s))\\
&=d(q\Z+r)+\bigcup_{0\leq s\neq r\leq q-1}\bdelta(\mathcal{M}_\mathscr{B}\cap (q\Z+s))=\bdelta(\mathcal{M}_\mathscr{B} \cup (q\Z+r)).
\end{aligned}
\end{multline*}
To show~\eqref{i2}, notice first that, for each $k\geq 1$, we have
$$
\un{d}(\mathcal{M}_\mathscr{B}\cap (q\Z+s)) \geq d(\mathcal{M}_{\{b_1,\dots,b_k\}}\cap (q\Z+s) ),
$$
whence
\begin{equation}\label{hol1}
\un{d}(\mathcal{M}_\mathscr{B}\cap (q\Z+s)) \geq \lim_{k\to \infty}d(\mathcal{M}_{\{b_1,\dots,b_k\}}\cap (q\Z+s) ).
\end{equation}
On the other hand, for each $k\geq 1$,
$$
\ov{\bdelta}(\mathcal{M}_\mathscr{B}\cap (q\Z+s))\leq d(\mathcal{M}_{\{b_1,\dots,b_k\}}\cap (q\Z+s))+{\bdelta}(\mathcal{M}_\mathscr{B} \setminus \mathcal{M}_{\{b_1,\dots,b_k\}}),
$$
whence, by \cref{da-er},
\begin{equation}\label{hol2}
\ov{\bdelta}(\mathcal{M}_\mathscr{B}\cap (q\Z+s)) \leq \lim_{k\to \infty}d(\mathcal{M}_{\{b_1,\dots,b_k\}}\cap (q\Z+s)).
\end{equation}
The claim follows from~\eqref{hol1} and~\eqref{hol2}.
\end{proof}

\begin{Lemma}\label{lm:help}
Assume that $\mathscr{B}$ is taut. Fix $k_0\geq 1$ and let $0<r<b_{k_0}$. Then
$$
\un{d}\left(\cm_\mathscr{B} \cup (b_{k_0}\Z+r)\right)> \un{d}\left(\cm_\mathscr{B}\right).
$$
\end{Lemma}
\begin{proof}
By \cref{hol:istnieje}, we have
\begin{equation}\label{Y1}
\un{d}(\cm_\mathscr{B}\cup(b_{k_0}\Z+r))=\bdelta(\cm_\mathscr{B}\cup(b_{k_0}\Z+r))=\bdelta(\cm_\mathscr{B})+\bdelta((b_{k_0}\Z+r)\setminus \cm_\mathscr{B}),
\end{equation}
where
\begin{equation}\label{Y2}
\bdelta((b_{k_0}\Z+r)\setminus \cm_\mathscr{B})=\bdelta((b_{k_0}\Z+r)\setminus \cm_{\mathscr{B}\setminus\{b_{k_0}\}}),
\end{equation}
since $(b_{k_0}\Z+r) \cap b_{k_0}\Z=\emptyset$. Moreover, since $(b_{k_0}\Z+r)\cup\cm_{\mathscr{B}\setminus\{b_{k_0}\}}$ is a disjoint union of $\cm_{\mathscr{B}\setminus\{b_{k_0}\}}$ and $(b_{k_0}\Z+r)\setminus \cm_{\mathscr{B}\setminus\{b_{k_0}\}}$ (and the logarithmic density of $(b_{k_0}\Z+r)\cup\cm_{\mathscr{B}\setminus\{b_{k_0}\}}$ and  $\cm_{\mathscr{B}\setminus\{b_{k_0}\}}$ exists by \cref{hol:istnieje} and \cref{da-er}, respectively), we obtain
\begin{equation}\label{Y3}
\bdelta((b_{k_0}\Z+r)\setminus \cm_{\mathscr{B}\setminus\{b_{k_0}\}})=\bdelta((b_{k_0}\Z+r)\cup\cm_{\mathscr{B}\setminus\{b_{k_0}\}})  -\bdelta(\cm_{\mathscr{B}\setminus\{b_{k_0}\}}).
\end{equation}
By the tautness of $\sB$, 
\begin{equation}\label{Y4}
\bdelta(\cm_\sB)>\bdelta(\cm_{\sB\setminus \{b_{k_0}\}}).
\end{equation}
Therefore, by~\eqref{Y1}, \eqref{Y2}, \eqref{Y3} and \eqref{Y4},
\begin{multline}\label{Y5}
\un{d}(\cm_\mathscr{B}\cup(b_{k_0}\Z+r))>\\
\bdelta(\cm_{\mathscr{B}\setminus\{b_{k_0}\}})+\bdelta((b_{k_0}\Z+r)\cup\cm_{\mathscr{B}\setminus\{b_{k_0}\}})  -\bdelta(\cm_{\mathscr{B}\setminus\{b_{k_0}\}})\\
=\bdelta((b_{k_0}\Z+r)\cup\cm_{\mathscr{B}\setminus\{b_{k_0}\}}).
\end{multline}
Moreover, applying consecutively \cref{hol:istnieje}, \cref{lerog} and \cref{da-er}, we obtain
\begin{align*}
\bdelta((b_{k_0}\Z+r)\cup \cm_{\sB\setminus \{b_{k_0}\}})&=\lim_{k\to\infty} d((b_{k_0}\Z+r)\cup \cm_{\{b_i : 1\leq i\leq k, i\neq k_0\}})\\
&\geq \lim_{k\to \infty}d (\cm_{\{b_i : 1\leq i\leq k\}})=\bdelta(\cm_\sB)=\un{d}(\cm_\sB).
\end{align*}
This, together with \eqref{Y5}, completes the proof.
\end{proof}

\begin{proof}[Proof of \cref{moreentropy}]
Fix $k_0\geq 1$ such that $s_{k_0}>1$. For $0<r<b_{k_0}$ let
$$
D_r:=\underline{d}\left(\cm_\mathscr{B} \cup (b_{k_0}\Z+r)\right)
$$
and $D:=\min_{0<r<b_{k_0}}D_r$. In view of Lemma~\ref{lm:help}, there exist $\vep>0$, $c>0$ such that
\begin{equation}\label{GH-2}
D-\un{d}(\cm_\mathscr{B})-2\vep>c>0.
\end{equation}
Let $K\geq k_0$ be large enough so that
\begin{equation}\label{GH-1}
d(\mathcal{M}_{\{b_1,\dots,b_K\}}\cup (b_{k_0}\Z+r)) \geq \un{d}(\mathcal{M}_\mathscr{B}\cup (b_{k_0}\Z+r))-\vep
\end{equation}
(such $K$ exists by \cref{hol:istnieje}). Finally, let $N_0\in\N$ be suffciently large, so that for $N>N_0$ we have
\begin{multline}
\begin{split}\label{GH0}
\frac{1}{N \cdot b_1\cdot\ldots \cdot b_K}\left|[0,N\cdot b_1\cdot\ldots\cdot b_K-1]\cap (\mathcal{M}_{\{b_1,\dots, b_K\}}\cup (b_{k_0}\Z+r)) \right|\\
\geq d(\mathcal{M}_{\{b_1,\dots, b_K\}}\cup (b_{k_0}\Z+r))-\vep.
\end{split}
\end{multline}
Fix $N> N_0$ and take $B$ which appears on $Y_{\geq s_1,\geq s_2,\dots}\cap \widetilde{X}_\eta$, with $|B|=N\cdot b_1\cdot\ldots\cdot b_K$. Then there exists $k\in\Z$ such that
\begin{equation}\label{GH1}
B+k\leq \eta[k,k+N\cdot b_1\cdot\ldots\cdot b_K-1].
\end{equation}
It follows by~\eqref{GH1} and by the choice of $k_0$ that there exists $0<r_0<b_{k_0}$ such that
\begin{equation}\label{GH2}
\text{supp }\eta \cap [k,k+N\cdot b_1\cdot\ldots\cdot b_K-1]\cap (b_{k_0}\Z+r_0)=\emptyset.
\end{equation}
Therefore, using \eqref{GH2}, \eqref{GH0}, \eqref{GH-1}, the definition of $D_{r_0}$ and $D$ and \eqref{GH-2}, we obtain
\begin{align*}
&\frac{|B|-|\text{supp }B|}{|B|}\\
&\geq \frac{1}{N\cdot b_1\cdot\ldots\cdot b_K}\left|[k,k+N\cdot b_1\cdot\ldots\cdot b_K-1]\cap (\cm_\mathscr{B}\cup (b_{k_0}\Z+r_0)) \right|\\
&\geq \frac{1}{N\cdot b_1\cdot\ldots\cdot b_K}\left|[k,k+N\cdot b_1\cdot\ldots\cdot b_K-1]\cap (\mathcal{M}_{\{b_1,\dots, b_K\}}\cup(b_{k_0}\Z+r_0)) \right|\\
&=\frac{1}{N\cdot b_1\cdot\ldots\cdot b_K}\left|[0,N\cdot b_1\cdot\ldots\cdot b_K-1]\cap (\mathcal{M}_{\{b_1,\dots, b_K\}}\cup(b_{k_0}\Z+r_0)) \right|\\
&\geq d(\mathcal{M}_{\{b_1,\dots, b_K\}}\cup (b_{k_0}\Z+r))-\vep\geq \un{d}(\mathcal{M}_\mathscr{B}\cup (b_{k_0}\Z+r_0))-2\vep\\
&= D_{r_0}-2\vep\geq D-2\vep>\un{d}(\cm_\mathscr{B})+c.
\end{align*}
Thus
\begin{equation}\label{eq:nosnik}
\frac{|\text{supp }B|}{|B|}< \ov{d}(\cf_\mathscr{B})-c.
\end{equation}
We will now proceed as in the proof of \cref{OOI}. For $n\in\N$, let
$$
\gamma^{s_1,s_2,\dots}(n):=|\{B\in \{0,1\}^n : B\text{ appears on }Y_{\geq s_1,\geq s_2,\dots}\cap \widetilde{X}_\eta\}|
$$
and, for $K\geq 1$,
$$
\gamma^{s_1,s_2,\dots,s_K}_K(n):=|\{B\in \{0,1\}^n : |\text{supp }B|\leq b_k-s_k\text{ for } 1\leq k\leq K\}|.
$$
Clearly,
$$
\gamma^{s_1,s_2,\dots}(n)\leq \gamma^{s_1,s_2,\dots,s_K}_K(n) \text{ for any $K\geq 1$}.
$$
Consider the following procedure of defining a block $B\in \{0,1\}^n$:
\begin{enumerate}[(a)]
\item choose $(r_1,\dots, r_K)\in \prod_{k\leq K}\Z/b_k\Z$, set $B(j):=0$ for $1\leq j\leq n$ such that $j\equiv r_k \bmod b_k$ for some $1\leq k\leq K$; choose $r_{k_0}'\not\equiv r_{k_0}\bmod b_{k_0}$ and set $B(j):=0$ for $1\leq j\leq n$ such that $j\equiv r_{k_0}\bmod b_{k_0}$,
\item complete the block by choosing arbitrarily $B(j)\in \{0,1\}$ for all other~$1\leq j\leq n$.
\end{enumerate}
Notice that all $B\in \{0,1\}^n$ satisfying
\begin{equation}\label{bloczkii}
|(\text{supp }B) \bmod b_k|\leq\begin{cases}
b_k-1,& \text{ for }k\neq k_0,\\
b_k-2,& \text{ for }k=k_0,
\end{cases}
\end{equation}
can be obtained this way. In particular, we obtain all blocks $B\in \{0,1\}^n$ such that 
$$
|(\text{supp }B) \bmod b_k|\leq b_k-s_k\text{ for }k\geq 1.
$$ 
Notice also that once the parameters $(r_1,\dots, r_K)$ and $r_{k_0}'$ in Step~(a) are fixed, the freedom in Step~(b) gives, for $n=N\cdot b_1\cdot\ldots\cdot b_K$, in view of~\eqref{eq:nosnik}, at most
$$
2^{N\cdot b_1\cdot\ldots\cdot b_K (\ov{d}(\cf_\mathscr{B})-c)}
$$
$N\cdot b_1\cdot\ldots\cdot b_K$-blocks. It follows that
\begin{multline*}
h_{top}(S,Y_{\geq s_1,\geq s_2,\dots}\cap \widetilde{X}_\eta)=\lim_{n\to\infty}\frac{1}{n}\log \gamma^{s_1,s_2,\dots}(n)\\
\leq \lim_{n\to \infty}\frac{1}{n}\log \gamma_K^{s_1,s_2,\dots,s_K}(n)\leq \ov{d}(\cf_\mathscr{B})-c=h_{top}(S,\widetilde{X}_\eta)-c,
\end{multline*}
which completes the proof.
\end{proof}

\begin{Cor}\label{gdziemax}
Suppose that $\mathscr{B}\subset \N$ is taut. Let $\nu\in \mathcal{P}(S,\widetilde{X}_\eta)$ be such that $h(S,\widetilde{X}_\eta,\nu)=h_{top}(S,\widetilde{X}_\eta)$. Then $\nu(Y\cap \widetilde{X}_\eta)=1$.
\end{Cor}
\begin{proof}
By considering the ergodic decomposition, we may restrict ourselves to $\nu\in\mathcal{P}^e(S,\widetilde{X}_\eta)$. Fix such $\nu$ and suppose that $h(S,\widetilde{X}_\eta,\nu)=h_{top}(S,\widetilde{X}_\eta)$ but $\nu(Y\cap \widetilde{X}_\eta)=0$ (by the ergodicity of $\nu$, we have $\nu(Y\cap \widetilde{X}_\eta)\in\{0,1\}$). Note that, for each $k\geq 1$, there exists $1\leq s_k<b_k$ such that
$
\nu(Y_{s_k}^k\cap \widetilde{X}_\eta)=1,
$
i.e., we obtain $(s_k)_{k\geq 1}$ such that
$
\nu(Y_{s_1,s_2,\dots}\cap \widetilde{X}_\eta)=1,
$
so, all the more,
$
\nu(Y_{\geq s_1,\geq s_2,\dots}\cap \widetilde{X}_\eta)=1.
$
Since $\nu(Y\cap\widetilde{X}_\eta)=0$, there exists $k\geq 1$ such that $s_k\geq 2$. But then, by \cref{moreentropy} and the variational principle,
$
h(S,\widetilde{X}_\eta,\nu)=h(S,Y_{\geq s_1,\geq s_2,\dots}\cap \widetilde{X}_\eta,\nu)\leq h_{top}(S,Y_{\geq s_1,\geq s_2,\dots}\cap \widetilde{X}_\eta)<h_{top}(S,\widetilde{X}_\eta).
$
This contradicts our assumption and we conclude.
\end{proof}

\section{Tautness and support of $\nu_\eta$ (proof of \cref{tautychar})}\label{s8}

We will now use \cref{TTC} and \cref{OOE} to prove \cref{tautychar}.
\begin{proof}[Proof of \cref{tautychar}]
Notice first that \eqref{tautycharA} $\Rightarrow $ \eqref{tautycharB} is an immediate consequence of \cref{gdziemax}. Now, we will show that also \eqref{tautycharB} $\Rightarrow$ \eqref{tautycharC} holds. We claim that
\begin{equation}\label{nm:1}
\nu_\eta(\varphi(\theta(Y\cap \widetilde{X}_\eta)))=1.
\end{equation}
Then, since by \cref{1.2.6} we have $\varphi(\theta(Y\cap \widetilde{X}_\eta))\subset Y$, it will follow that $\nu_\eta(Y)=1$. Moreover, since, by \cref{OOE}, we have $\nu_\eta(X_\eta)=1$, we obtain \eqref{tautycharC}. Thus, we are left to prove \eqref{nm:1}. Recall that by \cref{1.2.6}, we have $\theta_\ast(\nu)=\PP$.
Therefore,
\begin{multline*}
\nu_\eta(\varphi(\theta(Y\cap \widetilde{X}_\eta)))=\PP(\varphi^{-1}(\varphi(\theta(Y\cap \widetilde{X}_\eta))))\geq \PP(\theta(Y\cap \widetilde{X}_\eta))\\
=\theta_\ast\nu(\theta(Y\cap \widetilde{X}_\eta))=\nu(\theta^{-1}(\theta(Y\cap \widetilde{X}_\eta)))\geq \nu(Y\cap \widetilde{X}_\eta)=1,
\end{multline*}
i.e.\ \eqref{nm:1} indeed holds.

It remains to show that \eqref{tautycharC} implies \eqref{tautycharA}. Suppose that $\mathscr{B}$ is not taut. Let $\sB'$ be as in the proof of \cref{tautywszystko}. For simplicity, we assume that $\sB'$ is given by \eqref{KKindu6}, i.e.
$$
\sB'=( \mathscr{B}\setminus \bigcup_{n\geq 1}c_n\Z ) \cup \{c_n: n\geq 1\}=( \mathscr{B}\setminus \bigcup_{n\geq 1}c_n\mathscr{A}^n) \cup \{c_n: n\geq 1\},
$$
where $\sA^n$, $n\geq 1$, are Behrend sets. By \cref{tautywszystko}, $\sB'$ is taut and we have $\nu_\eta=\nu_{\eta'}$. Let
$$
Y':=\{x\in \{0,1\}^\Z : |\text{supp }x\bmod b_k'|=b_k'-1 \text{ for each }k\geq 1\}.
$$
By the first part of the proof, we have $\nu_{\eta'}(Y'\cap {X}_{\eta'})=1$. We will show that $\nu_\eta(Y\cap X_\eta)=0$. Since $\nu_\eta=\nu_{\eta'}$, it suffices to show that  $Y\cap Y'=\emptyset$. Take $a\geq 2$ such that $c_1 a \in \mathscr{B}$ and $c_1\in\mathscr{B}'$ and consider the natural projections
$$
\Z\xrightarrow{\pi_1} \Z/c_1a\Z \xrightarrow{\pi_2} \Z/c_1\Z
$$
($\pi_1(n)=n\bmod c_1a$ for $n\in\Z$ and $\pi_2(n)=n\bmod c_1$ for $n\in\Z/c_1a\Z$). Then, for any $A\subset \Z$, we have
$$
\pi_1(A) \subset \pi_2^{-1}(\pi_2(\pi_1(A))).
$$
Moreover, for any $B\subset \Z/c_1\Z$, we have $|\pi_2^{-1}(B)|=a |B|$. Therefore, for $x\in \{0,1\}^\Z$, we have
\begin{multline*}
|\text{supp }x \bmod c_1 a|=|\pi_1(\text{supp }x)|\leq |\pi_2^{-1}(\pi_2(\pi_1(\text{supp }x)))|\\
=a |\pi_2(\pi_1(\text{supp }x))|=a |\text{supp }x \bmod c_1|.
\end{multline*}
Therefore,
\begin{align*}
Y'&\subset \{x\in \{0,1\}^\Z : |\text{supp }x \bmod c_1|=c_1-1\}\\
&\subset \{x\in \{0,1\}^\Z : |\text{supp }x \bmod c_1 a|\leq c_1 a -a\}
\end{align*}
and, on the other hand, we have
$$
Y\subset \{x\in \{0,1\}^\Z : |\text{supp }x \bmod c_1 a|=c_1 a-1\}.
$$
Since $c_1 a -a<c_1 a-1$, we conclude that indeed $Y\cap Y'=\emptyset$. This completes the proof.
\end{proof}
\begin{Remark}
If $\mathscr{B}\subset \N$ has light tails then $\nu_\eta(Y\cap X_\eta)=1$ can be showed directly. Namely, fix $K\geq 1$ and let
$$
Y_K:=\{x\in \{0,1\}^\Z : |\text{supp }x\bmod b_k|=b_k-1 \text{ for }1\leq k\leq K \}.
$$
Then:
\begin{itemize}
\item $S(Y_K\cap X_\eta)=Y_K\cap X_\eta$,
\item $\eta\in Y_K\cap X_\eta$ (by \cref{gdzieeta}), in particular, $Y_K\cap X_\eta\neq\emptyset$,
\item $Y_K\cap X_\eta$ is open in $X_\eta$ (indeed, if $x\in Y_K\cap X_\eta$ and $M\in\N$ is such that $\text{supp }x \bmod b_k = (\text{supp }x \cap [0,M]) \bmod b_k$ for each $k\geq1$ then for each $y\in X_\eta$ with $y[0,M]=x[0,M]$, we have $y\in Y_K\cap X_\eta$).
\end{itemize}
In view of \cref{OOF}, since $Y_K\cap X_\eta$ is open and non-empty, we have $\nu_\eta(Y_K\cap X_\eta)>0$. By ergodicity and $S$-invariance, we obtain $\nu_\eta(Y_K\cap X_\eta)=1$. It follows that $\nu_\eta(Y\cap X_\eta)=\nu_\eta(\bigcap_{K\geq1}Y_K\cap X_\eta)=1$.

\end{Remark}

\section{Intrinsic ergodicity: taut case (\cref{OOH} -- first steps)}\label{s9}
Recall the following result:
\begin{Th}[\cite{MR3356811}]\label{interg}
If $\mathscr{B}\subset \N$ is infinite, coprime then $(S,\widetilde{X}_\eta)$ is intrinsically ergodic (in fact, $X_\eta=\widetilde{X}_\eta$).
\end{Th}
In this section we will extend \cref{interg} to the case when $\mathscr{B}$ is taut. The main ideas come from~\cite{MR3356811}. We will present the sketch of the proof only, referring the reader to~\cite{MR3356811} for the remaining details.\footnote{Another proof of \cref{jakzbenjim} will be presented in \cref{serevisited}.}

\begin{Th}\label{jakzbenjim}
Let $\mathscr{B}\subset\N$ and suppose that $\mathscr{B}$ is taut. Then $(S,\widetilde{X}_\eta)$ is intrinsically ergodic. In particular, if $X_\eta=\widetilde{X}_{\eta}$,\footnote{E.g. when $\mathscr{B}$ has light tails and satisfies \eqref{Au1}, see \cref{charakt}.} then $(S,X_\eta)$ is intrinsically ergodic.
\end{Th}

\begin{Remark}\label{RR:1}
If $\mathscr{B}\subset \N$ is finite, even though $X_\eta \subsetneq \widetilde{X}_\eta$, the subshift $(S,X_\eta)$ is intrinsically ergodic. Indeed, in view of \cref{B-skonczony}, $X_\eta$ is finite, with $|X_\eta|=\lcm(\mathscr{B})$ and $(S,X_\eta)$ is nothing but the rotation on $\lcm(\mathscr{B})$ points. It is uniquely ergodic, so, in particular, intrinsically ergodic.
\end{Remark}

\begin{proof}[Sketch of the proof of \cref{jakzbenjim}]
We will only consider the case when $\sB$ is infinite. Let $\nu$ be a measure of maximal entropy for $(S,\widetilde{X}_\eta)$. Then, by \cref{gdziemax}, $\nu(Y\cap \widetilde{X}_\eta)=1$. What we need to show is that the conditional measures $\nu_ g $ in the disintegration
$$
\nu=\int_ G  \nu_ g \ d\PP( g )
$$
(cf.\ \cref{1.2.6}) of $\nu$ over $\PP$ given by the mapping $\theta\colon Y\cap \widetilde{X}_\eta\to  G $ are unique $\PP$-a.e. In order to do it, we will show that for $A$ from some countable dense family of measurable subsets of $\widetilde{X}_\eta$,
\begin{equation}\label{doesnotdep}
\nu_ g (A) \text{ does not depend on }\nu,\text{ for }\PP\text{-a.e. } g \in G .
\end{equation}
\paragraph{Step 1.}
Let $Q=(Q_0,Q_1)$ be the partition of $Y\cap \widetilde{X}_\eta$ according to the value at the zero coordinate, i.e.\ $Q_i=\{y\in Y\cap \widetilde{X}_\eta : y(0)=i\}$, $i=0,1$ (this is a generating partition). Let
$$
Q^-:=\bigvee_{j\leq -1}S^jQ \text{ and }\mathcal{A}:=\theta^{-1}(\mathcal{B}( G )).
$$
Then, for $m\geq 0$, one can show that we have the following commuting diagram:
\begin{center}
\begin{tikzpicture}
	\node (U)  at (0,3) {$(Y\cap \widetilde{X}_\eta,\nu)$};
	\node (M) at (0,1.5) {$((Y\cap \widetilde{X}_\eta)/{S^{-m}Q^-},\ov\nu_m)$};
	\node (B) at (0,0) {$( G ,\PP)$};
	\draw[->] (U) edge node[auto] {$\pi_m$} (M)
		       (M) edge node[auto] {$\rho_m$} (B)
		       (U) edge[loop above] node[above] {$S$} (U)
		       (M) edge[loop right] node[left] {$S$} (M)
		       (B) edge[loop below] node {$T$} (B);
	\draw[<-] (U) to ++(3,0)   to node[right, midway]{$\varphi$}++(0,-3)  to(B) ;	
	\draw[->] (U) to ++(-3,0) to node[left,midway] {$\theta$} ++ (0,-3) to (B);	
\end{tikzpicture}
\end{center}
where $\pi_m\colon Y\cap \widetilde{X}_\eta\to ( Y\cap \widetilde{X}_\eta)/S^{-m}Q^-$ and $\rho_m\colon (Y\cap \widetilde{X}_\eta)/S^{-m}Q^-\to G $ are the natural quotient maps, $\overline{\nu}_m:=(\pi_m)_\ast(\nu)$, and $(\rho_m)_\ast(\overline{\nu}_m)=\PP$. In this diagram, $\theta$ is measure-preserving, while  $\varphi\colon  G \to Y\cap \widetilde{X}_\eta$ is defined $\PP$-a.e.\ and is not measure-preserving (notice that by \cref{tautychar}, we can treat $\varphi$ as a map with codomain $Y\cap X_\eta\subset Y\cap \widetilde{X}_\eta$).
\paragraph{Step 2.} Fix $m\geq 0$ and let, for $j=0,1$:
\begin{align*}
&C_m^j:=S^{-m}Q_j=\{x\in Y\cap \widetilde{X}_\eta : x(m)=j\},\\
&\widehat{C}_m^j:=\varphi^{-1}(C_m^j)=\{ g \in G  : \varphi( g )(m)=1\},\\
&B_m^j:=\rho_m^{-1}(\widehat{C}_m^j).
\end{align*}
This gives us the following diagram:
{\footnotesize
\begin{center}
\begin{tikzpicture}
	\node (U) at (0.3,3) {$Y\cap \widetilde{X}_\eta$};
	\node (Ur) at (2.5,3) {$=\theta^{-1}(\widehat{C}^0_{m})\cup \theta^{-1}(\widehat{C}^1_{m})$};
	\node (Ul) at (-1.7,2.98) {$S^{-m}Q_0\cup S^{-m}Q_1=$};
	\node (M) at (0.3,1.5) {$(Y\cap \widetilde{X}_\eta)/{S^{-m}Q^-}$};
	\node (Mr) at (4.2,1.5) {$=\rho^{-1}_{m}(\widehat{C}^0_{m})\cup\rho^{-1}_m(\widehat{C}^1_{m})=B^0_{m}\cup B^1_{m}$};
	\node (B) at (0.3,0) {$ G $};
	\node (Br) at (1.35,0) {$=\widehat{C}^0_{m} \cup \widehat{C}^1_{m}$};
	\draw[->] (U) edge node[auto] {$\pi_m$} (M)
		       (M) edge node[auto] {$\rho_m$} (B);
	\draw[<-] (Ur) to ++(4.5,0)   to node[right, midway]{$\varphi$}++(0,-3)  to(Br) ;	
	\draw[->] (Ul) to ++(-1.8,0) to node[left,midway] {$\theta$} ++ (0,-3) to (B);
\end{tikzpicture}
\end{center}
}

\paragraph{Step 3.}
Using $\theta^{-1}(\widehat{C}_m^0)\subset S^{-m}Q_0$, one can show that for each $\overline{y}\in B_m^0$,
$$
(\ov{\nu}_m(S^{-m}Q_0|S^{-m}Q^-),\ov{\nu}_m(S^{-m}Q_1|S^{-m}Q^-))=(1,0)=:(\lambda_0(Q_0),\lambda_0(Q_1)),
$$
whence
$$
H_\nu(S^{-m}Q|S^{-m}Q^-)(\ov{y})=0\text{ whenever }\ov{y}\in B_m^0.
$$
Therefore, using \cref{OOI} and \cref{quasi-gen}, we obtain
\begin{align*}
\bdelta&(\cf_\mathscr{B})=h_{top}(S,\widetilde{X}_\eta)=h(S,\widetilde{X}_\eta,\nu)
=\int_{Y/S^{-m}Q^-}H_\nu(S^{-m}Q | S^{-m}Q^-)\ d\ov{\nu}_m\\
&=\int_{B_m^1}H_\nu(S^{-m}Q | S^{-m}Q^-)\ d\ov{\nu}_m\leq \ov{\nu}_m(B_m^1)=(\rho_m)_\ast(\ov{\nu}_m)(\widehat{C}_m^1)\\
&=\PP(\widehat{C}_m^1)=\varphi_\ast\PP(C_m^1)=\nu_\eta(C_m^1)=\bdelta(\cf_\mathscr{B}).
\end{align*}
It follows that for $\ov{\nu}_m$-a.e.\ $\ov{y}\in B_m^1$,
\begin{multline}
(\ov{\nu}_m(S^{-m}Q_0|S^{-m}Q^-),\ov{\nu}_m(S^{-m}Q_1|S^{-m}Q^-))\\
=(1/2,1/2)=:(\lambda_1(Q_0),\lambda_1(Q_1)).
\end{multline}
\paragraph{Step 4.}
In view of Step 3., for $\ov{\nu}_m$-a.e.\ $\ov{y}$, we have
$$
\ov{\nu}_m(S^{-m}Q_{i_{m-r}}|S^{-m}Q^-)(\ov{y}i_{-m}\ldots i_{m-r-1})=\lambda_{j_r}(Q_{i_{m-r}}),
$$
where $j_r=\varphi(\rho_m(\ov{y}i_{-m}\ldots i_{m-r-1}))(m)=\varphi(\rho_m(\ov{y}))(m+r)$. Therefore, using the chain rule for conditional probabilities, one can show
\begin{multline}\label{thisformula}
\ov{\nu}_m(S^mQ_{i_m}\cap\ldots\cap Q_{i_0}\cap S^{-1}Q_{i_{-1}}\cap\ldots\cap S^{-m}Q_{i_{-m}}|S^{-m}Q^-)(\ov{y})\\
=\prod_{r=0}^{2m}\ov{\nu}_m(S^{-m}Q_{i_{m-r}}|S^{-m}Q^-)(\ov{y}i_{-m}\ldots i_{m-r-1})=\prod_{r=0}^{2m}\lambda_{j_r}(Q_{i_{m-r}}).
\end{multline}
For $A\in \bigvee_{t=-m}^{m}S^tQ$,
$$
\nu_ g (A)=\mathbb{E}^\nu(A| G )( g )=\mathbb{E}^\nu(\mathbb{E}^\nu(A|Y/S^{-m}Q^-)(\ov{y}_m)| G )( g ).
$$
Since~\eqref{thisformula} does not depend on $\ov{y}$ itself, but only on the values $\varphi(\rho_m(\ov{y}))(m+r)$, we obtain~\eqref{doesnotdep} for $A\in \bigvee_{t=-m}^{m}S^tQ$. The proof is complete as $m\geq 0$ was arbitrary.
\end{proof}

\section{Invariant measures (proof of \cref{OOG})}\label{s10}
In~\cite{MR3356811}, a description of $\mathcal{P}(S,X_\eta)$ was given in case of $\mathscr{B}$ infintie, coprime and thin (recall that in this case we have $X_\eta=\widetilde{X}_\eta$). Here, we extend this result by proving \cref{OOG}, which yields a description of $\mathcal{P}(S,\widetilde{X}_\eta)$ for all $\mathscr{B}$ (in particular, when $X_\eta=\widetilde{X}_\eta$, we obtain a description of $\mathcal{P}(S,X_\eta)$). 
\begin{Remark}
Notice that \cref{OOG} result is stated in a different, more compact form than in~\cite{MR3356811}. What coresponds directly to~\cite{MR3356811} is \cref{twY} in \cref{se:9.1} and \cref{miaryall} in \cref{se:9.3}. Notice that \cref{miaryall} is an immediate consequence of \cref{OOG} (it suffices to take $b_k'=b_k$ for all $k\geq 1$). The role of $b_k' \divides b_k$, $k\geq 1$ will become more clear later when we discuss the discrete rational part of the spectrum of $(S,\widetilde{X}_\eta,\nu)$, see \cref{sewymierne}.
\end{Remark}
We will present only sketches of the proofs, referring the reader to~\cite{MR3356811} for the remaining details (which can be repeated word by word). For the sake of simplicity, we will restrict to the case when $\mathscr{B}$ is infinite. For finite $\mathscr{B}$ the proofs go along the same lines (and are sometimes simpler).

\subsection{Invariant measures on $Y\cap \widetilde{X}_\eta$ (\cref{OOG} -- first steps)}\label{se:9.1}

\begin{Th}\label{twY}
For any $\nu\in\mathcal{P}^e(S,Y\cap \widetilde{X}_\eta)$, there exists $\widetilde{\rho}\in\mathcal{P}^e(S\times S,X_\eta\times \{0,1\}^\Z)$ such that $\widetilde{\rho}|_{X_\eta}=\nu_{\eta}$ and $M_\ast(\widetilde{\rho})=\nu$, where $M\colon X_\eta\times \{0,1\}^\Z\to \widetilde{X}_\eta$ stands for the coordinatewise multiplication.
\end{Th}
\begin{Remark}
Notice that in order to prove \cref{twY}, it suffices to find $\widetilde{\rho}\in\mathcal{P}(S\times S,X_\eta\times \{0,1\}^\Z)$ such that $\widetilde{\rho}|_{X_\eta}=\nu_{\eta}$ and $M_\ast(\widetilde{\rho})=\nu$ and use the ergodic decomposition.
\end{Remark}
We will first present the outline of the proof. Notice that if $\nu\in\mathcal{P}^e(S,Y\cap \widetilde{X}_\eta)$ then $\nu\neq \delta_{(\dots,0,0,0,\dots)}$.
\paragraph{Step 1.} We define $\widetilde{T}\colon  G \times\{0,1\}^\Z \to G \times\{0,1\}^\Z$ by
$$
\widetilde{T}( g ,x)=\begin{cases}
(T g ,x),& \text{ if }\varphi( g )(0)=0,\\
(T g ,Sx),& \text{ if }\varphi( g )(0)=1.
\end{cases}
$$
We will define $\Theta \colon Y\cap \widetilde{X}_\eta \to  G \times\{0,1\}^\Z$ a.e.\ with respect to any $\nu\in\mathcal{P}^e(S,Y\cap \widetilde{X}_\eta)$, and $\Phi\colon  G \times\{0,1\}^\Z\to X_\mathscr{B}$ a.e.\ with respect to any $\widetilde{T}$-invariant measure, so that the following diagram commutes:
\begin{equation*}\label{diagram1}
\begin{tikzpicture}[baseline=(current  bounding  box.center)]
\node (DL) at (0,0) {$ G \times\{0,1\}^{\Z}$};
\node (UL) at (0,1.5) {$Y\cap \widetilde{X}_\eta$};
\node (DR) at (3,0) {$ G \times\{0,1\}^{\Z}$.};
\node (UR) at (3,1.5) {$Y\cap \widetilde{X}_\eta$};
\node (DDL) at (0,-1.5) {$X_\mathscr{B}$};
\node (DDR) at (3,-1.5) {$X_\mathscr{B}$};
\draw[->] (DL) edge node[auto] {$\widetilde{T}$} (DR)
	       (UL) edge node[auto] {$S$} (UR)
	       (UL) edge node[auto] {$\Theta$} (DL)
	       (UR) edge node[auto] {$\Theta$} (DR)
	       (DL) edge node[auto] {$\Phi$} (DDL)
	       (DR) edge node[auto] {$\Phi$} (DDR)
	       (DDL) edge node[auto] {$S$} (DDR);
\end{tikzpicture}
\end{equation*}
Then, we will prove that
\begin{equation}\label{isidentity}
\Phi\circ \Theta=id \text{ a.e.\ with respect to any }\nu\in\mathcal{P}^e(S,Y\cap \widetilde{X}_\eta).
\end{equation}
This will give, for any $\nu\in\mathcal{P}^e(S,Y\cap \widetilde{X}_\eta)$, the equality
$$
\nu=\Phi_\ast\Theta_\ast\nu,\text{ with }\Theta_\ast\nu\in\mathcal{P}^e(\widetilde{T}, G \times \{0,1\}^\Z).
$$
\paragraph{Step 2.}
We will define $\Psi\colon  G \times\{0,1\}^\Z\to G \times\{0,1\}^\Z$ a.e.\ with respect to any $T\times S$-invariant measure, so that $\Psi$ is onto a.e.\ with respect to any $\widetilde{T}$-invariant measure, and so that the following diagram commutes:
\begin{equation}\label{diagram3}
\begin{tikzpicture}[baseline=(current  bounding  box.center)]
\node (DL) at (0,0) {$ G \times\{0,1\}^{\Z}$};
\node (UL) at (0,1.5) {$ G \times\{0,1\}^{\Z}$};
\node (DR) at (3,0) {$ G \times\{0,1\}^{\Z}$};
\node (UR) at (3,1.5) {$ G \times\{0,1\}^{\Z}$};
\draw[->] (DL) edge node[auto] {$\widetilde{T}$} (DR)
	       (UL) edge node[auto] {$T\times S$} (UR)
	       (UL) edge node[auto] {$\Psi$} (DL)
	       (UR) edge node[auto] {$\Psi$} (DR);
\end{tikzpicture}
\end{equation}
In fact, $\Psi$ will be defined on $G_0\times \{0,1\}^\Z$, where $G_0\subset G$ and $\PP(G_0)=1$.

We will define a system of measures $\{ \lambda_{(g,y)} : (g,y)\in G_0\times \{0,1\}^\Z\}$, where $\lambda_{(g,y)}(\Psi^{-1}( g ,y))=1$ for $(g,y)\in G_0\times \{0,1\}^\Z$  and such that
\begin{enumerate}[(a)]
\item\label{eq:Cz}
the map $F\colon ( g ,y)\mapsto \lambda_{( g ,y)}$ is measurable,
\item\label{eq:C}
$(T\times S)_\ast \lambda_{( g ,y)}=\lambda_{\widetilde{T}( g ,y)}$.
\end{enumerate}
Then for any $\rho\in \mathcal{P}^e(\widetilde{T}, G \times \{0,1\}^\Z)$, we will obtain
$$
\widetilde{\rho}:=\int \lambda_{( g ,y)}\ d\rho( g ,y)\in \mathcal{P}(T\times S, G \times \{0,1\}^\Z)\text{ with }\Psi_\ast\widetilde{\rho}=\rho.
$$
\paragraph{Step 3.}
We will show that
\begin{equation}\label{eq:rownosc}
M\circ (\varphi \times id_{\{0,1\}^\Z})=\Phi\circ \Psi.
\end{equation}
Then, for any $\nu\in\mathcal{P}^e(S,Y\cap \widetilde{X}_\eta)$,
$$
\nu=\Phi_\ast \Theta_\ast \nu=\Phi_\ast \Psi_\ast \widetilde{\Theta_\ast \nu}
=M_\ast (\varphi\times id_{\{0,1\}^\Z})_\ast \widetilde{\Theta_\ast \nu},
$$
with $\widetilde{\Theta_\ast \nu}\in\mathcal{P}( T\times S, G \times \{0,1\}^\Z)$. 
\paragraph{Step 4.}
To conclude it suffices to notice that
$$
\varphi\times id_{\{0,1\}^\Z} \colon  G \times \{0,1\}^\Z \to X_\eta \times \{0,1\}^\Z
$$
induces a map from $\mathcal{P}(T\times S, G \times \{0,1\}^\Z)$ to the simplex of probability $S\times S$-invariant measures on $X_\eta\times \{0,1\}^\Z$ whose projection onto the first coordinate is $\nu_\eta$.
\begin{Remark}
The above sketch can be summarized on the following commuting diagram:

$$
\begin{tikzpicture}[baseline=(current  bounding  box.center)]
%punkty 
\node (A1) at (0,0) {$Y\cap \widetilde{X}_\eta$};
\node (A2) at (3,0) {$Y\cap \widetilde{X}_\eta$};
\node (A3) at (6,0) {$G \times\{0,1\}^{\Z}$};
\node (A4) at (9,0) {$G \times\{0,1\}^{\Z}$};
\node (B1) at (0,-3) {$G \times\{0,1\}^{\Z}$};
\node (B2) at (3,-3) {$G \times\{0,1\}^{\Z}$};
\node (B3) at (6,-3) {$X_\eta \times \{0,1\}^\Z$};
\node (B4) at (9,-3) {$X_\eta \times \{0,1\}^\Z$};
\node (C2) at (3,-6) {$X_\sB$};
\node (C3) at (6,-6) {$X_\sB$};
%podpisy na strzałkach
\node (thetaL) at (-0.2,-1) {$\Theta$};
\node (thetaR) at (2.8,-1) {$\Theta$};
\node at (3.8,-0.8) {$\Psi$};
\node at (6.8,-0.8) {$\Psi$};
\node at (6.6,-1.7) {$\varphi\times Id$};
\node at (9.6,-1.7) {$\varphi\times Id$};
\node at (0.9,-4.2) {$\Phi$};
\node at (3.9,-4.2) {$\Phi$};
\node at (5.2, -4.2) {$M$};
\node at (8.2, -4.2) {$M$};
\draw[->] 
	(A1) edge node[auto] {$S$} (A2)
	(A3) edge node[auto] {$T\times S$} (A4)
	(B1) edge node[auto] {$\widetilde{T}$} (B2)
	(B3) edge node[auto] {$T\times S$} (B4)
	(C2) edge node[auto] {$S$} (C3)
%	(A3) edge node[auto] {$\Psi$} (B1)
%	(A4) edge node[auto] {$\Psi$} (B2)
	(A3) edge (B1)
	(A4) edge (B2)
%	(B1) edge node[below] {$\Phi$} (C2)
%	(B2) edge node[below] {$\Phi$} (C3);	
	(B1) edge (C2)
	(B2) edge (C3);
\draw[dashed,->] 
%	(A1) edge node[auto] {$\Theta$} (B1)
%	(A2) edge node[auto] {$\Theta$} (B2)
	(A1) edge (B1)
	(A2) edge (B2)
	(B3) edge (C2)
	(B4) edge (C3)
%	(B3) edge node[auto] {$M$} (C2)
%	(B4) edge node[auto] {$M$} (C3)
	(A3) edge (B3)
	(A4) edge (B4);
%	(A3) edge node[auto] {$\varphi\times Id$} (B3)
%	(A4) edge node[auto] {$\varphi\times Id$} (B4);
\end{tikzpicture}
$$
\end{Remark}

\begin{proof}[Proof of \cref{twY}]
Let
$$
Y_\infty:=\{y\in Y : |\text{supp }y \cap (-\infty,0)|=|\text{supp }y \cap (0,\infty)|=\infty\}.\footnote{Notice that the definition of $Y_\infty$ is different from the one in~\cite{MR3356811} -- we have changed the notation to simplify the proof.}
$$
Since $\nu\neq \delta_{(\dots,0,0,0,\dots)}$, we have $\nu(Y_\infty)=1$. For $x\in \{0,1\}^\Z$, $z\in Y_\infty$,
let $\widehat{x}_z$ be the sequence obtained by reading consecutive coordinates of $x$ which are in $\text{supp}\ z$, and such that
$$
\widehat{x}_z(0)=x(\min\{k\geq 0 : k\in \text{supp}\, z\}).
$$
\paragraph{Step 1.}
By \cref{1.2.6}, for $y\in Y_\infty\cap \widetilde{X}_\eta$, we have $\varphi(\theta(y))\in Y_\infty$.
Let $\Theta\colon Y_\infty\cap\widetilde{X}_\eta\to  G \times\{0,1\}^\Z$ be given by
$$
\Theta(y)=(\theta(y),\widehat{y}_{\varphi(\theta(y))}).
$$
One can show that
\begin{equation}\label{eq:wzorekes}
\widehat{Sx}_{Sz}=\begin{cases}
\widehat{x}_z,&\text{if } z(0)=0,\\
S\widehat{x}_z,&\text{if } z(0)=1.
\end{cases}
\end{equation}
Hence, in view of~\eqref{fikomut} and \cref{1.2.6},
it follows that $\Theta\circ S = \widetilde{T}\circ \Theta$ on~${Y_\infty}$.

Let $\Phi\colon\varphi^{-1}(Y_\infty)\times \{0,1\}^\Z \to X_\mathscr{B}$ be the unique element in $X_\mathscr{B}$ such that
$$
\Phi( g ,x)\leq \varphi( g )\text{ and }{\ \mkern65mu \widehat{\phantom{x}} \mkern-80mu(\Phi( g ,x))}_{\varphi( g )}=x.
$$
Since $\nu(Y\cap \widetilde{X}_\eta)=1$, by \cref{tautychar}, we have that $\nu_\eta(Y\cap \widetilde{X}_\eta)=1$, so, in particular, $\nu_\eta\neq \delta_{(\dots,0,0,0,\dots)}$. It follows that $\Phi$ is well-defined a.e.\ with respect to any $\widetilde{T}$-invariant measure. Moreover, using~\eqref{eq:wzorekes}, one can show that $S\circ \Phi=\Phi\circ \widetilde{T}$ on $\varphi^{-1}(Y_\infty)\times \{0,1\}^\Z$. It follows that also $\Phi\circ\Theta$ is well-defined a.e.\ with respect to any $\nu\in\mathcal{P}(S,Y\cap \widetilde{X}_\eta)$. Moreover, by the choice of $\Theta$ and $\Phi$, we obtain $\Phi\circ \Theta=id$ a.e.\ with respect to any $\nu\in\mathcal{P}(S,Y\cap \widetilde{X}_\eta)$.
\paragraph{Step 2.} Let $\Psi\colon \varphi^{-1}(Y_\infty)\times \{0,1\}^\Z\to \varphi^{-1}(Y_\infty)\times \{0,1\}^\Z$ be given by $\Psi( g ,x)=( g , \widehat{x}_{\varphi( g )})$. Using again~\eqref{eq:wzorekes}, one can show that diagram~\eqref{diagram3} commutes. Notice that $\emptyset\neq\Psi^{-1}( g ,y)\subset \{ g \}\times \{0,1\}^\Z$. Moreover, given $( g ,x)\in \Psi^{-1}( g ,y)$, all other points in $\Psi^{-1}( g ,y)$ are obtained by changing in an arbitrary way these coordinates in $x$ which are not in the support of $\varphi( g )$. In particular, each fiber $\Psi^{-1}( g ,y)$ is infinite.  For $k_1<\dots <k_s$ and $(i_1,\dots, i_s)\in \{0,1\}^s$, we define the following cylinder set:
\begin{equation}\label{cylinder}
C=C^{i_1,\dots,i_s}_{k_1,\dots,k_s}:=\{x\in \{0,1\}^\Z : x({k_j})=i_j, 1\leq j\leq s\}.
\end{equation}
For each such $C$ and for $A\in \mathcal{B}( G )$, we put
$$
\lambda_{( g ,y)}(A\times C):=\raz_A( g ) \cdot 2^{-m}, \text{ where }m=|\{1\leq j\leq s : \varphi( g )(k_j)=0\}|,
$$
if $ \Phi( g ,y)(k_j)=i_j\text{ whenever }\varphi( g )(k_j)=1$ (otherwise we set $\lambda_{( g ,y)}(A\times C):=0$). Conditions (a) and (b) required in Step 2. are proven in the same way as in~\cite{MR3356811}.
\paragraph{Step 3. and Step 4.}
Formula \eqref{eq:rownosc} follows directly by the choice of $\Phi$ and $\Psi$ and the proof is complete.
\end{proof}

%%%%%%%%%%%%%%%%%%%%%%%%%%%%%%%%%%%%%%%%%%%%%%%%%%%%%%%%%%%%%%%%%%%%%%%%%%%%%%%%%%
\subsection{Invariant measures on $\widetilde{X}_\eta$ (proof of \cref{OOG})}\label{se:9.3}
In this section we will prove the following:
\begin{Th}\label{miaryall}
For any $\nu\in\mathcal{P}^e(S,\widetilde{X}_\eta)$ there exist $b_k'\divides b_k$, $k\geq 1$,
 and $\widetilde{\rho}\in\mathcal{P}^e(S\times S,{X}_{\eta'}\times\{0,1\}^\Z)$ such that $\widetilde{\rho}|_{{X}_{\eta'}}=\nu_{\eta'}$ and $M_\ast(\widetilde{\rho})=\nu$, where $\eta'=\raz_{\cf_{\sB'}}$ ($\sB'=\{b_k' : k\geq 1\}$) and $M\colon X_{\eta'}\times \{0,1\}^\Z\to \widetilde{X}_{\eta'}$ stands for the coordinatewise multiplication.
\end{Th}

For the proof we will need several tools. Notice first that if $\nu=\delta_{(\dots,0,0,0,\dots)}$ then the above assertion holds true since $M_\ast(\delta_{(\dots,0,0,0,\dots)}\otimes \kappa)=\delta_{(\dots,0,0,0,\dots)}$ for any $\kappa\in \mathcal{P}(S,\{0,1\}^\Z)$, and $\delta_{(\dots,0,0,0,\dots)}=\nu_{\eta'}$ for $\eta'$ associated to $\mathscr{B}'=\{1\}$. Thus, we only need to cover the case $\nu\neq\delta_{(\dots,0,0,0,\dots)}$.

Recall that
$$
\widetilde{X}_\eta=\bigcup_{k\geq 1}\bigcup_{0\leq s_k\leq b_k}Y_{s_1,s_2,\dots}\cap \widetilde{X}_\eta
$$
is a partition of $\widetilde{X}_\eta$ into Borel, $S$-invariant sets.  Proceeding in a similar way as in~\cite{MR3356811}, we will now further refine this partition.

Fix $\underline{s}=(s_k)_{k\geq 1}$ with $1\leq s_k\leq b_k-1$, $\underline{a}=(a_1^{k},\dots,a_{s_k}^{k})_{k\geq 1}$ with $a_i^k\in \Z/b_k\Z$ for $1\leq i\leq s_k$ and $|\{a_1^k,\dots,a_{s_k}^k\}|=s_k$. Let
$$
Y_{k,s_k;a_1,\ldots,a_{s_k}}:=\{x\in \{0,1\}^\Z:
\text{supp }x\bmod b_k=\Z/b_k\Z\setminus\{a_1,\ldots,a_{s_k}\}\}.
$$
For each $k\geq 1$, any two sets of such form are either disjoint or they coincide. Since $\text{supp }Sx=\text{supp }x-1$, we have
\begin{equation}\label{her11}
SY_{k,s_k;a_1^k,\ldots,a_{s_k}^k}=Y_{k,s_k;a_1^k-1,\ldots,a_{s_k}^k-1}.
\end{equation}
Let
\begin{equation}\label{primy}
b_k':=\min\{j\geq 1 : \{a_1^k,\dots,a_{s_k}^k\}=\{a_1^k-j,\dots,a_{s_k}^k-j\}\}
\end{equation}
and note that $b_k'\geq 2$. Clearly, $S^{b'_k}Y_{k,s_k;a_1^k,\ldots,a_{s_k}^k}=Y_{k,s_k;a_1^k,
\ldots,a_{s_k}^k}$ and the sets
$$
Y_{k,s_k;a_1^k,\ldots,a_{s_k}^k},SY_{k,s_k;a_1^k,\ldots,a_{s_k}^k},\ldots,
S^{b'_k-1}Y_{k,s_k;a_1^k,\ldots,a_{s_k}^k}
$$
are pairwise disjoint. Finally, we define
$$
Y_{\underline{s},\underline{a}}:=\bigcap_{k\geq 1}\bigcup_{j=0}^{b_k'-1}S^j Y_{k,s_k;a_1^k,\dots,a_{s_k}^k}
$$
(notice that if $s_k=1$ for all $k\geq 1$, we have $Y_{\underline{s},\underline{a}}=Y$ for any choice of $\underline{a}$).

Fix $\underline{s},\underline{a}$ and suppose that $\mathcal{P}(S,{Y}_{\underline{s},\underline{a}}\cap \widetilde{X}_\eta)\neq \emptyset$. Let
$$
G_\sa:=\overline{\{\underline{n}_{\mathscr{B}'} : n\in\Z\}}\subset  G _{\mathscr{B}'}=\prod_{k\geq 1}\Z/b_k'\Z,
$$
where $b_k'$, $k\geq 1$, are as in \eqref{primy}, cf.\ \eqref{GIE}. Define $\varphi_{\underline{s},\underline{a}}\colon  G _{\underline{s},\underline{a}}\to \{0,1\}^\Z$ by
$$
\varphi_{\underline{s},\underline{a}}( g )(n)=\begin{cases}
1,&\text{if }  g_k-a_i^k+n\neq 0 \bmod b_k \text{ for all }k\geq 1, 1\leq i\leq s_k,\\
0,&\text{otherwise}
\end{cases}
$$
(cf.\ \eqref{dze8equi}). We also define $\theta_{\underline{s},\underline{a}}\colon Y_{\underline{s},\underline{a}}\cap \widetilde{X}_\eta \to  G_{\sB'}$  in the following way:
$$
\theta_{\underline{s},\underline{a}}(y)= g  \iff - g_k+a_i^k \not\in\text{supp}(y) \bmod b_k \text{ for all }1\leq i\leq s_k,
$$
cf.\ \eqref{THETA}.
Notice that given $y\in Y_\sa$ and $k_0\geq 1$, there exists $N\geq 1$ such that
\begin{equation}\label{wkskW}
(\text{supp }y)\cap [-N,N]\bmod b_k=\Z/b_k\Z\setminus\{-g_k+a_i^k : 1\leq i \leq s_k\} \text{ for }1\leq k\leq k_0
\end{equation}
\begin{Remark}
Notice that
\begin{equation}\label{jestwGW}
\theta_\sa(Y_\sa\cap \widetilde{X}_\eta)\subset G_\sa.
\end{equation}
Indeed, take $y\in Y_\sa\cap \widetilde{X}_\eta$. Given $k_0\geq 1$, let $N\geq 1$ be such that~\eqref{wkskW} holds and let $M\in\Z$ be such that $y[-N,N] \leq \eta[-N+M,N+M]$. It follows that $\theta(y)=(g_1,g_2,\dots)$, where $g_k\equiv -M\bmod b_k$ for $1\leq k\leq k_0$. This yields~\eqref{jestwGW}.
\end{Remark}
\begin{Remark}\label{cgtheta}
Note also that $\theta_\sa$ is continuous. Indeed, given $y\in Y_\sa$ and $k_0\geq 1$, let $N$ be such that~\eqref{wkskW} holds. Then, if $y'\in Y_\sa$ is sufficiently close to $y$ then~\eqref{wkskW} holds for $y'$ as well. Therefore, if $y_n\to y$ in $Y_\sa$ then $\theta_\sa(y_n)\to \theta_\sa(y)$.
\end{Remark}
Moreover, denote by $T_{\underline{s},\underline{a}}\colon G _{\underline{s},\underline{a}}\to G _{\underline{s},\underline{a}}$ the map given by
$$
T_\sa g = g +\underline1_{\mathscr{B}'}=(g_1+1,g_2+1,\dots),
$$
where $ g =(g_1,g_2,\dots)$.
\begin{Remark}[cf.\ \cref{1.2.6}]\label{1.2.6a}
We have:
\begin{itemize}
	\item%\label{F1}
	$T_\sa\circ \theta_\sa=\theta_\sa \circ S$,
	\item%\label{F2}
	for each $y\in Y_\sa\cap \widetilde{X}_\eta$, $y\leq \varphi_\sa( \theta_\sa (y))$,
	\item for any $\nu\in\mathcal{P}(S,Y_\sa\cap \widetilde{X}_\eta)$, $(\theta_\sa)_\ast(\nu)=\PP_\sa$.
\end{itemize}
\end{Remark}
\begin{Lemma}\label{notry}
Suppose that $\mathcal{P}(S,Y_\sa\cap \widetilde{X}_\eta)\neq\emptyset$. Then $(\varphi_\sa)_\ast(\PP_\sa)(Y_\sa)=1$. In particular, $(\varphi_\sa)_\ast(\PP_\sa)\neq \delta_{(\dots,0,0,0,\dots)}.$
\end{Lemma}
\begin{proof}
Take $\nu\in\mathcal{P}(S,Y_\sa\cap \widetilde{X}_\eta)$. It follows by \cref{1.2.6a} that
$$
(\varphi_\sa)_\ast(\PP_\sa)(Y_\sa)=(\varphi_\sa)_\ast(\theta_\sa)_\ast(\nu)(Y_\sa)\geq \nu(Y_{\sa})=1.
$$
Since $(\dots,0,0,0,\dots)\not\in Y_\sa$, we conclude.
\end{proof}
For $n\in\N$, let $M^{(n)}\colon (\{0,1\}^\Z)^{\times n}\to \{0,1\}^\Z$ be given by
$$
M^{(n)}((x^{(1)}_i)_{i\in\Z},\dots,(x^{(n)}_i)_{i\in\Z})=(x^{(1)}_i\cdot \ldots \cdot x^{(n)}_i)_{i\in\Z}.
$$
Moreover, we define $M^{(\infty)} \colon (\{0,1\}^\Z)^{\N}\to \{0,1\}^\Z$ as
$$
M^{(\infty)}((x^{(1)}_i)_{i\in\Z},(x^{(2)}_i)_{i\in\Z},\dots)=(x^{(1)}_i\cdot x^{(2)}_i\cdot \ldots  )_{i\in\Z}.
$$
\begin{Lemma}[cf.\ Lemma 2.2.22 in~\cite{MR3356811}]\label{may21a}
We have
$
(\varphi_\sa)_\ast (\PP_\sa)=M^{(\infty)}_\ast (\rho),
$
where $\rho$ is a joining of a countable number of copies of $(S,\{0,1\}^\Z,\nu_{\eta'})$.
\end{Lemma}
\begin{proof}
The proof is the same as in~\cite{MR3356811}.
\end{proof}

\begin{Lemma}[Lemma 2.2.23 in~\cite{MR3356811}]\label{may21b}
Let $\nu_1,\dots, \nu_n,\nu_{n+1}\in\mathcal{P}(S,\{0,1\}^\Z)$. Then for any joinings
\begin{itemize}
\item
$\rho_{1,n}\in J((S,\{0,1\}^\Z,\nu_1),\dots, (S,\{0,1\}^\Z,\nu_n))$,
\item
$\rho_{(1,n),n+1}\in J((S,\{0,1\}^\Z,M^{(n)}_\ast(\rho_{1,n})),(S,\{0,1\}^\Z,\nu_{n+1}))$
\end{itemize}
there exist:
\begin{itemize}
\item
${\rho}_{2,n+1}\in J((S,\{0,1\}^\Z,\nu_2),\dots,(S,\{0,1\}^\Z,\nu_n),(S,\{0,1\}^\Z,n+1))$,
\item
${\rho}_{1,(2,n+1)}\in J((S,\{0,1\}^\Z,\nu_1),(S,\{0,1\}^\Z,M^{(n)}_\ast({\rho}_{2,n+1})))$
\end{itemize}
such that $M^{(2)}_\ast(\rho_{(1,n),n+1})=M^{(2)}_\ast(\rho_{1,(2,n+1)})$.\footnote{We could write this property as $M^{(2)}_\ast(M^{(n)}_\ast(\nu_1\vee\dots\vee\nu_n)\vee \nu_{n+1})=M^{(2)}_\ast(\nu_1\vee M^{(n)}_\ast(\nu_2\vee \dots\vee\nu_n\vee\nu_{n+1}))$. However, until we say which joining we mean by each symbol $\vee$, this expression has no concrete meaning.}
\end{Lemma}
\begin{Remark}\label{may21c}
The above lemma remains true when we consider infinite joinings, i.e.\ instead of $\nu_1,\dots,\nu_n$ we have $\nu_1,\nu_2,\dots$, and instead of $M^{(n)}$ we consider $M^{(\infty)}$.
\end{Remark}

\begin{proof}[Proof of \cref{miaryall}]
Fix $\delta_{(\dots,0,0,0,\dots)}\neq\nu\in\mathcal{P}^e(S,\widetilde{X}_\eta)$ and let $\sa$ be such that $\nu(Y_\sa\cap\widetilde{X}_\eta)=1$. In view of Lemma~\ref{may21a},~Lemma~\ref{may21b} and \cref{may21c}, it suffices to show that there exists $\widetilde{\rho}\in\mathcal{P}(S\times S,\{0,1\}^\Z\times\{0,1\}^\Z)$ such that
the projection of $\widetilde{\rho}$ onto the first coordinate equals $(\varphi_\sa)_\ast(\PP_\sa)$ and $M_\ast(\widetilde{\rho})=\nu$.

By~\cref{notry}, we have $(\varphi_\sa)_\ast(\PP_\sa)\neq \delta_{(\dots,0,0,0,\dots)}$. The remaining part of the proof goes exactly along the same lines as the proof of \cref{twY}, with the following modification: we need to replace some objects related to $Y$ by their counterparts related to $Y_\sa$. Namely, instead of $ G ,\ \Theta,\ Y_\infty,\ \widetilde{T},\ \Phi$ and $\Psi$,
we use
$$
 G _\sa,\ \Theta_\sa,\ (Y_\sa)_\infty,\  \widetilde{T}_\sa,\ \Phi_\sa\text{ and } \Psi_\sa,
$$
where
\begin{itemize}
\item
$\Theta_\sa\colon Y_\sa\cap\widetilde{X}_\eta\to  G _\sa\times\{0,1\}^\Z$ is given by $\Theta_\sa(y):=(\theta_\sa(y),\widehat{y}_{\varphi_\sa(\theta_\sa y)})$,
\item
$(Y_\sa)_\infty:=\{y\in Y_\sa : |\text{supp }y \cap (-\infty,0)|=|\text{supp }y \cap (0,\infty)|=\infty\},$
\item
$\widetilde{T}_\sa\colon  G _\sa\times \{0,1\}^\Z\to  G _\sa\times \{0,1\}^\Z$ given by
$$
\widetilde{T}_\sa( g ,x)=\begin{cases}
(T_\sa g , x),&\text{if }\varphi_\sa( g )(0)=0,\\
(T_\sa g , Sx),&\text{if } \varphi_\sa( g )(0)=1,
\end{cases}
$$
\item
$\Phi_\sa( g , x)$ is the unique element in $X_\mathscr{B}$ such that
\begin{enumerate}[(i)]
\item\label{eq:C1}
$\Phi_\sa( g ,x)\leq \varphi_\sa( g )$,
\item
$\mkern85mu \widehat{\phantom{x}} \mkern-100mu {({\Phi_\sa( g ,x)})}_{\varphi_\sa( g )}=x$, i.e.\ the consecutive coordinates of $x$ can be found in $\Phi_\sa( g ,x)$ along $\varphi_\sa( g )$,
\end{enumerate}
\item
$\Psi_\sa( g ,x)=( g , \widehat{x}_{\varphi_\sa( g )})$.
\end{itemize}
\end{proof}

Repeating the proof of \cref{may21a}, we obtain the following:
\begin{Lemma}\label{lema10}
Fix $b_k' \divides b_k$ for $k\geq 1$. Then there exists $\rho\in \mathcal{P}(S\times S,X_\eta\times \{0,1\}^\Z)$ such that $\rho|_{X_\eta}=\nu_\eta$ and $M_\ast(\rho)=\nu_{\eta'}$.
\end{Lemma}

\cref{OOG} is a consequence of \cref{twY}, \cref{miaryall}, Lemma~\ref{lema10},~Lemma~\ref{may21b} and \cref{may21c}.

%%%%%%%%%%%%%%%%%%%%%%%%%%%%%%%%%%%%%%%%%%%%%%%%%%%%%%%%%%%%%%%%%%%%%%%%%%%%%%%%%%
\subsection{Rational discrete spectrum (proof of \cref{IZOIZO})}\label{sewymierne}

\begin{Remark}
Let $\un{s},\un{a}$ be such that $\mathcal{P}(S,Y_{\un{s},\un{a}})\neq\emptyset$ and fix $\nu\in \mathcal{P}(S,Y_{\un{s},\un{a}})$. Let $b_k'\divides b_k$, $k\geq1$, be as in the proof of \cref{miaryall}. Recall (from the proof of \cref{miaryall}) that there is an equivariant map $\Theta_\sa\colon Y_\sa\to  G _\sa\times\{0,1\}^\Z$. It follows that $(T_{\un{s},\un{a}},G_\sa,\PP_\sa)$ is a factor of $(S,Y_{\un{s},\un{a}},\nu)$. In particular, the rational discrete spectrum of $(S,Y_{\un{s},\un{a}},\nu)$ includes all $b_k'$-roots of unity.
\end{Remark}

\begin{Th}\label{notrynotry}
Suppose that $\mathcal{P}(S,Y_{\un{s},\un{a}}\cap \widetilde{X}_\eta)\neq\emptyset$. Then $\varphi_\sa$ yields an isomorphism of $(T_\sa,G_\sa,\PP_\sa)$ and $(S,Y_\sa\cap\widetilde{X}_\sa,(\varphi_\sa)_\ast(\PP_\sa))$.
\end{Th}
\begin{proof}
Since, by \cref{notry}, we have $(\varphi_\sa)_\ast(\PP_\sa)(Y_\sa)=1$, we obtain the following equivariant maps:
$$
(T_\sa,G_\sa,\PP_\sa)\xrightarrow{\varphi_\sa} (S,Y_\sa\cap \widetilde{X}_\eta,(\varphi_\sa)_\ast(\PP_\sa))
\xrightarrow{\theta_\sa}  (T_\sa,G_\sa,\PP_\sa).
$$
It follows by the coalescence of $(T_\sa,G_\sa,\PP_\sa)$ that $\varphi_\sa$ yields an isomorphism of $(T_\sa,G_\sa,\PP_\sa)$ and $(S,Y_\sa,(\varphi_\sa)_\ast(\PP_\sa))$.
\end{proof}
As an immediate consequence of the above and of \cref{gdziemax}, we obtain \cref{IZOIZO}.

\section{Tautness revisited}

%\subsection{Taut subsystem supporting all measures (proof of \cref{TTC})}\label{seatraktor}
\subsection{Tautness and combinatorics revisited (proof of \cref{OOJ})}\label{Jfirst1}

We will prove an extension of \cref{zawiewnio} and \cref{zawiewnio1}.
\begin{Cor}\label{PP:1}
Let $\sB,\sB'\subset \N$ and suppose that $\mathscr{B}$ is taut. Conditions \eqref{J:A} - \eqref{J:F} from \cref{zawiewnio} are equivalent to each of the following:
\begin{enumerate}[(a)]
\setcounter{enumi}{6}
\item\label{J:g}
$\nu_\eta\in\mathcal{P}(S,\widetilde{X}_{\eta'})$,
%$\nu_{\eta'}\in\mathcal{P}(S,\widetilde{X}_{\eta})$,
\item\label{J:h}
%$\mathcal{P}(S,\widetilde{X}_{\eta'})\subset\mathcal{P}(S,\widetilde{X}_{\eta})$.
$\mathcal{P}(S,\widetilde{X}_\eta)\subset\mathcal{P}(S,\widetilde{X}_{\eta'})$.
\end{enumerate}
\end{Cor}
\begin{proof}
Notice first that~\eqref{J:E} from \cref{zawiewnio} implies~\eqref{J:g}. Suppose now that~\eqref{J:g} holds. In view of \cref{OOG} and \cref{may21b}, this yields~\eqref{J:h}. Suppose that~\eqref{J:h} holds. By the variational principle, we have
\begin{equation}\label{P1}
h_{top}(S,\widetilde{X}_{\eta})=h_{top}(S,\widetilde{X}_\eta\cap\widetilde{X}_{\eta'}).
\end{equation}
Moreover, since $\widetilde{X}_\eta\cap\widetilde{X}_{\eta'}\subset X_{\mathscr{B}}\cap X_{\mathscr{B}'}=X_{\sB\cup \sB'}\subset X_{\mathscr{B}}$,
\begin{equation}\label{P2}
h_{top}(S,\widetilde{X}_\eta\cap\widetilde{X}_{\eta'})\leq h_{top}(S,X_{\mathscr{B}\cup\mathscr{B}'}) \leq h_{top}(S,X_{\mathscr{B}}).
\end{equation}
By \cref{OOI}, 
\begin{equation}\label{P3}
h_{top}(S,\widetilde{X}_{\eta})=h_{top}(S,X_{\mathscr{B}}).
\end{equation}
%Since, $X_\sB\subset X_{\sB\cup\{b'\}}\subset  X_{\sB\cup\sB'}$ for any $b'\in\sB'$, it follows that 
%$$
%h_{top}(S,X_{\sB})
%$$
Putting together \eqref{P1}, \eqref{P2} and \eqref{P3}, we obtain 
\begin{equation}\label{P4}
h_{top}(S,X_{\sB})=h_{top}(S,X_{\sB\cup \sB'}).
\end{equation}
Moroever, since $X_{\sB\cup\sB'}\subset X_{\sB\cup\{b'\}}\subset X_{\sB}$ for any $b'\in\sB'$,~\eqref{P4} yields
$$
h_{top}(S,X_{\mathscr{B}})=h_{top}(S,X_{\mathscr{B}\cup \{b'\}}) \text{ for any }b'\in\mathscr{B}'.
$$
It follows by \cref{OOI} that
$$
\bdelta(\mathcal{M}_{\mathscr{B}})=\bdelta(\mathcal{M}_{\mathscr{B}\cup \{b'\}}).
$$
In view of \cref{behnowy}, either $b'\in\mathcal{M}_{\mathscr{B}}$ or $\mathscr{B}$ is not taut. The latter is impossible, hence $b \divides b'$ for some $b\in\mathscr{B}$ and we conclude that \eqref{J:B} from \cref{zawiewnio} holds.
\end{proof}
\begin{Cor}
Suppose that $\mathscr{B},\mathscr{B}'\subset \N$ are taut. Conditions \eqref{J:A1} - \eqref{J:G1} from \cref{zawiewnio1} are equivalent to each of the following:
\begin{enumerate}[(a)]
\setcounter{enumi}{7}
\item\label{J:H1}
$\nu_\eta=\nu_{\eta'}$,
\item\label{J:I1}
$\nu_\eta\in\mathcal{P}(S,\widetilde{X}_{\eta'})$ and $\nu_{\eta'}\in\mathcal{P}(S,\widetilde{X}_{\eta})$,
\item\label{J:J1}
$\mathcal{P}(S,\widetilde{X}_\eta)= \mathcal{P}(S,\widetilde{X}_{\eta'})$.
\end{enumerate}
\end{Cor}
\begin{proof}
Clearly, \eqref{J:C1} from \cref{zawiewnio1} together with \cref{OOE} implies \eqref{J:H1}. Moreover, \eqref{J:H1} implies \eqref{J:I1} and, by \cref{PP:1}, \eqref{J:I1} implies \eqref{J:J1}. 
Suppose now that \eqref{J:J1} holds. Applying again \cref{PP:1}, we obtain that \eqref{J:B} from \cref{zawiewnio} holds. Moreover, \eqref{J:B} from \cref{zawiewnio} still holds when we exchange the roles of $\mathscr{B}$ and $\mathscr{B}'$. Therefore, using \eqref{J:A} from \cref{zawiewnio}, we conclude that $X_{\mathscr{B}'}=X_{\mathscr{B}}$, i.e.\ \eqref{J:A} from \cref{zawiewnio1} holds. This completes the proof. 
\end{proof}

%\subsection{tauty miary}
\subsection{Tautness and invariant measures (proof of \cref{TTC})}\label{seatraktor}
\cref{TTC} is an immediate consequence of \cref{tautywszystko}, \cref{OOG} and \cref{OOJ}. We will now prove \cref{TTC1}. For this, we will need the following standard lemma:
\begin{Lemma}\label{lemmaatractor}
Let $(T,X)$ be a topological dynamical system and let $X'\subset X$ be compact and $T$-invariant. Then the following are equivalent:
\begin{enumerate}[(a)]
\item\label{ecli1}
$\mathcal{P}(T,X)=\mathcal{P}(T,X')$,
\item\label{ecli2}
for each $x\in X$, we have $\lim_{n\to \infty, n\not\in E_x}d(T^nx,X')=0,\ \text{where }d(E_x)=0$.
\end{enumerate}
\end{Lemma}
\begin{proof}
We will show first \eqref{ecli1} $\Rightarrow$ \eqref{ecli2}. Assume that we have \eqref{ecli1}. Suppose that \eqref{ecli2} does not hold for some $x\in X$, i.e.\ there exist $\delta>0$ and $E_x\subset \Z$ with $\ov{d}(E_x)>0$ such that
\begin{equation}\label{wiosna3}
d(T^nx,X')\geq \delta \text{ for }n\in E_x.
\end{equation}
Let $f\in C(X)$ be such that $0\leq f\leq 1$, $f(x)=1$ if $d(x,X')\geq \delta$ and $f(x)=0$ if $x\in X'$. Let $(N_k)_{k\geq 1}$ and $\nu\in \mathcal{P}(T,X)$ be such that
\begin{equation}\label{kabel3}
d_{(N_k)}(E_x):=\lim_{k\to\infty}\frac{1}{N_k}|E_x\cap [0,N_k]|>0
\end{equation}
and
\begin{equation}\label{kabel4}
\frac{1}{N_k}\sum_{n\leq N_k} \delta_{T^nx} \to \nu.
\end{equation}
Then, using \eqref{kabel4} and~\eqref{ecli1}, we obtain
\begin{equation}\label{previous}
\frac{1}{N_k}\sum_{n\leq N_k}f(T^nx)\to \int_X f \ d\nu=\int_{X'} f \ d\nu=0.
\end{equation}
On the other hand, by \eqref{wiosna3}, \eqref{kabel3} and by the definitions of $E_x$ and $f$, we have
$$
\lim_{k\to\infty}\frac{1}{N_k}\sum_{n\leq N_k}f(T^nx)\geq d_{(N_k)}(E_x)>0,
$$
which contradicts~\eqref{previous}.

We will now show that \eqref{ecli2} implies \eqref{ecli1}. Suppose that for some $\nu\in\mathcal{P}^e(T,X)$, we have $\nu(X')=0$. Let $X'\subset U\subset X$ be an open set, such that $\nu(U)<\vep$. Let $f\in C(X)$ be such that $0\leq f\leq 1$, $f(x)=1$ for $x\in X'$ and $f(x)=0$ for $x\in X\setminus U$. By the ergodicity of $\nu$, there exists $x\in X$ such that
\begin{equation}\label{wiosna1}
\frac{1}{N}\sum_{n\leq N}\delta_{T^nx}\to \nu.
\end{equation}
Then, by the choice of $U$ and $f$, we have
\begin{equation}\label{wiosna2}
\int f\ d\nu \leq \vep.
\end{equation}
On the other hand, using \eqref{ecli2} and \eqref{wiosna1}, we obtain
$$
\int f\ d\nu=\lim_{N\to \infty}\frac{1}{N}\sum_{n\leq N}f(T^nx)=1,
$$
which yields a contradiction with~\eqref{wiosna2} and completes the proof.
\end{proof}
\begin{Def}\label{uwaatra}
When (b) of \cref{lemmaatractor} holds, we say that $X'$ is a \emph{quasi-attractor} in $(T,X)$.
\end{Def}
\cref{TTC1} follows immediately by \cref{TTC} and by \cref{lemmaatractor}. Moreover, \cref{TTC1} can be rephrased as follows:
\begin{Cor}\label{maxtautowy}
For any $\mathscr{B}\subset \N$, the subshift $(S,\widetilde{X}_\eta)$ has a quasi-attractor of the form $\widetilde{X}_{\eta'}$ for some taut set $\mathscr{B}'$ such that $\cf_{\sB'}\subset \cf_\sB$. Moreover, such $\mathscr{B}'$ is unique.
\end{Cor}

%%%%%%%%%%%%%%%%%%%%%%%%%%%%%%%%%%%%%%%%%%%%%%%%%%%%%%%%%%%%%%%%%%%%%%%%%%%%%%%%%%%%%%%%%%%%%%%%%%
\section{Intrinsic ergodicity revisited}
\subsection{Taut case revisited}\label{serevisited}
Now we present a second proof of \cref{jakzbenjim}.
\begin{proof}[Proof of \cref{jakzbenjim}]
We will use the objects introduced in course of the proof of \cref{twY}. There exists $C_0\subset C$ (recall that $C$ was defined in~\eqref{tujestC}) such that every point from $C_0$ returns to $C$ infinitely often under $T$ and $\PP(C_0)=\PP(C)$. It follows that every point from $C_0\times \{0,1\}^\Z$ returns to $C\times \{0,1\}^\Z$ infinitely often under $\widetilde{T}$ and
$\nu(C_0\times\{0,1\}^\Z)=\nu(C\times \{0,1\}^\Z)$ for every $\nu \in \mathcal{P}(\widetilde{T},G\times \{0,1\}^\Z)$. Thus, the induced transformation $\widetilde{T}_{{C}\times \{0,1\}^\Z}$ is well-defined. Recall that 
$$
\widetilde{T}(g,x)=\begin{cases}
(Tg,x),& \text{ if } g\not\in C,\\
(Tg,Sx),& \text{ if } g\in C.
\end{cases}
$$
It follows that $\widetilde{T}_{C\times \{0,1\}^\Z}=T_{C}\times S$ a.e.\ for any $\widetilde{T}$-invariant measure (cf.\ the definition of $\widetilde{T}$ and $C$).

We will show now that $\widetilde{T}$ has a unique measure of maximal (measure-theoretic) entropy. In view of Abramov's formula, for this, it suffices to show that $\widetilde{T}_{C\times \{0,1\}}=T_C\times S$ has a unique measure of maximal entropy. For any $T_C\times S$-invariant measure $\kappa$, by the Pinsker formula, we have
\begin{align}
\begin{split}\label{huhu}
 &h(S,\{0,1\}^\Z,\kappa |_{\{0,1\}^\Z})\leq h(T_C\times S, C\times \{0,1\}^\Z,\kappa )\\
&\quad\quad\leq h(T_C,C,\kappa |_C) +h(S,\{0,1\}^\Z,\kappa |_{\{0,1\}^\Z})=h(S,\{0,1\}^\Z,\kappa |_{\{0,1\}^\Z}).
\end{split}
\end{align}
Since $\kappa|_{\{0,1\}^\Z}$ can be arbitrary, it follows that the maximal entropy for $T_C\times S$ and for $S$ is the same. Moreover, the maximal entropy for $T_C\times S$ is achieved by $\kappa$ if and only if the maximal entropy for $S$ is achieved by $\kappa|_{\{0,1\}^\Z}$. In other words, this happens if and only if $\kappa|_{\{0,1\}^\Z}$ is the Bernoulli measure $B(1/2,1/2)$, i.e.\ when $\kappa$ is a joining of the unique invariant measure for $T_C$ and $B(1/2,1/2)$. Since the unique invariant measure for $T_C$ is of zero entropy, it follows by disjointness~\cite{MR0213508} that $\kappa$ is the product measure. In particular, $\kappa$ is unique.

It follows from~\eqref{isidentity} that $\Theta$ is 1-1. Hence,  $\Theta_\ast \colon\mathcal{P}(S,Y\cap\widetilde{X}_\eta)\to\mathcal{P}(\widetilde{T}, G \times\{0,1\}^\Z)$ is also 1-1 and for any $\nu\in\mathcal{P}(S,Y\cap \widetilde{X}_\eta)$, we have $h(S,Y\cap\widetilde{X}_\eta,\nu)=h(\widetilde{T}, G \times\{0,1\}^\Z,{\Theta_\ast\nu})$. The result follows now from \cref{gdziemax}.
\end{proof}

\begin{Remark}\label{forma}
Suppose that $\mathscr{B}\subset \N$ is taut. Notice that we have
$
\Psi_\ast(\PP \otimes B(1/2,1/2))=\PP\otimes B(1/2,1/2).
$
Moreover,
$$
(\PP\otimes B(1/2,1/2))_{C\times \{0,1\}^\Z} = \PP_C \otimes B(1/2,1/2).
$$
Since $h(T_C\times S,C\times \{0,1\}^\Z,{\PP_C \otimes B(1/2,1/2)})=\log 2$, it follows by the above proof of \cref{jakzbenjim} that
\begin{multline*}
\Phi_\ast\Psi_\ast(\PP\otimes B(1/2,1/2))=M_\ast (\varphi\times id)_\ast (\PP\otimes B(1/2,1/2))\\
=M_\ast (\nu_\eta\otimes B(1/2,1/2))
\end{multline*}
is the unique measure of maximal entropy for $(S,\widetilde{X}_\eta)$. 
\end{Remark}

\subsection{General case (proof of \cref{OOH})}\label{revi2}
\cref{OOH} is an immediate consequence of \cref{TTC} and \cref{jakzbenjim}.

%%%%%%%%%%%%%%%%%%%%%%%%%%%%%%%%%%%%%%%%%%%%%%%%%%%%%%%%%%%%%%%%%%%%%%%%%%%%%%%%%%%%%%%%%%%%%%%%%%

\section{Remarks on number theory}\label{s12}
\subsection{Consecutive gaps between $\mathscr{B}$-free numbers (proof of \cref{OOL})}\label{conse}
Fix $\sB\subset \N$ and denote by $(n_j)_{j\geq 1}$ the sequence of consecutive natural $\mathscr{B}$-free numbers. In~\cite{MR2414205}, the following was shown in case when $\mathscr{B}\subset \N$ satisfies~\eqref{settingerdosa}:
\begin{equation}\label{quant}
  \parbox{0.8\linewidth}{Let $\delta,\sigma>0$ be such that $20\sigma>9+3606\delta$. Then, for $N$ large enough there exists $j=j(N)\geq 1$ such that $n_j\in [N,N+N^\sigma]$ and $\min(n_{j+1}-n_{j},n_{j}-n_{j-1})>\Phi(N)$, where $\Phi(N)$ is the largest positive integer such that $\prod_{j=1}^{3\Phi(N)}b_j\leq N^\delta$.}
\end{equation}
In particular,
\begin{equation}\label{dwojka2}
\limsup_{j\to \infty} \inf(n_{j+2}-n_{j+1},n_{j+1}-n_{j})=\infty.
\end{equation}

\begin{proof}[Proof of \cref{OOL}]
It follows by \cref{TTD} that $X_\eta=\widetilde{X}_\eta$. Moroever, by \cref{OOF}, $X_\eta$ is the topological support of $\nu_\eta$. Since, by \cref{OOE}, $\eta$ is quasi-generic for $\nu_\eta$, the result follows.
\end{proof}
Even though, contrary to~\eqref{quant}, the result included in \cref{OOL} is not quantitative, it seems new and it strengthens~\eqref{dwojka2}.

\subsection{Abundant numbers}\label{seabu}
\begin{Def}
For $n\in \N$, consider the aliquot sum $s(n):=\sum_{d\divides n, d<n}d$. We say that $n\in \N$ is:
\begin{enumerate}[(i)]
\item
\emph{abundant} if $s(n)>n$,
\item
\emph{perfect} if $s(n)=n$,
\item
\emph{deficient} if $s(n)<n$.
\end{enumerate}
We will denote the set of abundant, perfect and deficient numbers by $\mathbf{A}$, $\mathbf{P}$ and $\mathbf{D}$, respectively.
\end{Def}
Notice that $\mathbf{A}$ is closed under taking multiples. It follows that
$$
\mathbf{A}=\N\cap\mathcal{M}_{\mathscr{B}_\mathbf{A}} \text{ and }\mathbf{P}\cup \mathbf{D} =\N\cap\mathcal{F}_{\mathscr{B}_\mathbf{A}}
$$
for some primitive ${\mathscr{B}_\mathbf{A}}\subset \N$. 
\begin{Lemma}\label{lm121}
$\sB_\mathbf{A}$ is thin. In particular, $\sB_\mathbf{A}$ has light tails and is Besicovitch.
\end{Lemma}
\begin{proof}
Erd\"os~\cite{MR1574879} showed 	that ${\mathscr{B}_\mathbf{A}} \cap [0,n]={\rm{o}}(n/\log^2n)$. Let $n_j$ be the $j$-th ${\mathscr{B}_\mathbf{A}}$-free natural number. Therefore, for $n$ sufficiently large, 
$
n\leq {j_n}({\log ^2 j_n})^{-1}.
$
It follows that, for large $n$, we have $n\log^2 n\leq n\log^2 j_n\leq j_n$, whence
\begin{equation}\label{abuthin}
\sum_{b\in{\mathscr{B}_\mathbf{A}}}\nicefrac{1}{b}=\sum_{n\geq 1}\nicefrac{1}{j_n}\leq \sum_{n\geq 1}\nicefrac{1}{n\log^2n}<\infty, 
\end{equation}
i.e.\ ${\mathscr{B}_\mathbf{A}}$ is thin. To complete the proof, it suffices to use the fact that thin sets are Besicovitch.
\end{proof}

\begin{Lemma}\label{lm676}
$d(\mathbf{P})=0$.
\end{Lemma}
\begin{proof}
Euclid in Proposition IX.36 in \emph{Elements} showed that  $\{2^k(2^{k+1}-1) : 2^{k+1}-1\in\mathcal{P}\}\subset 2\Z\cap\mathbf{P}$. In a posthumous 1849 paper, Euler proved the other inclusion, i.e., $2\Z\cap\mathbf{P}\subset \{2^k(2^{k+1}-1) : 2^{k+1}-1\in\mathcal{P}\}$, see~\cite{MR0245499}. Therefore,
$$
2\Z\cap\mathbf{P}=\{2^k(2^{k+1}-1) : 2^{k+1}-1\in\mathcal{P}\}.
$$
In particular, $d(2\Z\cap \mathbf{P})=0$. Moreover, $d((2\Z+1)\cap \mathbf{P})=0$ by \cite{MR0064805}\footnote{It is an open problem, whether $(2\Z+1)\cap \mathbf{P}=\emptyset$.} and we conclude.
\end{proof}
\begin{proof}[Proof of \cref{VFVF}]
By \cref{lm121} and \cref{OOK}, we have
\begin{equation}\label{kiki}
d(\{n\in\N: A+n \subset \mathbf{A} \text{ and }F+n \subset \mathbf{P}\cup\mathbf{D}\})>0.
\end{equation}
The assertion follows from \eqref{kiki} and \cref{lm676}.
 \end{proof}

\begin{proof}[Proof of \cref{VFVF1}]
Since $\{1,2,3,4,5\}\subset \mathbf{D}$, the assertion is an immediate consequence of \cref{VFVF}.
\end{proof}
\begin{Remark}
Notice that \cref{VFVF1} yields an indenpendent proof and strengthens the result from \cite{MR1540167} that there are infinitely many sequences of~5 consecutive deficient numbers.
\end{Remark}

\begin{Lemma}\label{infi}
$\sB_\mathbf{A}$ contains an infinite coprime subset.
\end{Lemma}
\begin{proof}
It follows from \cite{MR1574240} that $\bigcap_{1\leq k\leq K}(\mathcal{M}_{\mathscr{B}_\mathbf{A}}-k)\neq\emptyset$ for any $K\geq 1$, i.e., $(\dots,0,0,0,\dots)\in X_\eta$. 
To conclude, it suffices to use~\cref{TTB}.
\end{proof}
\begin{Remark}
Another way to prove the above lemma is this is to use the algorithm presented in~\cite{MR2134854}, outputting the smallest abundant number not divisible by the first $k$ primes.
\end{Remark}

\begin{proof}[Proof of \cref{VFVF2}]
The assertion is an immediate consequence of \cref{OOL}, \cref{lm121} and \cref{infi}.
\end{proof}

\begin{proof}[Proof of \cref{VFVF3}]
It follows by \cref{lm121}, \cref{infi} and \cref{TTD} that $X_\eta=\widetilde{X}_\eta$. In particular, by \cref{TTB}, $(S,X_\eta)$ is proximal. The intrinsic ergodicity of $(S,X_\eta)$ follows from the heredity of $X_\eta$ and from \cref{OOH}. Finally, the intrinsic heredity of $X_\eta$ and \cref{OOI} yields $h_{top}(S,X_\eta)=1-d(\mathbf{A})$. 
\end{proof}

\begin{Remark}
It remains open, whether we have $\widetilde{X}_\eta=X_{\mathscr{B}_\mathbf{A}}$. If the answer is positive, it would imply that given finite disjoint sets $A,B\subset \N$, one could always check in a finite number of steps whether for some $n\in\N$ we have $A+n\subset \mathbf{A}$ and $B+n\subset \mathbf{P}\cup \mathbf{D}$ (again, since $d(\mathbf{P})=0$, this is equivalent to the existence of $n\in\N$ such that $A+n\subset \mathbf{A}$ and $B+n\subset \mathbf{D}$). Indeed, it would be sufficient to check whether $|\text{supp }B \bmod b|<b$ for each $b\in\mathbf{A}\cap [1,\max B - \min B+1]$.
\end{Remark}

%%%%%%%%%%%%%%%%%%%%%%%%%%%%%%%%%%%%%%%%%%%%%%%%%%%%%%%%

%%%%%%%%%%%%%%%%%%%%%%%%%%%%%%%%%%%%%%%%%%%%%%%%%%%%%%%%

\footnotesize
\bibliography{cala.bib}

\bigskip
\footnotesize

\noindent
Aurelia Bartnicka\\
\textsc{Faculty of Mathematics and Computer Science, Nicolaus Copernicus University, Chopina 12/18, 87-100 Toru\'{n}, Poland}\par\nopagebreak
\noindent
\textit{E-mail address:} \texttt{aurbart@mat.umk.pl}

\medskip

\noindent
Stanis\l{}aw Kasjan\\
\textsc{Faculty of Mathematics and Computer Science, Nicolaus Copernicus University, Chopina 12/18, 87-100 Toru\'{n}, Poland}\par\nopagebreak
\noindent
\textit{E-mail address:} \texttt{skasjan@mat.umk.pl}

\medskip

\noindent
Joanna Ku\l aga-Przymus\\
\textsc{Institute of Mathematics, Polish Acadamy of Sciences, \'{S}niadeckich 8, 00-956 Warszawa, Poland}\\
\textsc{Faculty of Mathematics and Computer Science, Nicolaus Copernicus University, Chopina 12/18, 87-100 Toru\'{n}, Poland}\par\nopagebreak
\noindent
\textit{E-mail address:} \texttt{joanna.kulaga@gmail.com}

\medskip

\noindent
Mariusz Lema\'nczyk\\
\textsc{Faculty of Mathematics and Computer Science, Nicolaus Copernicus University, Chopina 12/18, 87-100 Toru\'{n}, Poland}\par\nopagebreak
\noindent
\textit{E-mail address:} \texttt{mlem@mat.umk.pl}

\end{document}